
\input amstex
\documentstyle{amsppt}
\magnification=\magstephalf
 \addto\tenpoint{\baselineskip 15pt
  \abovedisplayskip18pt plus4.5pt minus9pt
  \belowdisplayskip\abovedisplayskip
  \abovedisplayshortskip0pt plus4.5pt
  \belowdisplayshortskip10.5pt plus4.5pt minus6pt}\tenpoint
\pagewidth{6.5truein} \pageheight{8.9truein}
\subheadskip\bigskipamount
\belowheadskip\bigskipamount
\aboveheadskip=3\bigskipamount
\catcode`\@=11
\def\output@{\shipout\vbox{%
 \ifrunheads@ \makeheadline \pagebody
       \else \pagebody \fi \makefootline 
 }%
 \advancepageno \ifnum\outputpenalty>-\@MM\else\dosupereject\fi}
\outer\def\subhead#1\endsubhead{\par\penaltyandskip@{-100}\subheadskip
  \noindent{\subheadfont@\ignorespaces#1\unskip\endgraf}\removelastskip
  \nobreak\medskip\noindent}
\def\endremark{\par\revert@envir\endremark\vskip\postdemoskip}
\outer\def\enddocument{\par
  \add@missing\endRefs
  \add@missing\endroster \add@missing\endproclaim
  \add@missing\enddefinition
  \add@missing\enddemo \add@missing\endremark \add@missing\endexample
 \ifmonograph@ 
 \else
 \vfill
 \nobreak
 \thetranslator@
 \count@\z@ \loop\ifnum\count@<\addresscount@\advance\count@\@ne
 \csname address\number\count@\endcsname
 \csname email\number\count@\endcsname
 \repeat
\fi
 \supereject\end}
\catcode`\@=\active
\PSAMSFonts
\CenteredTagsOnSplits
\NoBlackBoxes
\nologo
\def\today{\ifcase\month\or
 January\or February\or March\or April\or May\or June\or
 July\or August\or September\or October\or November\or December\fi
 \space\number\day, \number\year}
\define\({\left(}
\define\){\right)}

\define\Aut{\operatorname{Aut}}
\define\CC{{\Bbb C}}

\define\End{\operatorname{End}}

\define\Hom{\operatorname{Hom}}
\define\Map{\operatorname{Map}}

\define\RR{{\Bbb R}}
\define\SS{\Bbb S}
\define\Spin{\operatorname{Spin}}

\define\Tr{\operatorname{Tr}}
\define\ZZ{{\Bbb Z}}
\define\[{\left[}
\define\]{\right]}

\define\chiup{\raise.5ex\hbox{$\chi$}}
\redefine\cir{S^1}

\define\dbar{{\bar\partial}}

\define\exertag #1#2{#2\ #1}

\define\inv{^{-1}}
\define\mstrut{^{\vphantom{1*\prime y}}}
\define\protag#1 #2{#2\ #1}
\define\rank{\operatorname{rank}}
\define\res#1{\negmedspace\bigm|\mstrut_{#1}}
\define\temsquare{\raise3.5pt\hbox{\boxed{ }}}

\define\theprotag#1 #2{#2~#1}

\define\xca#1{\removelastskip\medskip\noindent{\smc%
#1\unskip.}\enspace\ignorespaces }

\define\zmod#1{\ZZ/#1\ZZ}

\define\zt{\zmod2}

\define\rem#1{\marginalstar\begingroup\bf[{\eightpoint\smc{#1}}]\endgroup}
\def\strutdepth{\dp\strutbox} 
\def\marginalstar{\strut\vadjust{\kern-\strutdepth\specialstar}} 
\def\specialstar{\vtop to \strutdepth{ 
    \baselineskip\strutdepth 
    \vss\llap{$\bold{\Rightarrow}$ }\null}}

\define\conth#1{\medskip#1\smallskip}
\define\contsh#1#2{\par	\indent \S#1.\enspace#2\endgraf}
\subsubheadskip\bigskipamount

\NoRunningHeads 


\define\AR{\Cal{A}\mstrut _P}
\define\Ad{\operatorname{Ad}}
\define\BSU{\operatorname{BSU}}
\define\BTroth{\operatorname{B\Troth}}
\define\CJf#1{\CJ(H^*)_{\le #1}}
\define\CJ{\Cliff_{\scrJ}}
\define\Cliff{\operatorname{Cliff}}
\define\DAs{\DA^{\,2}}
\define\DA{\Dirac_A}
\define\Det{\operatorname{Det}}
\define\Dirac{D}

\define\Dzt{{\Dirac_0}^2}
\define\ESfin{\End\bigl(\Sfin(A) \bigr)}
\define\Efil#1{E(A)_{\le #1}}
\define\Er{\operatorname{E}}

\define\Adual{A^\vee}
\define\Fred{\operatorname{Fred}}
\define\Fr{\operatorname{F}}
\define\GR{G[P]}

\define\Gt{\tilde G}
\define\HfinsA{H_{\operatorname{fin}}^*(A)}
\define\IPER#1{\operatorname{IPER}(#1)}

\define\Ker{\operatorname{Ker}}
\define\LGRRt{(\widehat L\mstrut _PG)^\tau }
\define\LGRRz{(L\mstrut _PG)_{\text{rot}}}
\define\LGRR{\widehat L\mstrut _PG}
\define\LGRs{(\LGR)^\sigma }
\define\LGRts{(\LGR)^{\tms}}
\define\LGRt{(\LGR)^\tau }
\define\LGR{L\mstrut _PG}

\define\LSU{\operatorname{LSU}}
\define\LU{\operatorname{LU}}

\define\LgRRt{(\widehat L\mstrut _P\frak{g})^\tau }
\define\LgRR{\widehat L\mstrut _P\frak{g}}
\define\LgRfin{(\LgR\mstrut _\CC)\mstrut _{\operatorname{fin}}}
\define\LgRs{(\LgR)^\sigma }
\define\LgRts{(\LgR)^{\tms}}
\define\LgRt{(\LgR)^\tau }
\define\LgR{L\mstrut _P\frak{g}}
\define\Lie{\operatorname{Lie}}
\define\Lreg{(\Lambda ^\tau )^{\text{reg}}}
\define\Ltil{\Lambda^\tau}
\define\Ltt{U}
\define\Lt{\tilde\Lambda } 
\define\OGR{\Omega \mstrut _PG}
\define\OJ{\Or_{\scrJ}}
\define\ONA{\Omega _N(A)}
\define\Or{\operatorname{O}}
\define\PZ{\Pi \mstrut _Z} 
\define\Pin{\operatorname{Pin}}
\define\Pit{\tilde\Pi } 

\define\Rrot{\RR_{\text{rot}}}
\define\Rts{{R}^{\tau -\sigma }(\LGR)}
\define\Rt{{R}^\tau (\LGR)}
\define\SA{S_A}
\define\SO{\operatorname{SO}}
\define\SU{\operatorname{SU}}
\define\Sfin{\SS_{\operatorname{fin}}}
\define\Sym{\operatorname{Sym}}

\define\TT{\Bbb{T}}
\define\Troth{\widehat\TT_{\text{rot}}}
\define\Trot{\TT_{\text{rot}}}
\define\Tt{\tilde T}
\define\Ur{\operatorname{U}}
\define\Vfin{V_{\operatorname{fin}}}

\define\Wafft{\tilde W^e_{\text{aff}}}
\define\Waff{{W^e_{\text{aff}}}}
\define\Wfin{W_{\operatorname{fin}}}

\define\Wmin{W_{\emin}}
\define\ZAt{Z_A^\tau }
\define\ZA{Z_A}
\define\ad{\operatorname{ad}}
\define\aff#1#2{\Cal{A}_{#1}^{#2}}
\define\alcove{\frak{a}}
\define\am{\dot\pi _\mu }
\define\bul{^{\bullet}}
\define\calcove{\bar\frak{a}}
\define\cco{C_1^c}
\define\cd{\dot{\chi }}
\define\cent{\operatorname{center}}
\define\cj{\chi \mstrut _j}
\define\cpt{_{\text{cpt}}}
\define\dC#1{\check h(#1)}
\define\ejk{\eta _j^{(k)}}
\define\emin{e_{\text{min}}}
\define\epjk{{\eta '_j}^{(k)}}
\define\form{\ll\! \cdot,\cdot \!\gg}
\define\gC{\frak{g}_{\CC}}
\define\gn#1#2{\gamma ^{#1}(#2)}
\define\gpd{/\!/}

\define\hol{\operatorname{hol}}
\define\id{_1}
\define\ip#1#2{\ll\! #1,#2 \!\gg}
\define\ksth{k^*@!@!\th}
\define\ktb{\bar\kappa ^\tau }
\define\kts{\kappa ^{\tau -\sigma }}
\define\kt{\kappa ^\tau }
\define\lift#1{(#1)_A^\tau }
\define\lifts#1{(#1)_A^\sigma }
\define\liftso#1{(#1)_{A_0}^\sigma }
\define\longhookrightarrow{\lhook\joinrel\relbar\joinrel\rightarrow}

\define\nA{\frak{n}_A}
\define\opA{\Cal{F}_A}
\define\orbit{\Cal{O}}
\define\orb{\tilde\Cal{O}}
\define\pA{\frak{p}_A}
\define\pab{\pi \mstrut _{-\lambda ,-e}}
\define\pa{\pi \mstrut _{\lambda ,e}}
\define\plekb{\pi \mstrut _{-\lambda _k,-e_k}}
\define\plek{\pi \mstrut _{\lambda _k,e_k}}
\define\plepb{\pi \mstrut _{-\lambda ',-e'}}
\define\plep{\pi \mstrut _{\lambda ',e'}}

\define\rdA{\dot{\rho }_A}
\define\rtd{\dot{\rho}}
\define\rt{\hat{\rho }}

\define\sG{\Cal{G}}

\define\sH{\Cal{H}}
\define\sNh{\sigma (N)}
\define\sN{\Cal{N}}
\define\sP{\Cal{P}}
\define\scrJ{\Cal{J}}
\define\sh{\sigma }
\define\sign{\operatorname{sign}}
\define\sker{\SS_{\Ker}}

\define\smooth{\text{smooth}}
\define\spec{\operatorname{spec}}
\define\stab{Z_\mu }
\define\stl{\tilde{\sigma }}

\define\tKGG{K_G^{\tau + \dim G}(\GR)}
\define\tZ{\tau\mstrut _Z}

\define\th{\tau}

\define\tms{\tau -\sigma }
\define\torbit{\sigma (\orbit)}

\define\ttl{\tilde\tau }

\define\xA{\xi_A}
\define\xst{\xi ^*}
\define\zA{\frak{z}_A}
\define\zan#1{\zeta_{\lambda ,e}^{(#1)}}
\define\za{\zeta \mstrut _{\lambda ,e} }
\define\zet#1{\zeta \mstrut _{\lambda _#1,e_#1}}
\define\zlep{\zeta \mstrut _{\lambda ',e'}}
\define\zxct{(z^{-n}\,\bar{\xi}_A)^\tau }
\define\zxc{z^{-n}\,\bar{\xi}_A}
\define\zxt{(z^n\,\xi _A)^\tau }
\define\zx{z^n\,\xi _A}
\redefine\SS{S}
\input xy 
\xyoption{all}
\redefine\cir{S^1}

\refstyle{A}
\widestnumber\key{SSSSS}   
\document


	\topmatter
 \title\nofrills  Loop Groups and Twisted $K$-Theory II
\endtitle 
 \author Daniel S. Freed\\Michael J. Hopkins\\Constantin Teleman  \endauthor
 \thanks During the course of this work the first author was partially
supported by NSF grants DMS-0072675 and DMS-0305505, the second partially by
NSF grants DMS-9803428 and DMS-0306519, and the third partially by NSF grant
DMS-0072675. We also thank the KITP of Santa Barbara (NSF Grant PHY99-07949)
and the Aspen Center for Physics for hosting their summer programs, where
various sections of this paper were revised and completed.\endthanks
 \affil Department of Mathematics, University of Texas at Austin \\  
	Department of Mathematics, Harvard University \\  
	Department of Mathematics, University of Cambridge \endaffil 
 \address Department of Mathematics, University of Texas, 1 University
Station C1200, Austin, TX 78712-0257\endaddress 
 \email dafr\@math.utexas.edu \endemail
 \address Department of Mathematics, Harvard University, One Oxford Street,
Cambridge, MA 02138 \endaddress
 \email mjh\@math.harvard.edu\endemail 
 \address Department of Mathematics, University of California, 970 Evans
Hall, Berkeley, CA 94720-3840 \endaddress
 \email C.Teleman\@dpmms.cam.ac.uk\endemail
 \date November 16, 2012\enddate
 \abstract This is the second in a series of papers investigating the
relationship between the twisted equivariant $K$-theory of a compact Lie
group~$G$ and the ``Verlinde ring'' of its loop group.  We introduce the
Dirac family of Fredholm operators associated to a positive energy
representation of a loop group.  It determines a map from isomorphism classes
of representations to twisted $K$-theory, which we prove is an isomorphism if
$G$~is connected with torsion-free fundamental group.  We also introduce a
Dirac family for finite dimensional representations of compact Lie groups; it
is closely related to both the Kirillov correspondence and the equivariant
Thom isomorphism.  (In Part~III~\cite{FHT3} of this series we extend the
proof of our main theorem to arbitrary compact Lie groups~$G$ and provide
supplements in various directions.  In Part~I~\cite{FHT1} we develop twisted
equivariant $K$-theory and carry out some of the computations needed here.)
 \endabstract
	\endtopmatter
 \vskip -.3in

\document

 \comment
 lasteqno 0@  2
 \endcomment

Elliptic operators appear in different guises in the representation theory of
compact Lie groups.  The Borel-Weil construction~\cite{BW}, phrased in terms
of holomorphic functions, has at its heart the $\dbar$~operator on K\"ahler
manifolds.  The $\dbar$~operator and differential geometric vanishing
theorems figure prominently in the subsequent generalization by
Bott~\cite{B}.  An alternative approach using an algebraic Laplace operator
was given by Kostant~\cite{K1}.  The Atiyah-Bott proof~\cite{AB} of Weyl's
character formula uses a fixed point theorem for the $\dbar$~complex.  On a
spin K\"ahler manifold $\dbar$~can be expressed in terms of the Dirac
operator.  This involves a shift, in this context the $\rho
$-shift\footnote{$2\rho $~is the sum of the positive roots of a compact Lie
group~$G$; it may be identified with the first Chern class of the flag
manifold associated to~$G$.}  whose analog for loop groups appears in our
main theorem.  Dirac operators may be used instead of~$\dbar$ in these
applications to representation theory, and indeed they often appear
explicitly.

In this paper we introduce a new construction: the {\it Dirac family\/}
attached to a representation of a Lie group~$\sG$ which is either compact or
the loop group of a compact Lie group; in the latter case the representation
is restricted to have positive energy.  The Dirac family is a collection of
Fredholm operators parametrized by an affine space, equivariant for an affine
action of a central extension of~$\sG$ by the circle group~$\TT$.  For an
irreducible representation the support of the family is the coadjoint orbit
given by the Kirillov correspondence, and the entire construction is
reminiscent of the Fourier transform of the character~\cite{Ki,Rule~6}.  The
Dirac family represents a class in twisted equivariant $K$-theory, so we
obtain a map from representations to $K$-theory.  For compact Lie groups it
is a nonstandard realization of the twisted equivariant Thom homomorphism,
which was long-ago proved to be an isomorphism.\footnote{It is interesting to
note that the proof of this purely topological result uses Dirac operators
and index theory.}  Our main result, \theprotag{3.43} {Theorem}, is that this
map from representations to $K$-theory is an isomorphism when $\sG$~is a loop
group.  The existence of a construction along these lines was first suggested
by Graeme Segal; cf.~\cite{AS,\S5}.

The basic Dirac operator~\thetag{1.10} which appears in the Dirac family has
been termed the {\it cubic Dirac operator\/}.  Any Lie group~$G$ has a
distinguished line segment in its affine space of linear connections: the
endpoints are the flat connections~$\nabla _L, \nabla _R$ that are the
infinitesimal versions of the global parallelisms via left and right
translation.  If $G$~is compact then the Levi-Civita connection of a
bi-invariant metric is~$\frac 12\nabla _L + \frac 12\nabla _R$.  The
cubic\footnote{`Cubic' refers to the tangent vector to this line of
connections, which is the invariant 3-form on~$G$.  The Dirac operator
attached to {\it any\/} connection on this line has a cubic term in its local
formula.  The particular connection and Dirac operator of interest is
distinguished by the coefficient in front of that term.  The apparent
asymmetry between left and right is explained by our use of {\it left\/}
translation to trivialize the spin bundle on~$G$.} Dirac operator is
associated to the connection~$\frac 23\nabla _L + \frac 13\nabla _R$.  This
particular Dirac operator was introduced by Slebarski~\cite{Sl} and used by
Kostant in~\cite{K2}.  It also enjoys nice analytic properties, as
in~\cite{G}.  Alekseev and Meinrenken~\cite{AM} interpret this Dirac operator
as a differential on a quantized version of the Weil algebra.  It also
appears in the generalized geometry of Hitchin~\cite{H,\S10}.   The infinite
dimensional version appears in Taubes' work~\cite{T} on loop spaces.  It was
subsequently used by Landweber~\cite{L} in a generalization of Kostant's
paper to loop groups.  We later learned that many relevant formul\ae\ were
independently known in the physics literature, and indeed our Fredholm
family, which is a gauge-coupled Dirac-Ramond operator, was first flagged in
relation to twisted (non-equivariant) $K$-theory by Mickelsson~\cite{M}.

The finite dimensional case is developed in detail in~\S{1}.  Let $G$~be any
compact Lie group, $V$~a finite dimensional unitary representation, and
$\SS$~the spin space of the adjoint representation.  The Dirac family
consists of endomorphisms of~$V\otimes \SS$ parametrized by the
dual~$\frak{g}^*$ of the Lie algebra of~$G$.  It is equivariant for the
co-adjoint action of~$G$, possibly centrally extended to act on the spin
space, and the endomorphisms are invertible outside a compact set.  In fact,
if $V$~is irreducible we show in \theprotag{1.19} {Proposition} that the the
endomorphisms fail to be invertible only on a single coadjoint
orbit~$\orbit\subset \frak{g}^*$.  The interpretation in terms of the
equivariant Thom isomorphism is \theprotag{1.28} {Theorem}.
 
Our main application is to positive energy representations of loop groups.
In~\S{2} we review and develop the concept of energy in the theory of loop
groups.  From the beginning we work with twisted loop groups, that is, the
group of gauge transformations~$\LGR$ of a principal $G$-bundle $P\to\cir$.
We introduce the slightly larger group of gauge transformations which cover a
rigid rotation of~$\cir$.  Then a central extension~$\LGRt$ of~$\LGR$ by the
circle group~$\TT$ is said to be {\it admissible\/} if it extends over this
larger group---so is rotation-invariant---and if there is an invariant
bilinear form which pairs the Lie algebras of the center~$\TT$ and the
rotations~$\Trot$; see \theprotag{2.10} {Definition} for the precise
conditions.\footnote{We prove in an appendix that if the Lie algebra of
$G$~is semisimple, then any central extension is admissible.  On the other
hand, if $G$~is a torus of dimension at least two, then there exist
non-admissible central extensions.}  This bilinear form, central in Kac's
algebraic theory~\cite{Ka,\S6}, plays a crucial role here as it sets up a
correspondence between connections on~$P\to\cir$ and linear splittings of the
central extension of loop algebras (\theprotag{2.18} {Lemma}).  The space of
connections enters also in \theprotag{2.8} {Definition}, which associates to
the central extension~$\LGRt$ a {\it twisting\/} of the equivariant
$K$-theory of $G$~acting on itself by conjugation.  One novelty here is the
inclusion of a {\it grading\/}, a homomorphism $\LGR\to\zt$, in~$\tau $; it
affects a component of the associated twisting.  Infinitesimal rotations
measure energy, but the precise definition depends on a choice of connection.
Thus we obtain a family of energies parametrized by the space~$\AR$ of
connections.  Following~\cite{PS} in~\S{2.5} we introduce positive energy
representations, extending the standard definition to allow for gradings.
For a fixed admissible graded central extension~$\LGRt$ there is a finite set
of isomorphism classes of irreducibles.  They generate a free abelian group
we denote~$\Rt$.
 
The Dirac family construction is taken up again in~\S{3}, now in the infinite
dimensional setting of loop groups.  The adjoint spin representation~$\SS$ of
the loop group determines a distinguished central extension~$\LGRs$.  We form
a family of Dirac operators by tensoring a positive energy representation
of~$\LGRts$ with spinors.  This gives a family of Fredholm operators
parametrized by~$\AR$, equivariant for the central extension~$\LGRt$.  Here
we encounter the {\it adjoint shift\/}\footnote{For $G=\SU(2)$ it is the
famous $k\to k+2$ shift which occurs in low dimensional physics: conformal
field theory, Chern-Simons theory, etc.  It is a loop group analog of the
$\rho $-shift for compact Lie groups.  For an interpretation in terms of the
intrinsic geometry of the loop group see~\cite{F}.} by~$\sigma $.  This
Fredholm family represents an element of twisted $K$-theory in the model
developed in Part~I, and so the Dirac construction induces a homomorphism
  $$ \Phi \:\Rts\longrightarrow \tKGG. \tag{0.1} $$
Here $\GR$~is the union of components of~$G$ consisting of all holonomies of
connections on~$P\to\cir$, and $G$~acts on it by conjugation.  Our main
result, \theprotag{3.43} {Theorem}, asserts that $\Phi $~is an isomorphism.

The proof, presented in~\S{4} for the case when $G$~is connected with
torsion-free fundamental group, is computational: we compute both sides
of~\thetag{0.1} and check that $\Phi $~induces an isomorphism.  We deduce the
result for this class of groups from the special cases of tori and simply
connected groups.  The positive energy representations of the loop groups in
these cases---the left hand side of~\thetag{0.1}---are enumerated
in~\cite{PS}.  The twisted equivariant $K$-theory---the right hand side
of~\thetag{0.1}---is computed in~\cite{FHT1,\S4}.  The primary work here is
the analysis of the kernel of the Dirac family.  This is parallel to the
finite dimensional case in~\S{1}, and in fact reduces to it.  One of the main
points is a Weitzenb\"ock-type formula~\thetag{3.36} which relates the square
of Dirac to energy.  Specific examples are written out in~\S{4.4}.
 
In Part~III~\cite{FHT3} we complete the proof of \theprotag{3.43} {Theorem}
for any compact Lie group.  We also generalize to a wider class of central
extensions of loop groups which do not rely on energy.  The positive energy
condition on representations is replaced by an integrability condition at the
Lie algebra level~\cite{Ka}, and the whole treatment there relies much more
on Lie algebraic methods.  We also give a variation which incorporates energy
more directly into~\thetag{0.1}, and several other complements.  In many
cases the twisted equivariant $K$-theory~$K_G^{\tau+\dim G}(G)$ is a ring, in
fact a Frobenius ring.  In those cases the twisting~$\tau $ is derived from a
more primitive datum---a ``consistent orientation''---as we explain
in~\cite{FHT4}. 
 
The twistings of equivariant $K$-theory we encounter in this Part~II have a
special form relative to the general theory of Part~I.  Suppose a compact Lie
group~$G$ acts on a space~$X$.  Let $\sP$~be a space with the action of a
topological group~$\sG$, suppose the normal subgroup~$\sN\subset \sG$ acts
freely, and there are isomorphisms $\sP/\sN\cong X$ and $\sG/\sN\cong G$
compatible with the group actions.  In the language of Part~I, the quotient
maps define a {\it local equivalence of groupoids\/} $\sP \gpd\sG\to X\gpd
G$.  In this paper we use twistings of~$K_G(X)$ of the form
  $$ \th = (\sP @>\sN>> X, \sG^\tau , \epsilon ^\tau ). \tag{0.2} $$
These consist of $\sP\to X$ as above, a central extension~$\sG^\tau $
of~$\sG$ by the circle group~$\TT$, and a homomorphism $\epsilon ^\tau \:\sG
\to\zt$, termed a {\it grading\/}.  Together they define a graded central
extension of~$\sP\gpd\sG$; compare~\cite{FHT1,\S2.2}.

The authors which to express their gratitude to Graeme Segal for many helpful
conversations as well as admiration for his development of the geometric
approach to loop groups.

\bigskip
{\eightpoint\parindent0pt
 {\bf Contents}
 \conth{Section 1: The finite dimensional case}
 \contsh{1.1}{The spin representation (finite dimensional case)}
 \contsh{1.2}{The cubic Dirac operator}
 \contsh{1.3}{The kernel of~$\Dirac_\mu $ and the Kirillov correspondence}
 \contsh{1.4}{$K$-theory interpretation}
 \contsh{1.5}{Variation: projective representations}
 \contsh{1.6}{Examples}
 \conth{Section 2: Loop groups and energy}
 \contsh{2.1}{Basic definitions}
 \contsh{2.2}{Admissible central extensions}
 \contsh{2.3}{Special cases}
 \contsh{2.4}{Finite energy loops}
 \contsh{2.5}{Positive energy representations}
 \conth{Section 3: Dirac families and loop groups}
 \contsh{3.1}{The spin representation (infinite dimensional case)}
 \contsh{3.2}{The canonical 3-form on~$\LGR$}
 \contsh{3.3}{A family of cubic Dirac operators on loop groups; main theorem}
\newpage
 \conth{Section 4: Proofs} 
 \contsh{4.1}{The proof for tori}
 \contsh{4.2}{The proof for simply connected groups}
 \contsh{4.3}{The proof for connected~$G$ with $\pi _1$~torsion-free}
 \contsh{4.4}{Examples}
 \conth{Appendix: Central extensions in the semisimple case}
}\bigskip

 \head
 \S{1} The finite dimensional case
 \endhead
 \comment
 lasteqno 1@ 39
 \endcomment

 \subhead \S{1.1}.  The spin representation (finite dimensional case)
 \endsubhead

The basic reference is~\cite{ABS}.  Let $H$~be a finite dimensional real
vector space with an inner product~$\langle \cdot ,\cdot \rangle$.  Then
there are central extensions~$\Pin^{\pm}(H),\,\Pin^c(H)$ of the orthogonal
group~$\Or(H)$ which fit into the diagram
  $$ \CD
      1 @>>> \zt @>>> \Pin^{\pm}(H) @>>> \Or(H) @>>> 1\\
      @. @VVV @VVV @| \\
      1 @>>> \TT @>>> \Pin^c(H) @>>> \Or(H) @>>> 1 \endCD \tag{1.1} $$
There is a split extension of Lie algebras
  $$ 0 \longrightarrow i\RR \longrightarrow \frak{p}\frak{i}\frak{n}^c(H)
     \leftrightarrows \frak{o}(H) \longrightarrow 0, \tag{1.2} $$
the splitting induced from the first line of~\thetag{1.1}.  (Throughout
$i=\sqrt{-1}$.)  The Clifford algebra~$\Cliff^{\pm}(H^*)$ of~$H^*$ is the
universal unital associative algebra with a linear map $\gamma \:H^*\to
\Cliff^{\pm}(H^*)$ such that
  $$ \gamma (\mu )\gamma (\mu' ) + \gamma (\mu' )\gamma (\mu ) =
     \pm2\langle \mu ,\mu'  \rangle,\qquad \mu ,\mu' \in H^*. \tag{1.3} $$
Then $\Cliff^{\pm}(H^*)$~is $\zt$-graded, and the left hand side
of~\thetag{1.3} is the graded commutator ~$\bigl[\gamma (\mu ),\gamma (\mu
')\bigr]$.  Also, $\Pin^{\pm}(H)\subset \Cliff^{\pm}(H^*)$ and
$\Pin^c(H)\subset \Cliff^{\pm}(H^*)\otimes \CC$, and these groups inherit a
grading~$\epsilon $ from that of the Clifford algebra.  (A {\it grading\/} of
a Lie group~$G$ is a continuous homomorphism $\epsilon \:G\to\zt$.)  Let
$e^1,\dots ,e^n$ be an orthonormal basis of~$H^*$.  Then $\gamma (e^1)\cdots
\gamma (e^n)\in \Cliff(H^*)$ is a {\it volume form\/}, and is determined up
to sign.  If $H$~is oriented we use oriented orthonormal bases and so fix the
sign.  In this paper we consider {\it complex\/} modules
for~$\Cliff^{\pm}(H^*)$, and we denote $\Cliff^-(H^*)\otimes \CC \cong
\Cliff^+(H^*)\otimes \CC$ by~$\Cliff^c(H^*)$.  Also, if~$H=\RR^n$ with the
standard metric, then we use the notation~$C_n^{\pm}$ for the real Clifford
algebras, $C_n^c$~for the complex Clifford algebra, and $\Pin_n^{\pm}$~for
the pin groups.

The isomorphism classes of $\zt$-graded complex $\Cliff^c(H^*)$-modules form
a semigroup whose group completion is isomorphic to~$\ZZ$ if $\dim H$~is odd
and $\ZZ\oplus \ZZ$ if $\dim H$~is even.  If $\dim H$~is even, then the two
positive generators are distinguished by the action of a volume form on the
even component of an irreducible $\zt$-graded module.  If $\dim H$~is odd,
then $\Cliff^c(H^*)$~has a nontrivial center (as an ungraded algebra) which
is isomorphic to $\cco$; an element in the center has the form~$a+b\omega $
where $a,b\in \CC$ and $\omega $~is a volume form.  In the odd case we always
consider $\Cliff^c(H^*)$-modules with a commuting $\cco$-action which is not
necessarily that of the center, even though we do not explicitly mention that
action.  The group completion of isomorphism classes is isomorphic
to~$\ZZ\oplus \ZZ$.  The positive generators are distinguished by the action
of the product of a volume form in~$\Cliff^c(H^*)$ and a volume form in~$\cco$
on the even component.  For any~$H$ the equivalence relation of $K$-theory
($K^0$~if $\dim H$~is even, $K^1$~if $\dim H$~is odd) identifies the positive
generator of one summand~$\ZZ$ with the negative generator of the other.

Let~$\SS$ be an irreducible complex $\zt$-graded $\Cliff^c(H^*)$-module.
Then $\SS$~restricts to a graded representation of~$\Pin^c(H)$, so to a
projective representation of~$\Or(H)$.  At the Lie algebra level we use the
splitting~\thetag{1.2} to obtain a representation of ~$\frak{o}(H)$.  There
exist {\it compatible\/} metrics on~$\SS$---metrics such that
$\Pin^c(H)$~acts by unitary transformations and such that $\gamma (\mu )$~is
a skew-Hermitian transformation for each~$\mu \in H^*$.

 \subhead \S{1.2}.  The cubic Dirac operator
 \endsubhead

Let $G$~be a compact Lie group and $\langle \cdot ,\cdot \rangle$ a
$G$-invariant inner product on the Lie algebra~$\frak{g}$.  We apply~\S{1.1}
to~$V=\frak{g}$ using the adjoint representation $G\to \Or(\frak{g})$.  Thus by
pullback from~\thetag{1.1} with~$H=\frak{g}$ we obtain a graded central
extension
  $$ 1 \longrightarrow \TT \longrightarrow G^\sigma \longrightarrow
     G\longrightarrow 1 \tag{1.4} $$
and a canonically split extension of Lie algebras
  $$ 0 \longrightarrow i\RR \longrightarrow \frak{g}^\sigma \leftrightarrows
     \frak{g} \longrightarrow 0. \tag{1.5} $$
If $G$~is connected and simply connected, then \thetag{1.4}~is canonically
split as well; in general it is induced from an extension by~$\zt$.
Choose an {\it irreducible\/} $\zt$-graded complex
$\Cliff^c(\frak{g}^*)$-module~$\SS$ with a compatible metric.  Then
$\SS$~carries a unitary representation of~$G^\sigma $.
 
Fix a basis~$\{e_a\}$ of~$\frak{g}$ and let $\{e^a\}$~be the dual basis
of~$\frak{g}^*$.  Define tensors~$g$ and~$f$ by
  $$ \alignedat{2}
      \langle e_a,e_b \rangle &= g_{ab}\qquad &\langle e^a,e^b \rangle&=
     g^{ab} \\
      [e_a,e_b] &= f^c_{ab}e_c \qquad &\langle [e_a,e_b],e_c \rangle &=
     f_{abc}.\endaligned $$
The invariance of the inner product implies that $f_{abc}$~is skew in the
indices.  Define the invariant 3-form
  $$ \Omega = \frac 16 f_{abc}e^a\wedge e^b\wedge e^c \tag{1.6} $$
on~$\frak{g}$.  Let $\gamma ^a$~denote\footnote{We view~$\Cliff^c$ as the
complexification of~$\Cliff^-$ as reflected by the sign in the second
formula of~\thetag{1.7}.} Clifford multiplication by~$e^a$, and $\sigma
_a$~the action of~$e_a$ on spinors (via the splitting~\thetag{1.5}); both are
skew-Hermitian transformations of~$\SS$, the former odd and the latter even.
Then we have
  $$ \aligned
      \sigma _a &= \frac 14f_{abc}\gamma ^b\gamma ^c \\
      [\gamma ^a,\gamma ^b] &= -2g^{ab} \\
      [\sigma _a,\sigma _b] &= f^c_{ab}\sigma _c \\
      [\sigma _a,\gamma ^b] &= -f^b_{ac}\gamma ^c.\endaligned \tag{1.7} $$
Now any element of~$\frak{g}$ may  be identified with a left-invariant vector
field on~$G$, and so acts  on smooth functions by differentiation.  Let~$R_a$
denote the differentiation corresponding to~$e_a$.  It satisfies
  $$ [R_a,R_b] = f^c_{ab}R_c. \tag{1.8} $$
We use left translation to identify the tensor product $C^{\infty}(G)\otimes
\SS$ with spinor fields on~$G$, and let $R_a,\sigma _a,\gamma ^a$~operate on
the tensor product.  Then
  $$ [R_a,\sigma _b] = [R_a,\gamma ^b] =0. \tag{1.9} $$
 
Introduce the Dirac operator 
  $$ \Dirac_0 = i\gamma ^aR_a + \frac i3 \gamma ^a\sigma _a 
      \;=\; i\gamma ^aR_a + \frac i{12} f_{abc}\gamma ^a\gamma ^b\gamma ^c  
      \;=\; i\gamma ^aR_a + \frac i{2} \gamma (\Omega ). 
      \tag{1.10} $$
The last term is, up to a factor, Clifford multiplication by the
3-form~$\Omega $, accomplished using the canonical vector space isomorphism
between the exterior algebra and Clifford algebra.  The Dirac
operator~$\Dirac_0$~ is an odd formally skew-adjoint operator on smooth
spinor fields.  Observe that $\Dirac_0$~is $G^\sigma $-invariant.  It is not
the Riemannian Dirac operator, nor is it either Dirac operator associated to
the natural parallel transports on~$G$ using left or right
translation.\footnote{See the introduction for further discussion.  The
choice of coefficient in the second term of~\thetag{1.10} is made so that the
square of the Dirac operator has no first order term~$g^{ab}R_a\sigma _b$:
  $$ \Dzt = \frac 12[\Dirac_0,\Dirac_0 ] = g^{ab}(R_aR_b + \frac 13\sigma
     _a\sigma _b) .  \tag{1.11} $$
The first term in~\thetag{1.11} is the nonpositive Laplace operator on spinor
fields, formed using right translation as parallel transport.  Of course, the
entire operator $\Dzt$~is nonpositive as it is the square of a skew-adjoint
operator.  Formula~\thetag{1.11} follows from a straightforward computation
using~\thetag{1.7}, \thetag{1.8}, and~\thetag{1.9}, and it has an easy
extension to the computation of the square of~\thetag{1.12}.  We could use
these formulas to prove \theprotag{1.19} {Proposition}, but instead give an
argument which more directly extends to the infinite dimensional case; see
the proof of \theprotag{3.33} {Proposition}.}  Define the {\it Dirac
family\/} of skew-adjoint Dirac operators parametrized by~$\mu \in
\frak{g}^*$:
  $$ \Dirac_\mu = \gamma (\mu ) + \Dirac _0 = \mu _a\gamma ^a +
     \Dirac_0,\qquad \mu =\mu _ae^a\in \frak{g}^*. \tag{1.12} $$
The map $\mu \mapsto \Dirac_\mu $ is $G^\sigma $-equivariant.\footnote{There
is a real version of the Dirac family, which we have not completely
investigated.  Roughly speaking we replace the spin space~$S$ by the vector
space~$\Cliff^-(\frak{g}^*)$: the group~$G^\sigma $ acts on the left and
there is a commuting action of~$C^-_{\dim G}$ on the right.  Then in the
Dirac family $D(V)\:\frak{g}^*\to\End\bigl(V\otimes
\Cliff^-(\frak{g}^*)\bigr)$ we rotate both copies of~$\frak{g}^*$
to~$\sqrt{-1}\frak{g}^*$ thus obtaining a Real family of operators under the
involution of complex conjugation.  It represents an element of twisted
$KR$-theory.}

Recall that the Peter-Weyl theorem gives an orthogonal direct sum
decomposition of the $L^2$~spinor fields, where we trivialize the spin bundle
using left translation:
  $$ L^2(G)\otimes \SS \;\;\;\cong \bigoplus\limits_{\text{$V$ irreducible} }
     V^* \otimes V \otimes \SS.  $$
Here $V$~runs over (representatives of) equivalence classes of irreducible
representations of~$G$.  In each summand the group~$G$ operates via the
projective spin representation on~$\SS$, via right translation on~$V$, and
via left translation on~$V^*$.  Each summand is a finite dimensional space of
smooth spinor fields, and $\Dirac_\mu $~preserves the decomposition.  Since
$\Dirac_\mu $~only involves (infinitesimal) right translation, it operates
trivially on~$V^*$.  We study each summand separately and drop the
factor~$V^*$ on which $\Dirac_\mu $~acts trivially.  In other words, for each
finite dimensional irreducible unitary representation~$V $ of~$G$ we consider
the Dirac family of {\it finite dimensional\/} odd skew-adjoint operators
  $$ \aligned
      \Dirac(V):\frak{g}^* &\longrightarrow \End(V \otimes \SS) \\
      \mu &\longmapsto \Dirac_\mu \endaligned \tag{1.13} $$
given by formulas~\thetag{1.12} and~\thetag{1.10}, where now $R$~is the
infinitesimal action of~$\frak{g}$ on the representation~$V$.  Note that
$V$~need not be irreducible, though to analyze the operator it is convenient
to assume so.

 \subhead \S{1.3}.  The kernel of~$\Dirac_\mu $ and the Kirillov correspondence
 \endsubhead

Fix a maximal torus~$T$ in the identity component~$G\id\subset G$, and let
$\frak{t}\subset \frak{g}$ be its Lie algebra.  The lattice of characters
of~$T$, or {\it weights\/} of~$G$, is~$\Lambda =\Hom(T,\TT)$.  Write a
character of~$T$ as a homomorphism $e^{i\lambda }\:T\to\TT$ for some~$\lambda
\in \frak{t}^*$, and so embed~$\Lambda \subset \frak{t}^*$.  Under the
adjoint action of~$T$ the complexification of~$\frak{g}$ splits into a sum of
the complexification of~$\frak{t}$ and one-dimensional {\it root\/}
spaces~$\frak{g}_{\pm\alpha }$:
  $$ \frak{g}_{\CC} \cong \frak{t}_{\CC} \oplus \bigoplus\limits_{\alpha \in
     \Delta ^+} (\frak{g}_\alpha \oplus \frak{g}_{-\alpha }). \tag{1.14} $$
The roots form a finite set~$\Delta \subset \Lambda$ and come in
pairs~$\pm\alpha $ whose kernels form the infinitesimal diagram, a set of
hyperplanes in~$\frak{t}$.  The components of the complement of the
infinitesimal diagram are simply permuted by the {\it Weyl
group\/}~$W=N(T)/T$.  A {\it Weyl chamber\/} is a choice of component, and it
determines the set~$\Delta ^+$ of positive roots, roots which have positive
values on the Weyl chamber. The bracket~$i[\frak{g}_\alpha ,\frak{g}_{-\alpha
}]$ is a line in~$\frak{t}$, and the {\it coroot\/}~$H_\alpha $ is the
element in that line so that $\ad(H_\alpha )$~acts as multiplication by~2
on~$\frak{g}_\alpha $.  The coroot~$H_\alpha $ is positive if $\alpha $~is a
positive root.  The kernels of the coroots from a set of hyperplanes
in~$\frak{t}^*$, and the {\it dual Weyl chamber\/} is the component of the
complement which takes positive values on the positive coroots.  A weight is
{\it dominant\/} if it lies in the closure of the dual Weyl chamber.  A
weight~ $\mu \in \Lambda$ is {\it regular\/} if the Weyl group~$W$ acts
on~$\mu $ with trivial stabilizer.

Define
  $$ \rho =\frac 12\sum\limits_{\alpha \in \Delta ^+}\alpha .  $$
It is an element of~$\frak{t}^*$ which is not necessarily in the weight
lattice~$\Lambda $.  We follow~\cite{PS,\S2.7} and label an irreducible
representation by its {\it lowest\/} weight, which is {\it antidominant\/}.

        \proclaim{\protag{1.15} {Lemma}}
 As an ungraded representation $\SS$~is a sum of irreducible representations
of~$G\id^\sigma $ of lowest weight~$-\rho $.
        \endproclaim

        \demo{Proof}  
 Let $\frak{p}=\bigoplus_{\alpha \in \Delta ^+}\frak{g}_\alpha $ be the sum
of the positive root spaces, so that \thetag{1.14}~reads $\frak{g}_{\CC}\cong
\frak{t}_{\CC}\oplus \frak{p}\oplus \bar{\frak{p}}$.  Then as a $\zt$-graded
representation of the central extension~$T^\sigma $ of~$T$ we have the
$\zt$-graded tensor product decomposition
  $$ \split
      \SS =\SS(\frak{g}^*)&\cong \SS(\frak{t}^*) \otimes {\tsize\bigwedge}
     \bar{\frak{p}}^* \otimes (\Det \frak{p}^*)^{1/2} \\
      &\cong \SS(\frak{t}^*) \otimes {\tsize\bigwedge} {\frak{p}}\phantom{*}
     \otimes (\Det \bar\frak{p}\phantom{*})^{1/2},\endsplit \tag{1.16} $$
where $\SS(\frak{t}^*)$~is a fixed irreducible complex $\zt$-graded spin
module of~$\frak{t}^*$.  We claim that the weights of~$T^\sigma $ which occur
in $\SS'={\tsize\bigwedge} \frak{p}\otimes (\Det \bar\frak{p})^{1/2}$ are
those of the irreducible representation of~$G\id^\sigma $ of lowest
weight~$-\rho $.  The lemma follows directly from the claim as $T^\sigma
$~acts trivially on~$\SS(\frak{t}^*)$.  To verify the claim, note that the
character of~${\tsize\bigwedge} \frak{p}$ as an ungraded representation
of~$T$ is $\prod \limits_{\alpha \in \Delta ^+}(1+e^{i\alpha })$, so the
character\footnote{Characters of representations of~$T^\sigma $ on which the
central~$\TT$ acts as scalar multiplication form a torsor for the character
group of~$T$ and may be identified with a subset of~$\frak{t}^*$ via the
splitting~\thetag{1.5}.  The character~\thetag{1.17} may be expressed as a
sum of exponentials of elements of~$\frak{t}^*$ which are necessarily in this
subset.} of~$\SS'$ is
  $$ \prod\limits_{\alpha \in \Delta ^+} (e^{i\alpha /2} + e^{-i\alpha
     /2}) = \prod\limits_{\alpha \in \Delta ^+} \frac{e^{i\alpha } -
     e^{-i\alpha}}{e^{i\alpha /2} - e^{-i\alpha /2}}. \tag{1.17} $$
That this is the character of the irreducible representation of lowest
weight~$-\rho $ follows from the Weyl character formula; or, we observe that
\thetag{1.17}~is Weyl-invariant, has lowest weight~$-\rho $, and that $\dim
\SS'=2^{\#\Delta ^+}$ has the correct dimension by the Weyl dimension
formula.
        \enddemo

We now analyze the operator~\thetag{1.13}.  Fix~$\mu \in \frak{g}^*$.  Choose
a maximal torus~$T_\mu $ in the identity component of the stabilizer subgroup
~$Z_\mu \subset G$.  If $\mu $~is regular, the identity component of~$Z_\mu $
{\it is\/}~$T_\mu $.  Then $\mu $~lies in $\frak{t_\mu }^*\subset
\frak{g}^*$, i.e., $\mu $~annihilates the root spaces.  Furthermore, choose a
Weyl chamber so that $\mu $~is {\it antidominant\/}.\footnote{If $\mu $~is
regular, it is contained in a unique (negative) dual Weyl chamber and so
determines a notion of anti-dominance.}  Let $\rho =\rho (\mu )\in
\frak{t}_\mu ^*$~be half the sum of positive roots, defined relative
to~$T_\mu$.  Under $G\id\subset G$ there is a decomposition
  $$ V\cong \bigoplus V_{-\lambda } \tag{1.18} $$
with $V_{-\lambda }$~the $(-\lambda )$-isotypical component of~$V$, where
$\lambda $ runs over the dominant weights of~$G\id$.

        \proclaim{\protag{1.19} {Proposition}}
 Suppose $V$~is irreducible.  The operator~$\Dirac_\mu (V)$ is nonsingular
unless $\mu $~is regular and $-\lambda =\mu +\rho$~is a lowest weight in the
decomposition~\thetag{1.18} of~$V$ under~$G\id$.  In that case
  $$ \Ker \Dirac_{\mu } = K_{-\lambda } \otimes \SS_{-\rho }\subset
     V\otimes \SS, \tag{1.20} $$
where $K_{-\lambda }$ is the $(-\lambda )$-weight space of~$V$ in the
decomposition under~$T_\mu \subset G\id$ and $\SS_{-\rho }$~is the
$\zt$-graded $(-\rho )$-weight space of~$\SS$ in the decomposition
under~$T_\mu ^\sigma \subset G\id^\sigma $.  The stabilizer ~$\stab\subset G$
of~$\mu $, a group with identity component~$T_\mu $, acts irreducibly
on~$K_{-\lambda }$.  The Clifford algebra~$\Cliff^c(\frak{t}_\mu ^*)$ acts
irreducibly on~$\SS_{-\rho }$.
        \endproclaim

\flushpar
 If $G=G\id$~is connected, then $\stab=T_\mu $ and $\dim K_{-\lambda }=1$.
The latter is a standard fact in the representation theory of connected
compact Lie groups.  The corresponding fact in the disconnected case---that
$\stab$~acts irreducibly on~$K_{-\lambda }$---is less standard; it may be
found in~\cite{DK,\S4.13} for example.

        \demo{Proof}
 Introduce a shifted version $\am\:\frak{g}\to \End(V\otimes \SS)$ of the
infinitesimal Lie algebra action by the formula
  $$ \am(\xi ) = \xi ^a (R_a + \sigma _a - i\mu _a),\qquad \xi \in
     \frak{g}.  $$
Notice that $\am(\xi )$~is skew-adjoint.  The inner product~$\langle \cdot
,\cdot \rangle$ on~$\frak{g}$ gives isomorphisms $\xi \mapsto \xi ^*$ and
$\mu \mapsto \mu _*$ between~$\frak{g}$ and its dual.  Introduce also the
self-adjoint operator
  $$ E_\mu = i\am(\mu _*) - \frac{|\mu |^2}{2}. \tag{1.21} $$
Then since $\mu \in \frak{t}^*$~ it follows that $E_\mu $~is constant on
weight spaces of~$V\otimes \SS$---its value on the space with weight~$\omega $
is 
  $$ \langle \frac\mu 2- \omega ,\mu \rangle. \tag{1.22} $$
Since $\mu $~is antidominant $E_\mu $~has its minimum on the lowest weight
space.  Straightforward computations with \thetag{1.7}--\thetag{1.10} give
for~$\xi \in \frak{g}$ the second of the commutation relations
  $$ \aligned
      [\Dirac_\mu ,\am(\xi  )] &= [-iE_\mu ,\gamma(\xi ^* )], \\
      [\Dirac_\mu ,\gamma(\xi  ^*)] &= -2i\am(\xi );
     \endaligned   $$
the first  is the infinitesimal form  of the statement that  $\Dirac_\mu $ is
equivariant for the action of~$G^\sigma $.  Iterating we find
  $$ \aligned
      [\Dirac_\mu^2 ,\am(\xi )] &= [-2E_\mu ,\am(\xi )], \\
      [\Dirac_\mu^2 ,\gamma(\xi ^*)] &= [-2E_\mu ,\gamma(\xi ^*)].
     \endaligned  $$
Now $V\otimes \SS$ is generated from the lowest weight space by applying
operators~$\am(\xi )$ and~$\gamma (\xi ^*)$.  Hence $D_\mu ^2 + 2E_\mu $~is
constant on~$V\otimes \SS$, and it is easy to check from~\thetag{1.12}
and~\thetag{1.21} that this constant is independent of~$\mu $.  Therefore,
the self-adjoint nonpositive operator~$\Dirac_\mu ^2$ has its maximum when
$E_\mu $~achieves its minimum.  By~\thetag{1.22} and the remark which follows
this happens at~$\mu =-(\lambda +\rho )$ on the minimum weight
space~$K_{-\lambda }\otimes \SS_{-\rho }$.  It remains to check that
$\Dirac_\mu ^2$~does vanish on this space.

We use a complex basis compatible with~\thetag{1.14}.  Fix nonzero
elements~$e_{\pm \alpha }\in \frak{g}_{\pm\alpha }$, and choose a basis~$e_t$
of~$\frak{t}$.  Extend the inner products on~$\frak{g},\frak{g}^*$ to
bilinear forms on the complexifications.  Write
  $$ \Dirac_\mu = \gamma ^a (iR_a + \mu _a) + \frac i3 \gamma ^a\sigma
     _a.  $$
Then for~$\alpha \in \Delta ^+$ we have $R_{-\alpha } = \gamma ^{\alpha }=0$
on the lowest weight space $K_{-\lambda }\otimes \SS_{-\rho }$, so the
only contribution from the first term is from the~$e_t$, and these terms
contribute~$\gamma (\lambda +\mu )$.  For the last term
  $$ \frac i3\gamma ^a\sigma _a = \frac i3 \sigma _a\gamma ^a + \frac i3
     f^a_{ac}\gamma ^c.  $$
On the lowest weight space $\sigma _{-\alpha }=\gamma ^\alpha =0$ for~$\alpha
\in \Delta ^+$.  So the root indices contribute
  $$ \sum\limits_{\alpha \in \Delta ^+}\frac i3 f^\alpha _{\alpha c}\gamma ^c
      = \sum\limits_{\alpha \in \Delta ^+}\frac i3 f^\alpha _{\alpha t}\gamma
     ^t
      = \sum\limits_{\alpha \in \Delta ^+}\frac i3 \,(-i)\alpha (e_t)\gamma ^t
      = \frac 13 \gamma (2\rho ) 
      = \frac 23 \gamma (\rho ).  $$
For the sum over torus indices we find $\frac i3 \gamma ^t\sigma _t= \frac
13\gamma (\rho )$, so altogether we have
  $$ \Dirac_\mu = \gamma (\mu +\lambda +\rho )\qquad \text{on $K_{-\lambda
     }\otimes \SS_{-\rho }$}.  $$
Squaring,
  $$ \Dirac_\mu^2 = -|\mu +\lambda +\rho |^2\qquad \text{on $K_{-\lambda
     }\otimes \SS_{-\rho }$},   $$
from which the proposition follows.
        \enddemo

The kernels fit together into a $G^\sigma $-equivariant vector bundle  
  $$ \Ker\Dirac(V)\res\orbit\longrightarrow \orbit  $$
over a $G$-coadjoint orbit $\orbit\subset \frak{g}^*$.  That $\orbit$~is a
single $G$-orbit follows from the irreducibility of~$V$, since $G/G\id$~acts
transitively on the weights occurring in~\thetag{1.18}.  Then~\thetag{1.20}
gives a decomposition
  $$ \Ker\Dirac(V)\res\orbit \;=\; K \otimes \sker\;\subset\; \orbit \times
     (V\otimes \SS),  $$
The bundle $K\to\orbit$ is $G$-equivariant.  If $G$~is connected, then as a
nonequivariant bundle $\sker\to N$~is the spinor bundle of the normal bundle
of~$\orbit$ in~$\frak{g}^*$; see~\thetag{1.16}.  As an equivariant bundle it
differs by a character of~$T$, and if $G$~is not connected possibly by a
finite twist, as we now explain.
 
First, the normal bundle~$N\to\orbit$ is $G$-equivariant and is equivariantly
trivial as a $G\id$-bundle, since the identity component~$T_\mu $ of the
stabilizer of~$\mu $ acts trivially on the fiber~$\frak{t}_\mu $ of~$N$
at~$\mu $.  The adjoint representation
  $$ \stab \longrightarrow \stab/T_\mu \longrightarrow \Or(\frak{t}_\mu )
      $$
factors through the finite group~$\stab/T_\mu $ of components of the
stabilizer.  The pullback of~\thetag{1.1} for~$H=\frak{t}_\mu $ gives a
central extension~$\stab^{\sigma (N)}$ of~$\stab$ by~$\TT$.  Fix an
irreducible $\zt$-graded complex $\Cliff^c(\frak{t}_\mu
^*)$-module~$\SS(\frak{t}_\mu ^*)$ and compatible metric; it is then a
unitary representation of~$\stab^{\sigma (N)}$.  Note that $\SS_{-\rho (\mu
)}$~is also an irreducible $\zt$-graded complex $\Cliff^c(\frak{t}_\mu
^*)$-module, but the extension~$Z_\mu ^\sigma $ acts.  We fix the ``sign''
(\S1.1) of the $\Cliff^c(\frak{t}_\mu ^*)$-module~$\SS(\frak{t}_\mu ^*)$ by
asking that it agree with that of~$\SS_{-\rho }$.  Let $\Cliff^c(N)$~be the
bundle of Clifford algebras on the normal bundle.  Define\footnote{If $\rank
N$~is odd, then elements of~$L$ are required to commute with the
$\cco$~action.}
  $$ L = \Hom_{\Cliff^c (N)}\bigl(\SS(N),\sker \bigr). \tag{1.23} $$
Then $L\to\orbit$~is an even line bundle.  We summarize in the following
proposition, which interprets the twistings in $K$-theory language.

        \proclaim{\protag{1.24} {Proposition}}
 Let~$\sh$ be the twisting\footnote{In terms of the notation~\thetag{0.2} for
twistings
  $$ \sh = \bigl(\frak{g}^* @>\{1\}>> \frak{g}^*, G^\sigma , \epsilon ^\sigma
     \bigr),  $$
where $\epsilon ^\sigma $~is the composition $G\to \Or(\frak{g})\to \zt$ of the
adjoint map and the nontrivial grading of the orthogonal group.}
of~$K_G(\frak{g}^*)$ induced by the central extension $G^\sigma \to G$; it
restricts to a twisting of~$K_G(\orbit)$.  Over the orbit~$\orbit\cong
G/\stab$ the graded central extension~$\stab^{\sigma (N)}$ defines a second
twisting of~$K_G(\orbit)$,\footnote{In terms of our model~\thetag{0.2} for
equivariant twistings
  $$ \sNh = \bigl(G @>Z_\mu >> G/Z_\mu ,G\times (Z_\mu ^{\sigma
     (N)})^{\text{op}},\epsilon ^{\sigma (N)}\bigr),  $$
where $G\times Z_\mu ^{\text{op}}$ acts on~$G$ by $\text{left}\times
\text{right}$ multiplication and $\epsilon ^{\sigma (N)}$~is the grading
of~$Z_\mu $ induced from the adjoint action.}  denoted~$\sNh$, and the bundle
$\SS(N)\to\orbit$ induced from~$\SS(\frak{t}_\mu ^*)$ is a $\sNh$-twisted
bundle.  The equivariant tangent bundle $T\orbit\to\orbit$ determines the
twisting~$\torbit=\sh -\sNh$ of~$K_G(\orbit)$.  The bundle $L\to\orbit$, and
so $K\otimes L\to\orbit$, is $\torbit$-twisted and is even.  If $G=G\id$~is
connected, then $\sNh$~is canonically trivial, there is a canonical
isomorphism~$\torbit\cong \sh $, and $L$~is canonically the $G^\sigma
$-bundle given by the character~$-\rho $.
        \endproclaim

\flushpar
 The second statement follows from the $G$-equivariant isomorphism
of~$N\oplus T\orbit\to\orbit$ with the trivial bundle~$\orbit\times
\frak{g}^*\to\orbit$.  Notice that~$L\to\orbit$, a rank one bundle, may be
viewed as a trivialization of the twisting~$\torbit$.  The latter is
trivializable since the coadjoint orbit admits $\text{spin}^c$, in fact spin,
structures; see the last paragraph of~\S{1.4}.

\theprotag{1.19} {Proposition} is an inverse to the Kirillov correspondence:
starting with an irreducible representation we construct a coadjoint
orbit~$\orbit$, the support of the Dirac family.  Note that the orbit is
$\rho $-shifted from the lowest weight of the representation (in the
connected case) and so is a regular orbit.  Furthermore, we construct a
vector bundle~$K\otimes L\to\orbit$.  As we will see in the next section, we
can recover~$V$ from this bundle.

 \subhead \S{1.4}.  $K$-theory interpretation
 \endsubhead
 
 A $G$-equivariant family of odd skew-adjoint $\zt$-graded Fredholm operators
parametrized by a $G$-space~$X$ represents a class in~$K^0_G(X)$.  If there
is a commuting $\cco$-structure, then it represents an element of~$K^1_G(X)$.
We apply this to~$X=\frak{g}^*$ with two variations to obtain a class
  $$ \bigl[\Dirac(V)\bigr] \in K_G^{\sh + \dim G}(\frak{g}^*)\cpt,
      $$
where we use Bott periodicity to identify~$K^n_G$ with~$K^0_G$ or~$K^1_G$
according to the parity of~$n$.  First variation: the family of
operators~$\Dirac(V)$ is compactly supported (\theprotag{1.19}
{Proposition}), that is, invertible outside a compact subset of~$\frak{g}^*$.
Second variation: $V\otimes S$~is a representation of the central
extension~$G^\sigma $ of~$G$, which accounts for the twisting~$\sh $.  Note
that if $\dim G $~is odd there is a $\cco$-structure on~$\SS$, hence also one
on~$V\otimes \SS$, and it commutes with~$\Dirac_\mu (V)$ for all~$\mu \in
\frak{g}^*$.

The Thom isomorphism in equivariant $K$-theory for the inclusion
$j\:\{0\}\hookrightarrow \frak{g}^*$ is usually implemented as the map 
  $$ \aligned
      j_*\:K^0_G &\longrightarrow K_G^{\sh + \dim G}(\frak{g}^*)\cpt \\
      V &\longrightarrow (V\otimes \SS, 1\otimes \gamma )\endaligned
     \tag{1.25} $$
which assigns to each finite dimensional unitary representation the family of
Clifford multiplication operators.  This family is supported at the origin.
The Dirac family~\thetag{1.13} is a compactly supported perturbation
of~\thetag{1.25}, so implements the same map on $K$-theory.  
 
Let $i\:\orbit\hookrightarrow \frak{g}^*$ be the inclusion, and consider the
diagram
  $$ \xymatrix@!C{& K_G^{\torbit+\dim\orbit}(\orbit) \ar[dl]_{\pi
     ^{\orbit}_*} \ar[dr]^{i_*} \\ K^0_G \ar@<.5ex>[rr]^{j_*} & & K_G^{\sh
     +\dim G}(\frak{g}^*)\cpt \ar@<.5ex>[ll]^{(\pi ^{\frak{g}^*})_*}}
     \tag{1.26} $$
of twisted $G$-equivariant $K$-groups.\footnote{The (twisted) $K$-cohomology
groups in the diagram are isomorphic to untwisted $G$-equivariant
$K$-homology groups in degree zero.  Also, recall that we use Bott
periodicity so only keep track of the parity of the $K$-groups.}  Since
$\orbit$~is a regular orbit, $\dim\orbit=\dim G - \dim T$ for any maximal
torus~$T$.  (There is a similar diagram for any $G$-coadjoint orbit, but for
our purposes it suffices to consider regular orbits.)  For any $G$-space~$X$
we denote $\pi ^X\:X\to \operatorname{pt}$ the unique map.  Then~$\SS$,
$\SS(N)$ represent equivariant Thom classes and so define the indicated
pushforwards maps~$j_*,(\pi ^{\frak{g}^*})_*,i_*$; we choose a compatible
twisted equivariant Thom class to define~$\pi ^\orbit_*$.  Then
  $$ i_* = j_*\circ \pi ^{\orbit}_*,  $$
since by naturality and the compatibility of the Thom classes,
  $$ \gathered
      (\pi ^{\frak{g}^*})_* \circ i_* = \pi ^{\orbit}_* \\
      \text{$j_*$ is the inverse of $(\pi ^{\frak{g}^*})_*$} \endgathered
     \tag{1.27} $$

        \proclaim{\protag{1.28} {Theorem}}
 \rom(i\rom)\ For an irreducible representation~$V$ of~$G$, 
  $$ j_*[V] = \bigl[\Dirac(V)\bigr]. \tag{1.29} $$
 \rom(ii\rom)\ Let $[K\otimes L]\in K_G^{\sigma (\orbit)}(\orbit)$ be the
class of the vector bundle $K\otimes L\to\orbit$.  Then 
  $$ \align
      i_*[K\otimes L] &= \bigl[\Dirac(V)\bigr] \tag{1.30} \\
      \pi ^{\orbit}_*[K\otimes L] &= [V] \tag{1.31} \endalign $$
        \endproclaim

\flushpar
 The twisting of~$K\otimes L$ was described in \theprotag{1.24}
{Proposition}.

        \demo{Proof}
See the remarks following~\thetag{1.25} for~(i).  
 
The map~$i_*$ is multiplication by the Thom class of the normal bundle
to~$\orbit$ in~$\frak{g}^*$.  This gives a class in a tubular neighborhood
of~$\orbit$ which is compactly supported; we extend by zero to a compactly
supported class on~$\frak{g}^*$.  To see that
$\bigl[\Dirac(V)\bigr]$~represents this class, as asserted in~\thetag{1.30},
note first that $[\Dirac(V)]$~is supported on~$\orbit$ by \theprotag{1.19}
{Proposition}.  Hence it suffices to restrict to a tubular neighborhood~$U$
of~$\orbit$, which we identify with a neighborhood of the zero section in the
normal bundle~$N\to\orbit$ via orthogonal projection.  Write~$V\otimes
\SS\res U$ as the direct sum of $K\otimes L\otimes \SS(N)$, pulled back
from~$\orbit$ using the tubular neighborhood structure, and its orthogonal
complement.  We may ignore the latter since $\Dirac(V)$~is invertible on it
if $U$~is sufficiently small.  For any~$\mu \in \orbit$ and $\theta \in
\frak{g}^*$ perpendicular to~$\mu $ we have
  $$ \Dirac_{\mu +\theta } = \Dirac_\mu + \gamma (\theta ).  $$
Restricted to~$K\otimes L\otimes \SS(N)$ it acts as~$1\otimes 1\otimes \gamma
(\theta )$.  Since $\theta \mapsto \gamma (\theta )$ on~$\SS(N)$ is the Thom
class of the normal bundle, we have proved~\thetag{1.30}.
 
Equation~\thetag{1.31} is a formal consequence of~\thetag{1.29}
and~\thetag{1.30} using~\thetag{1.27}.  
        \enddemo

The Atiyah-Singer index theorem identifies~$\pi _*^{\orbit}$ as the
equivariant index of a Dirac operator on~$\orbit$.  Recall that the
$G$-coadjoint orbit $\orbit$~is a finite union of regular $G\id$-coadjoint
orbits, each of which is diffeomorphic to~$G\id/T$.  The latter carries a
canonical spin structure up to isomorphism;\footnote{Write $G\id/T \cong
\tilde{G}\id/\tilde{T}$ for $\tilde{G}\id$~the simply connected cover
of~$G\id$.  A choice of Weyl chamber gives a $\tilde{G}\id$-invariant complex
structure on~$\tilde{G}\id/\tilde{T}$ whose first Chern class
in~$H^2(\tilde{G}\id/\tilde{T})\cong \tilde{\frak{t}}^*$ is identified
with~$2\rho $.  Since $\rho $~is a character of~$\tilde{G}\id$, we conclude
that $\tilde{G}\id/\tilde{T}$~is spin, and since it is simply connected the
spin structure is unique up to isomorphism.  More precisely, the component of
any~$\mu \in \orbit$ is canonically~$G\id/T_\mu $ and $\frak{t}_\mu $~ has a
distinguished Weyl chamber which gives a distinguished complex structure
on~$\orbit$.  In particular, there is a canonical orientation of~$\orbit$.}
the twist~$\torbit$ in~\thetag{1.26} means that we consider bundles
on~$\orbit$ which are equivariant for a covering group of~$G$.  In case
$G$~is connected and simply connected, both~$\sh $ and~$\sNh$ are isomorphic
to the trivial twisting, whence $\pi _*^{\orbit}$~is the usual
$G$-equivariant Dirac index.  The map~$\pi _*^{\orbit}$ is termed {\it Dirac
induction\/} in representation theory.  Our construction in this section is
an explicit inverse to Dirac induction.  Equation~\thetag{1.31} is part of
the Borel-Weil-Bott theorem, but we do not prove the vanishing theorem
necessary to show that the kernel of the Dirac operator is precisely~$V$.
(Note that the usual version of this theorem uses holomorphic induction in
place of Dirac induction.)

 \subhead \S{1.5}.  Variation: projective representations
 \endsubhead

In the loop group version of the Dirac family we encounter the finite
dimensional construction considered here, but in a projective version which
we now sketch. 
 
Let $1\to \TT \to G^\tau \to G\to 1$ be a central extension
of~$G$.\footnote{Such extensions are classified by~$H^3(BG;\ZZ)$, which is
torsion for a compact Lie group~$G$.  Therefore, a central extension by~$\TT$
factors through a central extension by a finite cyclic group.}  A {\it
projective representation of~$G$ at level~$\tau $\/} is a representation
of~$G^\tau $ on which the center~$\TT$ acts by scalar multiplication.  The
free abelian group~$K_G^\tau $ generated by the irreducibles is a
$K_G$-module via tensor product.  In this situation $\frak{g}^*$~is replaced
by an affine space~$\aff G\tau $ for~$\frak{g}^*$ which has two equivalent
definitions: it is the space of linear splittings of the Lie algebra
extension~$0\to i\RR\to\frak{g}^\tau \to \frak{g}\to 0$, or equivalently it
is the subspace of~$(\frak{g}^\tau )^*$ consisting of functionals~$\mu
\:\frak{g}^\tau \to\RR$ with $\mu (i)=1,\;i\in i\RR$.  The group~$G$ acts
on~$\aff G\tau $ compatibly with the coadjoint action on~$\frak{g}^*$.  Note
that the action of~$G$ on~$\aff G\tau $ has a fixed point, the center of mass
of a $G$-orbit, from which $\aff G\tau \cong \frak{g}^*$ as $G$-spaces.
(This is in contrast to the infinite dimensional situation we study later.)
 
Let $T\subset G$ be a maximal torus.  By restriction we obtain a central
extension $T^\tau \to T$, and $T^\tau$~is a maximal torus of~$G^\tau $.  The
affine space~$\aff T\tau $ sits inside~$\aff G\tau $ as the subspace of
functionals which vanish on the root spaces.  The infinitesimal characters of
projective representations of~$T$ at level~$\tau $ form a subset $\Lambda
^\tau \subset \aff T\tau $ which is a torsor for the weight lattice~$\Lambda
\subset \frak{t}^*$.  The Weyl group acts on~$\aff T\tau $
preserving~$\Lambda ^\tau $; the action is conjugation by the normalizer
of~$T^\tau $ in~$G^\tau $.  Choose a positive Weyl chamber in~$\frak{t}$, and
so by pullback a positive Weyl chamber in~$\frak{t}^\tau $.  (The Weyl groups
of~$G$ and~$G^\tau $ are canonically isomorphic.)  
 
Given central extensions~$G^{\tau _1}\to G$ and $G^{\tau _2}\to G$ we can
form their product~$G^{\tau _1+\tau _2}\to G$, which is also a central
extension.  Similarly, we define an inverse to a central extension.  Now the
Lie algebra extension~\thetag{1.5} associated to the canonical spin central
extension~$G^\sigma \to G$ is canonically split,\footnote{In the infinite
dimensional case (\S{3}) there is no such canonical splitting.} so for
any~$\tau $ we identify~$\aff G\tau $ with~$\aff G{\tms}$.
 
The projective version of our Dirac family proceeds as follows.  Fix a
central extension~$G^\tau \to G$ and a projective representation~$V$ at
level~$\tms$.  So $V\otimes \SS$~is a projective representation at
level~$\tau $.  The Dirac family~\thetag{1.12} is parametrized by~$\mu \in
\aff G\tau \cong \aff G\tms$:
  $$ D_\mu = i\gamma ^a(R_a)_{\mu } + \frac i2\gamma (\Omega ), \tag{1.32} $$
where $(R_a)_\mu $~denotes the infinitesimal action of the basis
element~$e_a$ on~$V$ defined using the splitting $\mu \:\frak{g}\to
\frak{g}^{\tms}$.  Then 
  $$ D_{\mu +\nu }= D_\mu + \gamma (\nu ),\qquad\nu \in
     \frak{g}^*. \tag{1.33} $$
The discussion proceeds as before, and in particular \theprotag{1.19}
{Proposition} holds, but with the understanding that $\mu \in \aff G\tau $,
$\rho \in \aff G\sigma $, and $\lambda \in \aff G\tms$.  Also, in the
$K$-theory discussion we use the fact that the Lie algebra
extension~$\frak{g}^\tau \to \frak{g}$ is split to find a fixed point of the
$G$-action on~$\aff G\tau $ and so define the Thom isomorphism.

\newpage
 \subhead \S{1.6}.  Examples
 \endsubhead

        \example{\protag{1.34} {Example}}
 ($G=\TT$)\quad The irreducible representations are all one-dimensional and
are labeled by an integer~$n$: the representation~$L_n=\CC$ is~$\lambda
\mapsto \lambda ^n$.  Identify~$\frak{g}=i\RR$ and set\footnote{This inner
product generalizes for~$G=U(n)$ to the inner product $\langle A,A' \rangle =
-\Tr(AA')$ on skew-Hermitian matrices.} $\langle ia,ia' \rangle = aa'$,
thereby identifying $\frak{g}^*\cong i\RR$ as well.  Set $\SS=\CC\oplus \CC$
with $\gamma =\left(\smallmatrix 0&i\\i&0 \endsmallmatrix\right)$.  Then for
$\mu =ia\in \frak{g}^*$ the operator~$D_\mu $ on~$L_n\otimes \SS$ is the
matrix $D_{ia} = \left(\smallmatrix 0&i(a-n)\\i(a-n)&0
\endsmallmatrix\right)$.  The kernel of the family~$D(L_n)$ is supported
at~$\mu =in$.  Note that the coadjoint action is trivial, since $G$~is
abelian, so the coadjoint orbits are points.
        \endexample

        \example{\protag{1.35} {Example}}
 ($G=\Or_2$)\quad Now $G$~is not connected.  We
identify~$\frak{g},\frak{g}^*$ as in the previous example.
Reflections---elements of~$G$ not in the identity component---act as
multiplication by~$-1$ on~$\frak{g}^*$, so the coadjoint orbits are~$\{0\}$
and~$\{\pm ia\}_{a\not= 0}$.  There is a one-dimensional irreducible
representation~$V_\epsilon $ which sends a matrix to multiplication by its
determinant, and there are two-dimensional irreducibles~$V_n$ labeled by {\it
positive\/} integers: the restriction of~$V_n$ to~$\SO_2\subset \Or_2$
is~$L_n\oplus L_{-n}$ in the notation of \theprotag{1.34} {Example}.
For~$V_\epsilon $ the infinitesimal representation vanishes, and
$D_{ia}=\left(\smallmatrix 0&ia\\ia&0 \endsmallmatrix\right)$ acting
on~$V_\epsilon \otimes \SS$.  The kernel is supported at~$\mu =0$ and is
obviously~$V_\epsilon \otimes \SS$.  We can identify~\thetag{1.23}
with~$V_\epsilon $, and the stabilizer group~$Z_0=\Or_2$ acts via the
representation~$V_\epsilon $.  For~$V_n$ we identify
  $$ D_{ia} = \pmatrix 0&0&0&i(a-n)\\
      0&0&i(a+n)&0\\
      0&i(a-n)&0&0\\
      i(a+n)&0&0&0 \endpmatrix.  $$
So $\Ker D(V_n)$~is supported on the coadjoint orbit~$\mu =\pm in$. 
        \endexample

        \example{\protag{1.36} {Example}}
 ($G=\zt\times \TT$)\quad Again $G$~is not connected, but now $G$~is abelian so
the irreducibles are one-dimensional and the coadjoint orbits are points.
Let $K$~denote the sign representation of~$\zt$; then the irreducibles of~$G$
are~$L_n$ and~$L_n\otimes K$, where $n$~ranges over all integers.  As in
\theprotag{1.34} {Example} the kernel of the Dirac family is supported
at~$\mu =in$ for both~$D(L_n)$ and~$D(L_n\otimes K)$.  They are distinguished
by the action of the stabilizer~$Z_\mu =G$ on the kernel. 
        \endexample

        \example{\protag{1.37} {Example}}
 ($G=\SU_2$)\quad There is a single irreducible representation in each positive
dimensions, and in a notation adapted to the quotient group~$\SO_3$ we label
the representation~$V_j$ of dimension~$2j+1$ by~$j=0,\frac 12,1,\dots $.  The
Lie algebra is $\frak{g}=\left\{ \left(\smallmatrix ia&b+ic\\-b+ic&-ia
\endsmallmatrix\right) \right\} $, and using the inner product $\langle A,A'
\rangle = -\frac 12\Tr(AA')$ we identify $\frak{g}\cong \frak{g}^*\cong
i\RR^3$.  We take as maximal torus is the group of diagonal matrices with Lie
algebra~$i\RR$.  Then the positive root is identified as~$2i$ from
which~$\rho =i$.  The coadjoint orbits are spheres centered at the origin.
\theprotag{1.19} {Proposition} states that $\Ker D(V_j)$~is supported on the
sphere of radius~$2j+1$.  Notice that even for~$j=0$ the support is a {\it
regular\/} coadjoint orbit, due to the $\rho $-shift. 
        \endexample

        \example{\protag{1.38} {Example}}
 ($G=\SO_3$)\quad Only the representations of~$\SU_2$ with $j$~integral pass to
the quotient~$\SO_3$.  Take as maximal torus~$T$ the set of rotations which
fix a particular axis.  The inner product is~$\langle A,A' \rangle = -\frac
12\Tr(AA')$.  With the identification~$\frak{t}=i\RR$ we find~$\rho =\frac
i2$.  Note that $\rho $~is {\it not\/} a weight of~$G$, which is possible as
$G$~is not simply connected.  Also, $\Ker D(V_j)$~is supported on the sphere
of radius~$j+\frac 12$, where now~$j\in \ZZ^{\ge 0}$.
        \endexample

        \example{\protag{1.39} {Example}}
 ($G^\tau =\Ur_2\to G=\SO_3$)\quad Notice that $\Ur_2$~has double
cover~$\SU_2\times \TT$, and the projective representations of interest pull
back to representations of~$\SU_2\times \TT$ which are of half-integral spin
on~$\SU_2$ and standard on~$\TT$.  Let~$T$ be the maximal torus of~$\SO_3$ as
in the previous example.  We may identify the affine space~$\aff T\tau $ with
~$i\RR$ and the nontrivial element of the Weyl group acts by~$ia\mapsto -ia$.
It has no fixed points on the torsor~$\Lambda ^\tau $ of imaginary
half-integers.  The Dirac family for the spin~$j$ representation has kernel
supported on the sphere of radius~$j+\frac 12$ in~$\aff {\SO_3}\tau $.
        \endexample

 \head
 \S{2} Loop groups and energy
 \endhead
 \comment
 lasteqno 2@ 51
 \endcomment

\subhead \S{2.1}.  Basic definitions
 \endsubhead

Let $G$~be a compact Lie group and $P\to\cir$ a principal $G$-bundle with
base the standard circle~$\cir=\RR/2\pi \ZZ$.  The bundle~$P$ is classified
up to isomorphism by a conjugacy class in~$\pi _0(G) = G/G\id$.  Let
$\GR\subset G$ denote the corresponding union of components, which is a
$G$-space under conjugation.  Note that if $G$~is connected, then $P$~is
necessarily trivializable and $\GR=G$.

        \example{\protag{2.1} {Notation}}
 $\LGR$~is the group of smooth gauge transformations of $P\to\cir$.  Its Lie
algebra is denoted~$\LgR$.  Let $\Trot$~be the group of rigid rotations
of~$\cir$ and $\LGRRz$~be the group of smooth automorphisms of $P\to\cir$
which cover elements of~$\Trot$.  For any connected finite covering
$\Troth\to\Trot$ we define $\LGRR$ by pullback: 
  $$ \CD
      \LGRR @>>> \Troth \\
      @VVV @VVV \\
      \LGRRz @>>> \Trot\endCD \tag{2.2} $$
Finally, $\AR$~is the space of smooth connections on $P\to\cir$.
        \endexample

\flushpar
 For the trivial bundle $P=\cir\times G$ we omit the subscript~`$P$'.  Then
$LG$~is naturally identified with the smooth loop group~$LG=\Map(\cir,G)$.
For topologically nontrivial bundles $\LGR$~is often termed a {\it twisted
loop group\/}.  Its Lie algebra~$\LgR$ may be identified with~$\Omega
^0_{\cir}(\frak{g}\mstrut _P)$, the space of smooth sections of the adjoint
bundle $\frak{g}\mstrut _P\to\cir$, or equivalently as the space of smooth
$G$-invariant vertical vector fields on $P\to\cir$.  The space~$\AR$ of
connections is affine with associated vector space~$\Omega
^1_{\cir}(\frak{g}\mstrut _P)$.
 
Fix a basepoint~$p\in P$ in the fiber over~$0\in \RR/2\pi \ZZ$.  Then for
each connection~$A\in \AR$ parallel transport around~$\cir$ maps~$p$ to a
point~$p\cdot \hol(A)$ of the same fiber for some~$\hol(A)\in G$.  The {\it
holonomy\/} is a surjective map
  $$ \hol\:\AR\longrightarrow \GR. \tag{2.3} $$
Furthermore, the holonomy map is equivariant for the natural action of~$\LGR$
on~$\AR$ and the conjugation action of~$G$ on~$G[P]$ as follows.  For
$\varphi \:P\to P$ in~$\LGR$ and~$A\in \AR$ we have a new
connection~$(\varphi \inv )^*(A)$ and a group element~$g_{\varphi }\in G$
defined by~$\varphi (p)=p\cdot g\inv _{\varphi }$.  (Recall $A$~is a 1-form
on~$P$.)  Then
  $$  \hol\bigl((\varphi \inv )^*(A) \bigr) = g_\varphi
     \hol(A)g_{\varphi }\inv.   $$
The subgroup $\OGR\subset \LGR$ of gauge transformations which fix~$p\in P$
acts freely on~$\AR$, and \thetag{2.3} is a (left) $\OGR$-principal bundle.
 
There is an exact sequence of groups 
  $$ 1\longrightarrow \LGR \longrightarrow \LGRR \longrightarrow
     \Troth\longrightarrow 1  $$
and a corresponding sequence of Lie algebras 
  $$ 0\longrightarrow \LgR \longrightarrow \LgRR \longrightarrow
     i\Rrot\longrightarrow 0. \tag{2.4} $$
(The Lie algebras of~$\Troth$ and~$\Trot$ are canonically isomorphic.)
Splittings of~\thetag{2.4} are in 1:1~correspondence with~$\AR$.  For~$A\in
\AR$ we denote by~$d_A$ the lift to~$\LgRR$ of the canonical generator $i\in
i\Rrot$.  Geometrically, $d_A$~is a $G$-invariant horizontal vector field
on~$P$ which projects to the standard vector field on~$\cir$.  We remark that
the extended loop group~$\LGRR$ may be used to eliminate the basepoint~$p$ in
the map~\thetag{2.3} above.\footnote{Namely, if we include the choice of
basepoint in the domain of the holonomy map, then
  $$ \hol\:\AR\times P\longrightarrow \GR  $$
is a principal $\LGRR$-bundle and admits a commuting $G$-action.}

The Lie algebra~$\frak{g}$ of a compact Lie group~$G$ has a canonical
splitting 
  $$ \frak{g} = \cent(\frak{g})\oplus [\frak{g},\frak{g}] = \frak{z}\oplus
     \frak{g}' \tag{2.5} $$
into its abelian and semisimple pieces; the decomposition is orthogonal with
respect to any invariant inner product on~$\frak{g}$.  This decomposition is
Ad-invariant, so leads to a decomposition of the adjoint bundle
$\frak{g}_P\to\cir$ of any principal $G$-bundle $P\to\cir$, whence to a
decomposition of the loop algebra~$\LgR$.

 \subhead \S{2.2}.  Admissible central extensions
 \endsubhead

The representations of~$\LGR$ of interest are projective, and so we study
central extensions
  $$ 1 \longrightarrow \TT \longrightarrow \LGRt \longrightarrow \LGR
     \longrightarrow 1 \tag{2.6} $$
of loop groups with corresponding central extensions 
  $$ 0 \longrightarrow i\RR \longrightarrow \LgRt \longrightarrow \LgR
     \longrightarrow 0 \tag{2.7} $$
of loop algebras.  If $P$~is trivial we write $LG^\tau $ for the central
extension.  Recall that a {\it graded\/} central extension of~$\LGR$ is a
central extension~\thetag{2.6} together with a homomorphism $\epsilon
\:\LGR\to\zt$.

        \definition{\protag{2.8} {Definition}}
 Let $(\LGRt,\epsilon )$ be a graded central extension of~$\LGR$.  The {\it
associated twisting\/} of the $G$-equivariant $K$-theory of~$\GR$ is
  $$ \th = (\AR@>\OGR>> \GR,\LGRt,\epsilon), \tag{2.9} $$
where $G$~acts on~$\GR$ by conjugation.
        \enddefinition

\flushpar
 Here $\AR\to\GR$ is the principal $\OGR$-bundle~\thetag{2.3} and the
twisting has the form~\thetag{0.2}.  This construction of a twisting from a
graded central extension of the loop group is fundamental to our work.

Let\footnote{Our~$K$ and $d_A$ (below) are $i=\sqrt{-1}$~times real multiples
of the $K$ and~$d$ in~\cite{Ka}.  (Compare~\cite{Ka,(7.2.2)} to~\thetag{2.21}
and~\thetag{2.22} below.)  This motivates the sign in~\thetag{2.12}.}
$K$~denote the central element~$i\in i\RR$.

        \definition{\protag{2.10} {Definition}}
 A central extension~$\LGRt$ is {\it admissible\/} if
 \roster
 \item There exists a central extension~$\LGRRt$ of~$\LGRR$ which fits into
a commutative diagram 
  $$ \CD
      @. @. 1 @. 1 \\
      @. @. @VVV @VVV \\
      1 @>>> \TT @>>> \LGRt @>>> \LGR @>>> 1\\
      @. @| @VVV @VVV \\
      1 @>>> \TT @>>> \LGRRt @>>> \LGRR @>>> 1\\
      @. @. @VVV @VVV \\
      @. @. \Troth @= \Troth \\
      @. @. @VVV @VVV \\
      @. @. 1 @. 1 \endCD \tag{2.11} $$

 \item There exists an $\LGRRt$-invariant symmetric bilinear form~$\form_\tau
$ on~$\LgRRt$ such that
  $$ \ll\! K,{d}\!\gg_\tau  \;=\; -1 \tag{2.12} $$
for all ${d}\in \LgRRt$ which project to~$i\in i\Rrot$.
 \endroster
        \enddefinition

\flushpar
 It follows immediately from~\thetag{2.12} that
  $$ \ll\! K,\LgRt\!\gg_\tau  \;=0, \tag{2.13} $$
and in particular $\ll\! K,K \!\gg_\tau \;=0$.  Hence $\ll\!\cdot ,\cdot
\!\gg_\tau $~induces an $\LGR$-invariant symmetric bilinear form on~$\LgR$,
which we also denote with double angle brackets.  For $A\in \AR$ the
element~$d_A\in \LgRR$ has a unique lift to a null element of~$\LgRRt$ which
we denote~$d_A^\tau $:
  $$ \ll\! d_A^\tau ,d_A^\tau  \!\gg_\tau  \;=0. \tag{2.14} $$

In the sequel the phrase ``$\LGRt$~is an admissible central extension''
implies that we have fixed~$\LGRRt$ and~$\ll\!\cdot ,\cdot \!\gg_\tau $.  If
there is no possibility for confusion, we write $\form $ for the invariant
bilinear form.  We remark that the finite covering map $\Troth\to\Trot$ is
often the identity; see the discussion of tori in~\S{2.3} for a case in which
it is a double cover.  The choice of~$\LGRRt$ and a bilinear form is
additional structure, as reflected in \theprotag{2.18} {Lemma} below, for
example.

If $\frak{g}$~is semisimple, then every central extension of~$\LGR$ is
admissible.

        \proclaim{\protag{2.15} {Proposition}}
 Let $G$~be a compact Lie group with $[\frak{g},\frak{g}]=\frak{g}$, i.e.,
$\frak{g}$~semisimple.  Let $P\to\cir$ be a principal $G$-bundle.  Then any
central extension~$\LGRt$ is admissible.  Furthermore, for each~$\LGRRt$
which satisfies~\thetag{2.11} there exists a unique $\LGRRt$-invariant
symmetric bilinear form $\form\;=\;\form_\tau $ which
satisfies~\thetag{2.12}.
        \endproclaim

\flushpar
 For connected and simply connected groups there is a direct construction,
which we give in~\S{2.3}.  The general proof is more complicated and is
deferred to the appendix.

The Lie algebra~$\frak{g}$ of any compact Lie group~$G$ carries a canonical
$G$-invariant symmetric bilinear form~$\langle \cdot ,\cdot   \rangle_\sigma
$ defined by 
  $$ \langle \xi ,\eta \rangle_\sigma =-\frac 12\Tr(\ad\xi \circ
     \ad\eta ),\qquad \xi ,\eta \in \frak{g}.  $$
Note that $\langle \cdot ,\cdot \rangle_\sigma $~is positive semidefinite,
and is positive definite if $\frak{g}$~is semisimple.  Now if $P\to\cir$ is a
principal $G$-bundle, then we define the $L^2$~metric
  $$ \ll\!\beta_1 ,\beta_2 \!\gg_\sigma =\int_{\cir}\langle \,\beta_1
     (s),\beta_2 (s)\, \rangle_\sigma \;\frac{|ds|}{2\pi },\qquad \beta_1
     ,\beta_2 \in \LgR. \tag{2.16} $$
We will prove (\theprotag{3.13} {Proposition}) that \thetag{2.16}~is
associated to a canonical admissible graded central extension of~$\LGR$.

Let $\LGRt\to\LGR$ be any admissible central extension.  The space of linear
splittings of~\thetag{2.7} carries an affine action of~$\LgR$.  Namely, if
  $$ \aligned
      s\:\LgR &\longrightarrow \LgRt \\
      \beta &\longmapsto \beta ^\tau _s\endaligned  $$
is a splitting and~$\beta' \in \LgR$, then we define a new splitting~$\beta
\mapsto \beta ^\tau _{s+\beta' }$ by
   $$ \beta ^\tau _{s+\beta' } = \beta ^\tau_s - \ll\! \beta ,\beta'
     \!\gg_\tau \,K.  $$
Recall that we identify~$\AR$ as the space of $G$-invariant vector
fields~$d_A$ on~$P\to \cir$ which project to the standard vector field
on~$\cir$.  In this form $\AR$~ is affine for the action of~$\LgR$ by
subtraction: for $\xi \in \LgR$ we have
  $$ d_{A+\xi ds} = d_A - \xi .  \tag{2.17} $$

        \proclaim{\protag{2.18} {Lemma}}
 Let~$\LGRt$ be an admissible central extension of~$\LGR$.  Then the form
$\form_{\tau }$ determines an $\LgR$-equivariant map
  $$ \aligned
      \AR &\longrightarrow \{\text{linear splittings of $\LgRt \to \LgR$}\} \\
      A &\longmapsto (\beta \mapsto \beta _A^\tau )\endaligned  $$
which is an isomorphism if $\form_\tau $~is nondegenerate.
        \endproclaim

\flushpar
 Therefore, a connection simultaneously splits all admissible central
extensions.

        \demo{Proof}
 Characterize~$\beta _A^\tau $ by the condition 
  $$ \ll\! \beta _A^\tau ,d_A^\tau \!\gg_\tau \;=0 \tag{2.19} $$
for $d_A^\tau \in \LgRRt$ defined by~\thetag{2.14}.
        \enddemo

Suppose $(\LGR)^{\tau _i}\to \LGR,\,i=1,2$ are admissible central extensions.
Then the  product $(\LGR)^{\tau _1+\tau _2} = (\LGR)^{\tau _1}\times
_{\LGR}(\LGR)^{\tau _2}$ is also admissible: the bilinear forms on~$\LgR$
satisfy
  $$ \ll\! \beta_1 ,\beta_2 \!\gg_{\tau _1+\tau _2} \;= \;\ll\! \beta_1
     ,\beta_2 \!\gg_{\tau _1} +\; \ll\! \beta_1 ,\beta_2 \!\gg_{\tau
     _2},\qquad \beta_1 ,\beta_2 \in \LgR.  $$
There is also an (admissible) inverse central extension~$(\LGR)^{-\tau }$ to an
(admissible) central extension~$(\LGR)^{\tau }$; in the admissible case the
form changes sign when passing to the inverse.

 \subhead \S{2.3}. Special cases
 \endsubhead

We treat the three prototypical classes of compact Lie group:  connected
simply connected groups, tori, and finite groups.

 \subsubhead Simply connected groups\endsubsubhead
 If $G$~is simply connected, then every central extension~\thetag{2.6} is
admissible and $\form$~is uniquely determined by the extension.  This is a
special case of \theprotag{2.15} {Proposition}.  To verify this we may as
well assume $P$~is trivial so that $\LGR=LG$~is the standard loop group.
Then~\cite{PS,\S4.2} central extensions of the loop algebra~$L\frak{g}$
correspond to $G$-invariant symmetric bilinear forms~$\langle \cdot ,\cdot
\rangle$ on~$\frak{g}$ as follows. For~$\xi \in \gC$ and~$n\in \ZZ$ we
write~$z^n\xi$ for the loop $s\mapsto e^{ins} \xi$ in~$L\frak{g}_{\CC}$.  The
algebraic direct sum $\oplus _{n\in \ZZ}z^n\frak{g}_{\CC}$~is dense
in~$L\frak{g}_{\CC}$.  Define
  $$ (\widehat L\frak{g})^\tau =  \RR K\oplus L\frak{g}\oplus \RR d
     \tag{2.20} $$
as a vector space.  As for the Lie bracket, $K$~is central, $d$~acts as the
derivation
  $$ [d,z^n\xi ]=in z^n\xi \tag{2.21} $$
on~$L\frak{g}_\CC$, and the bilinear form on~$\frak{g}$ enters into the
bracket
  $$ [z^n\xi, z^{m}\eta] = z^{n+m}[\xi ,\eta ]+ n\langle \xi ,\eta \rangle
     \delta _{n+m=0}\;\frac Ki,\qquad \xi ,\eta \in
     \frak{g}_{\CC}. \tag{2.22} $$
Any $\widetilde{L\frak{g}}^\tau $-invariant bilinear form~$\ll\! \cdot ,\cdot
\!\gg$ satisfies
  $$ \ll\! [ z^n\xi, z^{-n}\eta],d\!\gg \;=\; \ll\! z^n\xi,[
     z^{-n}\eta,d]\!\gg\;=\; in\!\ll\! z^n\xi, z^{-n}\eta\!\gg, \tag{2.23} $$
and from~\thetag{2.22} and~\thetag{2.12} we deduce that for~$n\not= 0$ 
  $$ \ll\! z^n\xi, z^{-n}\eta\!\gg\;=\langle \xi ,\eta \rangle. \tag{2.24} $$
For~$n=0$ we derive~\thetag{2.24} by consideration of~$\ll\![z\zeta _1,z\inv
\zeta _2],\eta \!\gg$ and the fact that~$[\frak{g},\frak{g}]=\frak{g}$.
Similarly, we find $\ll\!z^n\xi ,z^m\eta \!\gg\;=0$ if~$n+m\not= 0$.  Thus
the bilinear form on~$L\frak{g}$ is necessarily the $L^2$~metric
  $$ \ll\!\beta_1 ,\beta_2 \!\gg\;=\int_{\cir}\langle\, \beta_1 (s),\beta_2
     (s)\, \rangle\;\frac{|ds|}{2\pi },\qquad \beta_1 ,\beta_2 \in
     L\frak{g}. \tag{2.25} $$
Notice that this form is symmetric and is also $LG$-invariant. 
 
For $G$~connected and simply connected 
  $$ H^2(LG;\RR)\cong \Sym^2(\frak{g}^*)^G,  $$
the correspondence given by associating to the $G$-invariant symmetric
form~$\langle \cdot ,\cdot   \rangle$ the left-invariant 2-form 
  $$ \omega (\beta_1 ,\beta_2 ) = \int_{\cir}\langle \beta_1 ,d\beta_2
     \rangle,\qquad \beta_1 ,\beta_2 \in L\frak{g},  $$
on~$LG$.  Note that $\omega $~is the cocycle~\thetag{2.22} which defines the
central extension of~$L\frak{g}$, and for a central extension~\thetag{2.6}
of~$LG$ we interpret it as $1/i$~times the curvature of the left-invariant
connection given by the splitting~\thetag{2.20}.  In particular, it lies in a
full lattice of forms determined by~$H^2(LG;2\pi \ZZ)\subset H^2(LG;\RR)$,
and this lattice classifies central extensions of~$LG$ (see~\cite{PS,\S4.4}).
Any such extension extends to~$\widehat{L}G=LG\rtimes \Trot$.

 \subsubhead Tori\endsubsubhead
 Let $T$~be a torus group.  Write $\frak{t}=\Lie(T)$ and $\Pi =\exp\inv
(1)/2\pi \subset \frak{t}$, so that $T\cong \frak{t}/\Pi $.  The loop group
has a canonical decomposition
  $$ LT \cong T\times \Pi \times U, \tag{2.26} $$
where 
  $$ U = \exp V,\qquad V=\left\{\beta \:\cir\to\frak{t}:\int_{\cir}\beta
     (s)\,ds=0\right\}.  $$
Note that the vector space $V\subset L\frak{t}$~is isomorphic
to~$L\frak{t}/\frak{t}$.  We can interpret~$T\times \Pi $ as the Lie group of
parametrized closed geodesics on~$T$ with respect to an invariant metric.  An
element~$X \in \Pi $ corresponds to the one-parameter group~$\varphi _X
(s)=\exp(sX )$.  Let $\Lambda =\Hom(T,\TT)$ be the $\ZZ$-dual of~$\Pi $.

        \proclaim{\protag{2.27} {Proposition}}
 Suppose $LT^\tau $~is an admissible central extension of~$LT$.  Then 
 	\roster
 \item"(i)" $LT^\tau $ is the product of central extensions~$(T\times \Pi
)^{\tau }$ and~$\Ltt^{\tau}$.
 \item"(ii)" Write $[\beta _1,\beta _2]=\omega (\beta _1,\beta _2)K$ for
$\beta _1,\beta _2\in L\frak{t}$.  Then 
  $$ \omega (\beta _1,\beta _2)=\ll\!\beta _1,\dot\beta _2
     \!\gg, \tag{2.28} $$ 
where $\dot\beta $ is the derivative~$d\beta /ds$.
 \item"(iii)" The commutator in~$(T\times \Pi )^{\tau }$ determines a
homomorphism $\kappa \:\Pi \to\Lambda $ and a bilinear form~$\langle \cdot
,\cdot   \rangle$ on~$\frak{t}\times \frak{t}$.  The form~$\langle \cdot
,\cdot   \rangle$~is symmetric.
	\endroster
        \endproclaim

The central extension~$U^{\tau }$ of the vector space~$U$ is a {\it
Heisenberg group\/} if the commutator pairing~\thetag{2.28} is nondegenerate.
Conditions~(i) and~(ii) are stated as {\it analytic regularity\/} in
~\cite{FHT3,\S2}; here they follow from the existence of~$\form$.  For~(iii)
observe that the restriction $T^\tau \to T$ of $LT^\tau \to LT$ is
necessarily split.  Let $t\mapsto \tilde{t}$ be a splitting.  Choose a lift
$\tilde \varphi _X\in LT^\tau $ of $\varphi _X\in LT$ for each~$X\in \Pi $.
Then $\kappa_X\in \Lambda= \Hom(T,\TT)$ is defined by
  $$ \tilde{\varphi }_X\tilde{t}\tilde{\varphi }_X\inv \tilde{t}\inv
     =\kappa_X (t). \tag{2.29} $$
Let ~$\dot\kappa_X\in \Hom(\frak{t},\RR)$ be $1/i$~times the derivative
of~$\kappa_X$; it defines a bilinear form 
  $$ \langle X ,\xi \rangle =-\dot\kappa_X(\xi )  $$
which extends to all of~$\frak{t}\times \frak{t}$ and takes integral values
on~$\Pi \times \Pi $.  (The sign is for convenience.)

        \demo{Proof}
 The assertion in~(ii) follows immediately from the invariance of~$\form$, as
in~\thetag{2.23}.  This implies that the constant loops are orthogonal to
loops which are derivatives, i.e., to elements of~$V$.

 By admissibility there is a central extension $\widehat LT^\tau \to
LT\rtimes \Troth$ for a finite cover $\Troth\to\Trot$ of degree~$\delta $.
Choose a splitting $R_u\mapsto \tilde{R}_u$ for~$R_u\in \Troth,\,0\le u\le
2\pi \delta $.  Then
  $$ \tilde{R}_u\tilde{\varphi }_X\tilde{R}_{-u}\tilde{\varphi }_X\inv
     =\widetilde{\exp(uX)}c_u(X) \tag{2.30} $$
in $\widehat LT$, where $c_u(X)$~is central.  Write $c_u(X) =
\exp\bigl(iu\dot c(X)\bigr)$.  Imposing the group law we find for~$X_1,X_2\in
\Pi $ that
  $$ \frac{c_u(X_1+X_2)}{c_u(X_1)c_u(X_2)} = \kappa(X_2)(\exp uX_1),
      $$
from which $\langle \cdot ,\cdot \rangle$~is symmetric and~$\dot c (X)=
-\frac 12 \langle X,X-\eta _0 \rangle$ for some~$\eta _0\in \frak{t}$.  This
proves~(iii).  Furthermore, from the condition $c_{2\pi \delta}(X)=1$ we
deduce~$\frac \delta 2\langle X,X-\eta _0 \rangle\in \ZZ$ for all~$X\in \Pi
$.  Below we prove that $\eta _0=0$ and so this condition is always satisfied
if~$\delta =2$ or if $\delta =1$ and $\Pi $~is {\it even\/}: $\langle X,X
\rangle\in 2\ZZ$ for all~$X\in \Pi $.

Let $U^\tau \to U$ be the restriction of~$LT^\tau \to LT$ to~$U$.  Choose an
embedding $L\:V\to \Lie (L\frak{t}^\tau )$ which splits the projection
$\Lie(L\frak{t}^\tau )\to\Lie(LT)\to V$, and for $v\in V$ set $\tilde
v=\exp(Lv)$.  For~$X\in \Pi $ define the character $a_X\:V\to\TT$ by
  $$ \tilde{\varphi }_X\tilde{v}\tilde{\varphi }_X\inv \tilde{v}\inv =a_X
     (v). \tag{2.31} $$
Let $\dot a_X\:V\to \RR$ be $\frac 1i$ times the derivative of~$a_X$.  Then
use~\thetag{2.29}, \thetag{2.30}, and~\thetag{2.31} to see that for~$X \in
\Pi $, $\xi \in \frak{t}$, and $v\in V$ we have
  $$ \align
      \Ad_{\tilde{\varphi} \mstrut _X }(d) &= d-X+\frac 12\langle X,X+\eta _0
     \rangle K \tag{2.32}
      \\
      \Ad_{\tilde{\varphi} \mstrut _X }(\xi ) &= \xi + \dot\kappa_X(\xi )K
     \tag{2.33}
      \\
      \Ad_{\tilde{\varphi} \mstrut _X }(v ) &= v + \dot a_X(v)K
     \tag{2.34}
      \endalign $$  
We apply the $LT^\tau $-invariance of~$\ll\!\cdot ,\cdot \!\gg$
and~\thetag{2.12}, \thetag{2.13} to \thetag{2.32}--\thetag{2.34}.  Thus
  $$ \ll\! d,v \!\gg\;=\;\ll\! \Ad_{\tilde{\varphi} \mstrut _X
     }(d),\Ad_{\tilde{\varphi} \mstrut _{X }}(v )\!\gg\;=\; \ll\! d,v
     \!\gg\;-\;\dot a_X(v),  $$
from which $\dot a_X(v)=0$ and so~$a_X(v)=1$.  Note that
$\ll\!X,v\!\gg\;=\;0$ by the argument in the first paragraph of the proof.
Then \thetag{2.31}~implies~(i).  Next,
  $$ \ll\! d,\xi \!\gg\;=\;\ll\! \Ad_{\tilde{\varphi} \mstrut _X
     }(d),\Ad_{\tilde{\varphi} \mstrut _{X }}(\xi )\!\gg\;=\; \ll\! d,\xi
     \!\gg\;-\;\ll\!X ,\xi \!\gg \;-\;\dot\kappa_X(\xi),  $$
from which $\ll\!X ,\xi \!\gg \;=-\dot\kappa_X(\xi)$ and so 
  $$ \ll\!\xi ,\eta \!\gg \;=\;\langle \xi ,\eta \rangle \quad \text{for
     all~$\xi ,\eta \in \frak{t}$}.   $$
From~\thetag{2.32} and \thetag{2.14} we deduce 
  $$ 0=\;\ll\!d,d\!\gg\;=\;\ll\!\Ad_{\tilde{\varphi} \mstrut _X
     }(d),\Ad_{\tilde{\varphi} \mstrut _X }(d)\!\gg = - \langle X,X+\eta _0
     \rangle + \ll\!X,X\!\gg,  $$
and so~$\eta _0=0$ as claimed earlier.
        \enddemo

 \subsubhead Finite groups\endsubsubhead
 For $G$~finite and $P\to \cir$ a principal $G$-bundle, any central extension
of~$\LGR$ is admissible.   To see this, fix a  basepoint~$p\in P$ and suppose
the  holonomy is~$h\in G$.   Then we  identify~$P$ with  $(\RR\times G)/\ZZ$,
where $1\in \ZZ$~acts as
  $$ (x,g)\longmapsto (x+1,hg). \tag{2.35} $$
Then the ``loop group''~$\LGR$, or group of gauge transformations of~$P$, is
identified with $Z_h\subset G$, the centralizer of~$h$.  The
extension~$\LGRR$ (with $\Troth=\Trot$) is identified with $(\RR\times
Z_h)/\ZZ$, where $1\in \ZZ$~acts as in~\thetag{2.35}.  Given a central
extension~$Z_h^\tau $ of~$Z_h$, we fix a lift~$\tilde{h}\in Z_h^\tau $ of~$h$
and define $\LGRRt$ as the quotient $(\RR\times Z_h^\tau )/\ZZ$, with
$\tilde{h}$ replacing~$h$ in~\thetag{2.35}.  (Different choices
of~$\tilde{h}$ lead to isomorphic central extensions.)

 \subhead \S{2.4}.  Finite energy loops
 \endsubhead

First, given~$A\in \AR$ we construct a decomposition of~$\LgR$ into finite
dimensional subspaces adapted to~$d_A$.  Denote by~$\ZA\subset  \LGR$ the
stabilizer of~$A\in \AR$.  Fix a $G$-invariant inner product~$\langle \cdot
,\cdot \rangle$ on~$\frak{g}$.  Using the basepoint~$p\in P$ we identify the
fiber~$\frak{g}\mstrut _P$ at~$0\in \cir$ of the adjoint bundle
with~$\frak{g}$.  The holonomy~$\hol(A)$ acts as an orthogonal transformation
of~$\frak{g}$, so as a unitary transformation of~$\frak{g}_{\CC}$.  Write the
latter as~$\exp(-2\pi i\SA)$ for some self-adjoint~$\SA$ with eigenvalues
strictly between~$-1$ and~$1$.  Notice that $\SA$~is not uniquely defined,
nor can it be made to depend continuously on~$A$.  Decompose
  $$ \frak{g}_\CC \cong (\zA)_{\CC} \oplus \pA \oplus \nA,
     \tag{2.36} $$
where $\frak{z}_A$~is the Lie algebra of~$Z_A$ and $\SA$~is positive
on~$\pA$, negative on~$\nA$.  It is useful to write
  $$ \pA = \bigoplus\limits_{0<\epsilon <1}(\pA)_\epsilon , \tag{2.37} $$
where $S_A=\epsilon $ on~$(\pA)_\epsilon $; there is a similar decomposition
of~$\nA$.  By definition $\SA$~vanishes on the Lie algebra~$\zA$.  Define
  $$ \aligned
      \frak{g}_\CC &\longrightarrow \LgR_{\CC} \\
      \xi &\longmapsto \xA\endaligned \tag{2.38} $$
by letting $\xA(s)$ equal~$\exp(is\SA)$ applied to the parallel transport
of~$\xi $ using~$A$.  Then 
  $$ [d_A,\xA] = i@,\SA(\xi )_A.  $$
For~$n\in \ZZ$ define 
  $$ (z^n\xA)(s) = e^{ins}\xA(s), \tag{2.39} $$
so that 
  $$ [d_A,z^n\xA] = iz^n(S_A + n)(\xi )_A.  $$
The algebraic direct sum 
  $$ \LgRfin(A)=\bigoplus\limits _{n\in \ZZ}z^n\frak{g}_\CC \tag{2.40} $$
is dense in~$\LgR_\CC$.  Elements of~$\LgRfin(A)$ are said to have {\it
finite energy\/}.
 
Let $\LGRt\to\LGR$ be an admissible central extension.  Recall from
\theprotag{2.18} {Lemma} that the connection defines a splitting $\beta
\mapsto\beta ^\tau _A$ of the loop algebra extension. 

        \proclaim{\protag{2.41} {Lemma}}
 If $[d_A,\beta ]=iE\beta $ in~$\LgR$, then the lift of~$\beta $ to the
central extension~$\LgRt$ also has definite energy: $[d_A^\tau ,\beta
^\tau _A]= iE\beta ^\tau _A$.
        \endproclaim

        \demo{Proof}
Use the characterization~\thetag{2.19} and the invariance of the
bilinear form: 
  $$ \ip{d_A^\tau }{[d_A^\tau ,\beta ^\tau_A]}_\tau =
      \ip{[d_A^\tau ,d_A^\tau  ]}{\beta ^\tau_A}_\tau = 0.
      $$
Since $[d_A^\tau , \beta ^\tau_A]$~is a lift of~$iE\beta $ the result
follows.
        \enddemo

 \subhead \S{2.5}.  Positive energy representations
 \endsubhead

Following Segal~\cite{PS,\S9} we define a distinguished class of projective
representations of loop groups.  

        \definition{\protag{2.42} {Definition}}
 Fix $\LGRt$~ an admissible graded central extension of~$\LGR$.  Let $\rho
\:\LGRt\to U(V)$~be a unitary representation on a $\zt$-graded complex
Hilbert space.  Assume that $\rho $~is even and the center~$\TT$ acts by
scalar multiplication.  We say $V$~has {\it positive energy\/} if
 \roster
 \item $\rho $~extends to a unitary representation $\rt\:\LGRRt\to U(V)$.
For~$A\in \AR$ denote the skew-adjoint operator~$\dot{\rt}(d_A^\tau )$ by~$iE_A$.
 \item For all~$A\in \AR$ the {\it energy operator\/} $E_A$~is self-adjoint
with discrete spectrum bounded below.
 \endroster A positive energy representation is said to be {\it finitely
reducible\/} if it is a finite sum of irreducible representations.
        \enddefinition

\flushpar
 We regard ungraded representations as $\zt$-graded representations with odd
part zero.  For an irreducible positive energy representation the
extension~$\rt$ is unique up to a character of~$\Troth$, so the energy is
determined up to a shift by~$n/\delta $, where $n\in \ZZ$ and $\delta $~is
the degree of the cover $\Troth\to\Trot$.
 
Let $V$~be a positive energy representation of~$\LGRt$, and choose an
extension to $\LGRRt$.  Fix a connection~$A$ and decompose~$V$ as the Hilbert
space direct sum of closed $\zt$-graded subspaces~$V_e(A),\, e\in \RR$ on
which $E_A$~acts as multiplication by~$e$.  The positive energy condition
asserts that~$V_e(A)\not= 0$ for a discrete set of~$e$ which is bounded
below.  The algebraic direct sum
  $$ \Vfin(A) = \bigoplus\limits_{e\in \RR}V_e(A) \tag{2.43} $$
consists of vectors of {\it finite energy\/}; it is dense in~$V$.  Now $A$~
defines a splitting of the Lie algebra central extension (see
\theprotag{2.18} {Lemma}), and we use it to define the infinitesimal
representation $\rdA$ on elements of~$\LgR$; namely, $\rdA(\beta )$~denotes
the infinitesimal action of~$\beta^\tau _A$, where~ $\beta \in \LgR$.  The
image of~$\rdA$ consists of skew-Hermitian operators on~$V$, which in general
are unbounded.  Note from \theprotag{2.41} {Lemma} that the lift of $\beta
\in \LgR$ with definite energy has the same energy.  Suppose~$\xi \in
\frak{g}_\CC$ is an eigenvector of~$\SA$ of eigenvalue~$\epsilon $.  It
follows that
  $$ \rdA(z^n\xi _A)\bigl(V_e(A) \bigr) \subset V_{e+n+\epsilon }(A),
      $$
so in particular $\rdA(z^n\xi _A)\bigl(V_e(A) \bigr) = 0$ if $e+n+\epsilon
$~is less than the minimal energy~$\emin$.  A positive energy
representation~$V$ has a minimal $E_A$-energy subspace~$V_{\emin}(A)$ which
is a representation of~$\ZAt$, the restriction of the central
extension~$\LGRt$ over the stabilizer~$\ZA$ of~$A\in \AR$.  Also, $\nA \oplus
\bigoplus\limits_{n>0}z^{-n}\frak{g}_{\CC}$ acts trivially on~$V_{\emin}(A)$.
If $V$~is irreducible, then $\Vfin(A)$~is the subspace spanned by vectors 
  $$ \rtd_{A}\bigl(z^{n_r}(\xi _r)_A\bigr)\dots \rtd_{A}\bigl(z^{n_1}(\xi
     _1)_A\bigr)\Omega, \tag{2.44} $$
where $\Omega \in V_{\emin}(A)$, $r\ge0$, $\xi _i$~is an eigenvector of~$S_A$
with eigenvalue~$\epsilon _i$, and~$n_i\ge 0$ with equality only if~$\epsilon
_i>0$.  The energy of~\thetag{2.44} is $\emin + \sum (n_i+\epsilon _i)$.
This proves the following.

        \proclaim{\protag{2.45} {Lemma}}
 Let $V$~be a finitely reducible positive energy representation.  Then for
any connection~$A$ and~$C>0$ the sum of the $E_A$-eigenspaces for eigenvalue
less than~$C$ is finite dimensional. 
        \endproclaim

\flushpar
 In other words, $E_A$~has the spectral characteristics of a first-order
elliptic operator on the circle.

We use~\thetag{2.14} to compute the dependence of the energy operator on the
connection.  Let $A'$~denote the connection whose horizontal vector field
is\footnote{The corresponding 1-form is~$A+\beta \,ds$;
see~\thetag{2.17}.}~$d_A-\beta $ for some $\beta \in \LgR$.  Then
  $$ E_{A'} = E_A + i\rdA(\beta ) + \frac{\ll\! \beta ,\beta
     \!\gg}{2}. \tag{2.46} $$

The following proposition relates the existence of positive energy
representations to the bilinear form on the central extension (\theprotag{2.10}
{Definition}). 

        \proclaim{\protag{2.47} {Proposition}}
 Let $\LGRt$~be an admissible central extension of~$\LGR$ which admits a
nonzero positive energy representation~$V$.  Then the inner
product~$\form_\tau $ on~$\LgR$ is positive semidefinite.
        \endproclaim

        \demo{Proof}
 Fix a connection~$A$.  It suffices to prove that 
  $$ \ll \zx,\zxc\gg\;\;\ge\;\;0\qquad \text{for all $\xi \in (\frak{z}_A)_\CC
     \oplus \frak{p}_A,\;n\ge0$}. \tag{2.48} $$
(Recall the notation from~\thetag{2.36}, \thetag{2.38}, and~\thetag{2.39};
also, we drop the subscript~`$\tau $' on the form, which is bilinear on the
complexified loop algebra.)  By \theprotag{2.18} {Lemma} $A$~defines a
linear splitting $\LgR\to\LgRt$, which is characterized by~\thetag{2.19}; we
denote it~$\beta \mapsto \beta ^\tau $.  Then
  $$ [\beta_1 ^\tau ,\beta_2 ^\tau ] = [\beta_1 ,\beta_2 ]^\tau + \omega
     (\beta_1 ,\beta_2 )\frac Ki  $$
for a real skew form~$\omega $, and 
  $$ \ll d_A^\tau ,[\beta_1 ^\tau ,\beta_2 ^\tau ]\gg \;\;=\;\;-\,\frac{\omega
     (\beta_1 ,\beta_2 )}{i}.  $$
For $\xi ,n$ as in~\thetag{2.48} with $\xi $~an eigenvector of~$\SA$ with
eigenvalue~$\epsilon >0$ (see~\thetag{2.37}), we find 
  $$ \split  
      \ll d_A^\tau ,\left[ \zxt,\zxct \right]\gg \;\;&=\;\; \ll[d_A^\tau ,\zxt],\zxct\gg\\
     &=\;\; i(n+\epsilon )\;\ll\zx,\zxc \gg, \endsplit $$
whence
  $$ \omega (\zx,\zxc) = (n+\epsilon )\;\ll\zx,\zxc\gg.  $$
Then for~$\Omega\in V_{\emin}(A)$ and $n>0$, since $\rdA\bigl(\zxt
\bigr)^*=-\rdA\bigl(\zxct \bigr)$ we have
  $$ \split
      0\le \|\,\rdA\bigl(\zxt\bigr)\Omega \,\|_V^2 &= -\,\langle\;
     \rdA\bigl([\zxct,\zxt]\bigr)\Omega \,,\,\Omega \;\rangle_V\\
      &= (n+\epsilon )\,\ll \zx,\zxc\gg\;\;-\;\;\langle\; \rdA\bigl([\bar\xi
     _A,\xA]^\tau \bigr)\Omega \,,\,\Omega \;\rangle_V. \endsplit\tag{2.49} $$

Now we use the orthogonal decomposition~\thetag{2.5}.  If $\xi \in
\cent(\frak{g})$, then the last term of~\thetag{2.49} vanishes and
\thetag{2.48}~follows immediately.  On the other hand, if $\xi =[\eta ,\eta
']$, then $\xA = [\eta _A,\eta _A']$ and
  $$ \split
      \ll\zx\,,\,\zxc\gg\;\;&=\;\;\ll[\eta _A,z^n\eta _A']\,,\,\zxc\gg \\
      &=\;\;\ll\eta _A\,,\,[z^n\eta '_A,\zxc]\gg \\
      &=\;\;\ll\eta _A\,,\,[\eta '_A,\bar\xi _A]\gg \\
      &=\;\;\ll[\eta _A,\eta '_A]\,,\,\bar\xi _A\gg \\
      &=\;\;\ll\xA\,,\,\bar\xi _A\gg.  \endsplit  $$
The argument applies to sums of terms~$[\eta ,\eta ']$.  If
$\ll\xA,\bar\xi _A\gg\;<0$, then for $n$~large the right hand side
of~\thetag{2.49} is negative, which is a contradiction.  This completes the
proof of~\thetag{2.48}.
        \enddemo

We state without proof some basic facts about positive energy
representations.  To avoid problems with abelian factors, we assume that the
inner product~$\form$ is positive definite when restricted to~$L_P(\frak{z})$
for~$\frak{z}=\operatorname{center}(\frak{g})$.
  $$ \align
      &\text{Positive energy representations are completely
     reducible~\cite{PS,\S11.2}}  \\
      &\text{For any irreducible positive energy representation,}
     \\ \vspace{-6pt}
      &\text{$V_e(A)$~is finite dimensional~\cite{PS,(9.3.4)}}\\
      &\text{If $V$~is finitely reducible and $\rt$ in \theprotag{2.42(1)}
     {Definition} is given, then if}\\ \vspace{-6pt}
      &\text{\theprotag{2.42(2)} {Definition} is satisfied for one~$A\in\AR$
     it is satisfied for all~$A\in \AR$} \\
      &\text{There is a finite number of isomorphism classes of irreducible
     positive}\\ \vspace{-6pt}
      &\text{energy representations of~$\LGRt$}\\ \endalign $$
See~\cite{PS,(9.3.5)} and the remarks which follow for the last statement.

Following the usual procedure we make an abelian group out of positive energy
representations.

        \definition{\protag{2.50} {Definition}}
 Let $\Rt$~denote the abelian group generated by isomorphism classes of
finitely reducible $\zt$-graded positive energy representations of~$\LGRt$
under direct sum, modulo the subgroup of isomorphism classes of
representations which admit a commuting action of the Clifford
algebra~$\cco$.
        \enddefinition

\flushpar
 We extend~$\Rt$ to a $\zt$-graded free abelian group.  The aforementioned
representations are even.  Odd representations in addition carry a (graded)
commuting action of~$\cco$.  An odd representation is equivalent to zero if
this action extends to a commuting action of $C_2^c$.

        \example{\protag{2.51} {Example}}
 To illustrate odd representations consider the cyclic group~$\zt$ of order
two with generator~$x$.  With the trivial grading and no central extension
there are two inequivalent irreducible even representations (of graded
dimension~$1|0$) and no nontrivial odd representations.  Thus $K^0_{\zt}$~has
rank two and $K^1_{\zt}=0$.  With the nontrivial grading~$\epsilon $ and no
central extension there is up to isomorphism one even irreducible
representation on the $(1|1)$-dimensional space~$\CC\oplus \CC$; the
element~$x$ acts as~$\left(\smallmatrix 0&1\\1&0 \endsmallmatrix\right)$.
But this represents the trivial element of $K$-theory, as there is a
commuting action of~$\cco$ where the generator acts as~$\left(\smallmatrix
0&i\\i&0 \endsmallmatrix\right)$.  Together with this $\cco$-action this
representation is nontrivial in odd $K$-theory, as there is no extension to
an action of~$C_2^c$.  Therefore, $K^\epsilon _{\zt}=0$ and
$K^{\epsilon+1}_{\zt}$~has rank one.
        \endexample

 \head
 \S{3} Dirac Families and Loop Groups 
 \endhead
 \comment
 lasteqno 3@ 44
 \endcomment

 \subhead \S{3.1}.  The spin representation (infinite dimensional case)
 \endsubhead

For an algebraic approach, see~\cite{KS}; our geometric approach
follows~\cite{PS,\S12}.  Let $H$~be an infinite dimensional real Hilbert
space.  We define the Clifford algebra~$\Cliff^c(H^*)$ as in finite
dimensions~\thetag{1.3}, but now the construction of an irreducible Clifford
module, and so of the spin representation, depends on a {\it polarization\/}.
Recall that a complex structure on~$H$ is an orthogonal map $J\:H\to H$ with
$J^2=-1$; it follows that $J$~is skew-symmetric.

        \definition{\protag{3.1} {Definition}}
   A {\it polarization\/} is a set~$\scrJ$ of complex structures maximal
under the property that any two elements differ by a Hilbert-Schmidt
operator.  The {\it restricted orthogonal group\/} of~$(H,\scrJ)$ is
  $$ \OJ(H) = \{T\in \Or(H): TJT\inv \in \scrJ \text{ for all $J\in
     \scrJ $}\}.  $$
        \enddefinition

\flushpar
 The group $\OJ(H)$~has two components, each simply connected with second
homology free of rank one; in fact, $\OJ(H)$~has the homotopy type
of~$\Or_{\infty }/\Ur_{\infty }$.  In particular, there is a unique nontrivial
grading $\epsilon \:\OJ(H)\to\zt$.  Also, $\OJ(H)$~acts transitively
on~$\scrJ$ with contractible isotropy groups, so $\scrJ$~has the same
homotopy type as~$\OJ(H)$.  There is a distinguished central extension
  $$ 1 \longrightarrow \TT\longrightarrow \Pin^c_{\scrJ}(H)\longrightarrow
     \Or_{\scrJ}(H)\longrightarrow 1, \tag{3.2} $$
which, together with the grading, forms a graded central extension
of~$\OJ(H)$.  We remark that there is no infinite dimensional analog of the
Pin group~\thetag{1.1}.

If $J_0\:H\to H$ is a real skew-adjoint Fredholm operator with nonzero
eigenvalues~$\pm i$ of infinite multiplicity---in other words, $J_0$~has
finite dimensional kernel and is a complex structure on the orthogonal
complement---then it determines a polarization 
  $$ \scrJ =\bigl\{J: J\text{ a complex structure},\;J-J_0\text{
     Hilbert-Schmidt}\bigr\} \tag{3.3} $$
if $\dim\ker J_0$ is even.  (Recall that the parity of $\dim\ker J_0$---the
mod~2 index of~$J_0$---is invariant under deformations.)  If $\dim\ker
J_0$~is odd, then define 
  $$ \multline \scrJ = \bigl\{J: J\text{ real skew-adjoint
     Fredholm},\;\dim\ker J=1,\\ \text{$J$ a complex structure on $(\ker
     J)^\perp$},\;J-J_0 \text{
     Hilbert-Schmidt}\bigr\}.\endmultline\tag{3.4} $$
(We require that the eigenvalues~$\pm i$ have infinite multiplicity.)  Call
\thetag{3.4}~an {\it odd\/} polarization; \thetag{3.3}~ is {\it even\/}.
The kernels lead to an action of~$\cco$ on~$\SS$ which commutes with the
$\Pin^c_{\scrJ}(H)$ and $\Cliff^c(H^*)$-actions.  As usual, we carry this extra
structure implicitly.

The spin representation is a distinguished irreducible unitary representation
  $$ \chi \:\Pin^c_{\scrJ}(H)\to\Aut(\SS), \tag{3.5} $$
where $\SS=\SS^0\oplus \SS^1$ is a $\zt$-graded complex Hilbert space.  It is
constructed in~\cite{PS,\S{12}} as the space of sections of a holomorphic
line bundle over~$\scrJ$.  There is a Clifford multiplication
  $$ \gamma \:H^* \longrightarrow \End(\SS)  $$
which is compatible with~\thetag{3.5} in the sense that
  $$ \gamma (g\cdot \mu )=\chi (g)\gamma (\mu )\chi (g)\inv ,\qquad g\in
     \Pin^c_{\scrJ}(H),\quad \mu \in H^*. \tag{3.6} $$
We arrange that $\gamma (\mu )$~ be skew-Hermitian.  Clifford
multiplication~$\gamma $ is odd, and $\chi$~is compatible with the gradings
on~$\Pin^c_{\scrJ}(H)$ and~$\SS$.
 
Let $P\to\cir$ be a principal $G$-bundle.  Fix\footnote{Different choices
lead to isomorphic constructions.  In fact, we can use any invariant inner
product~$\form$ on~$\LgR$ with respect to which $d_A$~is skew-Hermitian.} a
$G$-invariant inner product~$\langle \cdot ,\cdot \rangle$ on~$\frak{g}$, and
so an $L^2$~inner product~$\form$ on~$\LgR$ (as in ~\thetag{2.25}).  Take
$H$~be the $L^2$~completion of~$\LgR$.  For~$A\in \AR$ the elliptic
operator~$d_A$ extends to a skew-adjoint operator on~$H$ with discrete
spectrum.  There is an orthogonal decomposition of the
complexification\footnote{Of course, \thetag{3.7}~is the $L^2$~completion
of~\thetag{2.40}.  Below we denote the real points in the algebraic direct
sum as~$H_{\operatorname{fin}}$.  The finite dimensional subspace~$H_e(A)$
contains smooth elements, i.e., lies in~$\LgR_{\CC}$.}
  $$ H_{\CC} \cong \bigoplus\limits_{e\in \RR}H_e(A)\qquad \text{($L^2$
     completed)}, \tag{3.7} $$
where $d_A=ie$ on~$H_e(A)$.  Elements of~$H_e(A)$ are said to have energy~$e$.
Define the {\it real\/} skew-adjoint Fredholm operator
  $$ J_A = \cases 0 ,&\text{on $H_0(A)$};\\\dfrac {d_A}{|e|},&\text{on
     $H_e(A)$}.\endcases \tag{3.8} $$
The proofs of the next two lemmas are similar to~\cite{PS,(6.3.1)}, so are
omitted.  The loop group~$\LGR$ acts on the loop algebra~$\LgR$ by the
pointwise adjoint action of~$G$ on~$\frak{g}$.

        \proclaim{\protag{3.9} {Lemma}}
 The polarization~$\scrJ$ defined by~$J_A$ is independent of~$A$.
        \endproclaim

        \proclaim{\protag{3.10} {Lemma}}
 Any~$\varphi \in \LGR$ extends to a bounded orthogonal operator on~$H$ which
preserves the polarization~$\scrJ$. 
        \endproclaim

\flushpar
 Hence we have a homomorphism 
  $$ \LGR\longrightarrow \OJ(H), \tag{3.11} $$
and so by pullback of~\thetag{3.2} a distinguished graded central extension 
  $$ 1 \longrightarrow \TT \longrightarrow \LGRs \longrightarrow \LGR
     \longrightarrow 1 \tag{3.12} $$
and a graded unitary representation of~$\LGRs$ on~$\SS$. 

        \proclaim{\protag{3.13} {Proposition}}
 The central extension~$\LGRs$ is admissible with respect to the symmetric
positive semidefinite bilinear form~\thetag{2.16}. 
        \endproclaim

\flushpar 
 Therefore, by \theprotag{2.18} {Lemma} a connection~$A\in \AR$ defines a
splitting 
  $$ \aligned
      \LgR &\longrightarrow \LgRs \\
      \beta &\longmapsto \beta _A^\sigma \endaligned  $$
of the spin central extension of the loop algebra. 

        \demo{Proof}
 The adjoint action factors through the semisimple adjoint group~$\Ad G$,
which leads to a factorization of~\thetag{3.11} through a twisted loop group
for~$\Ad G$.  Thus we may as well assume that $G$~is semisimple.  Fix a
connection~$A$ on~$P\to\cir$ and let $\varphi _t$~be the one-parameter group
in~$\LGRR$ generated by~$d_A$.  (We take $\Troth=\Trot$ in~\thetag{2.2}.)
Then $\LGRR$~is generated by~$\LGR$ and~$\{\varphi _t\}_{t\in \RR}$.  Now
by~\thetag{3.8} we have $\varphi _t^*J_A=J_A$, and so $\varphi _t$~acts
on~$\LgR$ by an element of~$\OJ(H)$.  Therefore~\thetag{3.11}
and~\thetag{3.12} extend to~$\LGRR$.  By \theprotag{2.15} {Proposition} there
is a unique bilinear form $\form_\sigma $ on~$(\LgRR)^\sigma $ compatible
with the lifted rotation action.
 
We use an explicit model of~$\SS$ to compute~$\form_\sigma $.  Dual
to~\thetag{3.7} is a decomposition
  $$ H_{\CC}^* \cong \bigoplus\limits_{e\in \RR}H_e^{*}(A)\qquad \text{($L^2$
     completed)},\qquad H_e^*(A) = \bigl(H_{-e}(A) \bigr)^*, \tag{3.14} $$
into finite dimensional subspaces of energy~$e$.  Then the dense finite
energy subspace of the spin representation may be realized as
(compare~\thetag{1.16})
  $$ \Sfin(A) = \SS_0(A)\;\otimes\; {\tsize\bigwedge} ^{\bullet}
     \Bigl(\bigoplus\limits_{e>0}H_e^*(A) \Bigr),  \tag{3.15} $$
where $\SS_0(A)$~is an irreducible $\zt$-graded module for the
(finite-dimensional) Clifford algebra of~$\frak{z}_A^*$.  (Recall that
$\frak{z}_A$~is the Lie algebra of the stabilizer~$Z_A$ of~$A$.)  Normalize
the energy on~$\SS$ so elements of~$\SS_0(A)\otimes {\tsize\bigwedge} ^0$
have zero energy.  Now the inner product~$\form$ on~$H$ induces a linear map
$H_{\CC}\to H^*_{\CC}$, denoted~$\beta \mapsto\beta ^*$, which sends elements
of energy~$e$ to elements of energy~$-e$.  Then for $\beta \in H_{-e}(A)$ we
have $\beta ^*\in H_{e}^*(A)$ and Clifford multiplication on~$\Sfin(A)$ is
given as
  $$ \gamma (\beta ^* ) = \cases 1\otimes \epsilon (\beta ^* )
     ,&e>0;\\\gamma _0(\beta ^* )\otimes 1,&e=0;\\-1\otimes \iota (\beta
     ),&e<0.\endcases \tag{3.16} $$
Here $\epsilon $~is exterior multiplication, $\iota $~is interior
multiplication, and $\gamma _0$~is Clifford multiplication on~$\SS_0(A)$.
This satisfies the usual Clifford relation
  $$ \bigl[ \gamma (\beta_1 ^*),\gamma (\beta_2 ^*)\bigr] = -2\ll\! \beta_1
     ,\beta_2 \!\gg = -2\ll\! \beta_1^* ,\beta_2^* \!\gg \tag{3.17} $$
for $\beta_1 ,\beta_2 $ of finite energy.  On~$\Sfin(A)$ the
operator~$\cd\bigl(\lifts{\beta}\bigr)$ corresponding to~$\beta $ of finite
energy is given by the first formula of~\thetag{1.7}, where now the sum is
infinite as the indices range over a basis of finite energy vectors
in~$\LgR$; see~\thetag{3.19} below.  But only a finite number of terms are
nonzero when acting on a fixed element of~$\Sfin(A)$.  (We formalize such
infinite sums in~\S{3.2}.)

For the explicit computation we first treat the case $P\to\cir$ trivial and
choose $A=A_0$ the trivial connection.  Fix a basis~$e_a$
of~$\frak{g}_{\CC}$, which we identify with constant loops
in~$L\frak{g}_{\CC}$.  Then $\{z^ne_a\}_{n\in \ZZ}$~is an (algebraic) basis
of~$L\frak{g}_{\CC}$; let $\{z^{-n}e^a\}_{n\in \ZZ}$~be the dual basis.  For
any~$n\in \ZZ$ let $\gamma ^a(n)$~denote Clifford multiplication
by~$z^{n}e^a$; it has energy~$n$.  Now let $\xi ,\eta \in \frak{g}_\CC$ be
constant loops, and as usual $\liftso{z\xi },\liftso{z\inv \eta }$ the lifts
of the indicated loops to the central extension~$(L\frak{g})^\sigma _{\CC}$.
Then from the invariance of the bilinear form we see that it suffices to
compute $\ip{z^n\xi }{z^m\eta }_\sigma $ for~$n=-m=1$ and also that
  $$ [\liftso{z\xi },\liftso{z\inv \eta }] = [z\xi ,z\inv \eta ]^\sigma
     _{A_0} \;\;+ \; \ip{\xi}{\eta }_\sigma \frac Ki. \tag{3.18} $$
(The proof of \theprotag{2.47} {Proposition} contains these assertions.)  We
evaluate the commutator using the infinitesimal spin representation, and it
suffices to evaluate on a (vacuum) vector~$\Omega \in \SS_0(A_0)\otimes
{\tsize\bigwedge} ^0$.  Now the infinitesimal spin action is
(c.f.~\thetag{1.7}) 
  $$ \cd\bigl(\liftso{z^ne_a} \bigr) = \frac 14f_{abc}\sum\limits_{k+\ell
     =n}\gamma ^b(k)\gamma ^c(\ell ), \tag{3.19} $$
and as stated above the sum is finite on finite energy vectors.  For example,
  $$ \cd\bigl(\liftso{z\xi} \bigr)\Omega = \frac 14f_{abc}\xi ^a[\gamma
     ^b(1)\,\gamma ^c(0) + \gamma ^b(0)\,\gamma ^c(1)]\,\Omega = \frac
     14f_{abc}\xi ^a[2\gamma ^b(1)\,\gamma ^c(0)]\Omega.  $$
Hence
  $$ \split
      \bigl[ \cd(\liftso{z\xi }),\cd(\liftso{z\inv \eta })\bigr]\Omega
       &= - \cd\bigl(\liftso{z\inv \eta } \bigr)\cd\bigl(\liftso{z\xi }
     \bigr)\Omega \\
      &= -\frac 14f_{abc}f_{a'b'c'}\xi ^a\eta ^{a'}\gn{b'}0\,\gn{c'}{-1}\,\gn
     b1\,\gn c0\,\Omega \\
      &= \frac 12f_{abc}f_{a'b'c'}g^{bc'}\xi ^a\eta ^{a'}\gn{b'}0\,\gn
     c0\,\Omega \\
      &=\left[ \frac 14 f^b_{aa'}f\mstrut _{bb'c}\gn{b'}0\,\gn c0 - \frac
     12f^c_{ab}f^b_{a'c} \right]\xi ^a\eta ^{a'}\Omega \\
      &= \left[ \cd\bigl([\xi ,\eta ]^\sigma _{A_0} \bigr) - \frac
     12\Tr(\ad\xi \circ \ad\eta ) \right]\Omega . \endsplit \tag{3.20} $$
To pass from the third to fourth line we use the Jacobi identity and
rearrange the terms.  Recalling that $K$~acts as multiplication by~$i$ in a
representation, we see from ~\thetag{3.18} that
  $$ \ip\xi \eta _\sigma = -\,\frac 12\Tr(\ad\xi \circ \ad\eta ), \tag{3.21} $$
which agrees with~\thetag{2.16}.

For the general (twisted) loop group the computation is similar, but now we
use a generic connection~$A$ which satisfies~\thetag{A.4}, \thetag{A.5} in
the appendix.  (We freely use the notation of that appendix in this
paragraph.)  Then from~\thetag{A.8} we have the basis $\{z^n\xi
_i,z^n\cj,z^n\pa,z^n\pab\}\mstrut _{n\in \ZZ}$ of~$\LgR_{\CC}$, but now the
energies are not necessarily integers.  Write this basis as~$\{e_p\}$, let the
structure constants be $F\mstrut _{pqr}=\ip{[e_p,e_q]}{e_r}$, and let $\gamma
^p$~denote the Clifford multiplication.  Then the usual formula
$\cd\bigl(\lifts{e_p} \bigr)= \frac 14 F\mstrut _{pqr}\gamma ^q\gamma ^r$
holds for basis elements~$e_p$ of {\it nonzero\/}\footnote{For the trivial
connection \thetag{3.19}~{\it is\/} valid for zero energy~($n=0$).  This may
be verified from the $n\not= 0$~formula by computing
  $$ \bigl[\cd\bigl(\liftso{z^2\xi }\bigr),\cd\bigl(\liftso{z^{-2}\eta
     }\bigr)\bigr] - \bigl[\cd\bigl(\liftso{z\xi
     }\bigr),\cd\bigl(\liftso{z^{-1}\eta }\bigr)\bigr],  $$
analogous to \thetag{3.22} and~\thetag{3.23} below.} energy, so excludes
the~$\xi _i$.  As before, it suffices to compute~$\ip\xi \xi _\sigma
,\,\ip{\cj}{z\inv \chi \mstrut _{j'}}_\sigma $, and $\ip{\pa}{\pab}_\sigma $,
where $\xi \in \frak{h}_0$ and~$E(\cj) + E(\chi \mstrut _{j'})=1$.
Since~$E(\xi )=0$ and~$[\xi ,\xi ]=0$, invariance of the inner product
implies 
  $$ \bigl[\cd\bigl(\lifts{z\xi }\bigr),\cd\bigl(\lifts{z\inv \xi
     }\bigr)\bigr] = i\ip\xi \xi _\sigma K.  $$
The computation of the bracket is similar to~\thetag{3.20}.  For the
remaining cases~$\ip{\cj}{z\inv \chi \mstrut _{j'}}_\sigma $
and~$\ip{\pa}{\pab}_\sigma $, which involve vectors of non-integral energy, we
compute the left hand sides of 
  $$ \bigl[\cd\bigl(\lifts{z\cj}\bigr),\cd\bigl(\lifts{z^{-2}\chi \mstrut
     _{j'}}\bigr)\bigr] \;-\;
     \bigl[\cd\bigl(\lifts{\cj}\bigr),\cd\bigl(\lifts{z^{-1}\chi \mstrut
     _{j'}}\bigr)\bigr] = i\ip{\cj}{z\inv \chi \mstrut _{j'}}_\sigma K
     \tag{3.22} $$
and
  $$ \bigl[\cd\bigl(\lifts{z\pa}\bigr),\cd\bigl(\lifts{z^{-1}\pab
     }\bigr)\bigr] \;-\;
     \bigl[\cd\bigl(\lifts{\pa}\bigr),\cd\bigl(\lifts{\pab }\bigr)\bigr] =
     i\ip{\pa}{\pab}_\sigma K \tag{3.23} $$
In all cases we find~\thetag{3.21}, but we omit the details as the
manipulations are similar to~\thetag{3.20}.
        \enddemo

We record some specific facts about the spin representation.  Recall that
$Z_A\subset \LGR$ is the stabilizer of a connection~$A\in \AR$ and
$\frak{z}_A$ its Lie algebra.

        \proclaim{\protag{3.24} {Lemma}}
	\roster
 \item The spin representation~$\SS$ of~$\LGRs$ satisfies the positive energy
condition.  For any fixed $A\in \AR$ we arrange that the minimal energy be
zero and the minimal energy space~$\SS_0(A)$ be an irreducible $\zt$-graded
module for~$\Cliff^c(\frak{z}_A^*)$.  The dense subspace~$\Sfin(A)$ of finite
energy vectors is generated by $\gamma (\beta _k^*)\cdots\gamma (\beta
_2^*)\gamma (\beta _1^*)s_0$, for $\beta _i\in \LgRfin(A)$ and~$s_0\in
\SS_0(A)$.  Also, the restriction $Z_A^\sigma \to Z\mstrut _A$ of
$\LGRs\to\LGR$ acts on~$\SS_0(A)$ compatibly with the
$\Cliff^c(\frak{z}_A^*)$-action, and $Z_A^\sigma \to Z\mstrut _A$ is isomorphic
to the graded central extension~\thetag{1.4} \rom(with~$G=Z_A$\rom).  The
splitting $\frak{z}\mstrut _A\to\frak{z}_A^\sigma $ of the corresponding Lie
algebra extension induced from~\theprotag{2.18} {Lemma} is a homomorphism of
Lie algebras.
 \item Suppose $P\to\cir$ is trivial and $A_0$~is the trivial connection.
Then $Z_A=G$ and the splitting $\frak{g}\to\frak{g}^\sigma $ from
\theprotag{2.18} {Lemma} agrees with the splitting in~\thetag{1.5}.
	\endroster
        \endproclaim

        \demo{Proof}
 Most of the statements in~\therosteritem{1} are immediate
from~\thetag{3.15}.  For the penultimate assertion, note that $Z_A^\sigma
$~preserves energy, so acts on~$\SS_0(A)$, and the isomorphism class of the
graded central extension is determined by the action.  As stated
after~\thetag{3.15}, $Z_A^\sigma $~acts via the finite dimensional spin
representation, so the extension is isomorphic to~\thetag{1.4}.  For the last
assertion we use the characterization~\thetag{2.19} of the splitting and note
that if $\lifts{\zeta _1},\lifts{\zeta _2}$~are lifts in~$\frak{z}_A^\sigma
$, then
  $$ \ip{d_A^\sigma }{[\lifts{\zeta _1},\lifts{\zeta _2}]}_{\sigma } =
     \;\ip{[d_A^\sigma ,\lifts{\zeta _1}]}{\lifts{\zeta _2}}_{\sigma } = 0,
      $$
so that $[\zeta _1,\zeta _2]^\sigma _A = [\lifts{\zeta _1},\lifts{\zeta
_2}]$. Assertion~\therosteritem{2} is the statement that \thetag{3.19}~holds
for~$n=0$; cf. the first equation of~\thetag{1.7}.
        \enddemo

 \subhead \S{3.2}.  The canonical 3-form on~$\LGR$
 \endsubhead

Consider the 3-form
  $$ \Omega (\beta _1,\beta _2,\beta _3) = \;\ll\![\beta _1,\beta
     _2],\beta _3\!\gg_\sigma ,\qquad \beta _1,\beta _2,\beta _3\in
     \LgR. \tag{3.25} $$
Since $\form_\sigma $, defined in~\thetag{2.16}, is symmetric and
$\LGR$-invariant, the form~$\Omega $ is indeed skew-symmetric.  Furthermore,
$\Omega $~is $d_A$-invariant for all~$A$, since $\form_\sigma $~is
$(\LGRR)^\sigma $-invariant.  Thus if $\beta _1,\beta _2 ,\beta _3$~have
finite energy---i.e., live in the algebraic direct sum~\thetag{2.40}---and
are eigenvectors of~$d_A$, then $\Omega (\beta _1 ,\beta _2,\beta _3)$~is
nonzero only if the sum of the eigenvalues vanishes.  Finally, $\Omega $~is
$\LGR$-invariant as well.  In infinite dimensions there is no natural notion
of Clifford multiplication by a 3-form.  Rather, to define the action
of~$\Omega $ on the spin representation~$\SS$ we use the energy operator of a
connection~$A$ and work with finite energy spinor fields to define an
operator~$Q_A$ on finite energy spinors.  (See~\cite{L,\S7} for another
discussion.)

The subspace of finite energy vectors~$\HfinsA\subset H^*$ has
complexification the {\it algebraic\/} direct sum $\HfinsA_\CC = \bigoplus
\limits_ e H^*_e(A)$.  Clifford multiplication on finite energy
1-forms~\thetag{3.16} extends to an injective map
  $$ \Cliff^c(\HfinsA)\longrightarrow \ESfin. \tag{3.26} $$
Following~\cite{KS,\S7} we endow $\ESfin$~with the weak topology relative to
the discrete topology on~$\Sfin(A)$: a sequence $\{T_N\}\mstrut _N\subset
\ESfin$ converges if and only if $T_N(s)$ stabilizes for all~$s\in \Sfin(A)$
as~$N\to\infty $.  Furthermore, the standard filtration on~$\Cliff^c(\HfinsA)$
induces a filtration on the image of~\thetag{3.26}, and we let $\Efil
p\subset \ESfin$ denote the weak closure of the image of
$\Cliff^c(\HfinsA)_{\le p}$.  The $\zt$-grading on~$\ESfin$ induces a
$\zt$-grading $\Efil p^{\bullet} = \Efil p^0 \oplus \Efil p^1$.  Let
$\End_e\bigl(\Sfin(A) \bigr)$ denote the endomorphisms of energy~$e$, i.e.,
those which raise energy on~$\Sfin$ by~$e$, and set $E_e(A)_{\le p}\bul=
\Efil p\bul\cap \End_e\bigl(\Sfin(A) \bigr)$.
 
We now define an operator~$Q_A$ on~$\Sfin(A)$ which plays the role of
Clifford multiplication by~$\Omega $.  Approximate~\thetag{3.25} by
truncation: let $\ONA \in {\tsize\bigwedge} ^3\HfinsA$ be given
by~\thetag{3.25} on elements~$\beta_1 ,\beta_2 ,\beta_3 $ of definite
energies not exceeding~$N$ in absolute value; $\ONA$~vanishes if some
absolute value of energy exceeds~$N$.  As in finite dimensions there is a
vector space isomorphism\footnote{In fact, an isomorphism of left
$\Cliff^c\bigl(\HfinsA \bigr)$-modules~\cite{KS,\S4}.}  ${\tsize\bigwedge}
^\bullet\bigl(\HfinsA\bigr)\cong \Cliff^c\bigl(\HfinsA \bigr)$, and the image
of 3-forms lies in the odd Clifford algebra at filtration level~$\le3$.
Noting that $\ONA$~has energy zero, we see that its image~$(Q_A)_N$ under
\thetag{3.26} lies in~$E_0(A)_{\le3}^1$.  It is
straightforward\footnote{Because of the skew-symmetry, there is no normal
ordering necessary.} to check that for any~$s\in \Sfin(A)$ the sequence
$\{(Q_A)\mstrut _Ns\}\mstrut _N$~stabilizes, so that $(Q_A)\mstrut _N\to
Q_A'\in E_0(A)^1_{\le3}$ for some~$Q_A'$.  Also, $E_0(A)^1_{\le 1} $ is
canonically the finite dimensional space~$H_0^* = (\frak{z}_A)^*_\CC$, acting
by Clifford multiplication on~$\SS_0(A)\otimes {\tsize\bigwedge} ^0\subset
\Sfin(A)$.
 
        \proclaim{\protag{3.27} {Proposition}}
 There exists a unique $Q_A\in Q_A' + E_0(A)_{\le 1}^1\subset
E_0(A)^1_{\le3}$ such that for $\beta $~of finite energy,
  $$ -[\frac 14 Q_A,\gamma (\beta ^*)] =\cd(\beta ^\sigma _A),\qquad
     \beta \in H_{\operatorname{fin}}(A), \tag{3.28} $$
where $\cd\:(\LgR)^\sigma _{\operatorname{fin}}\to\ESfin$ is the
infinitesimal spin representation.  Furthermore, $A\mapsto Q_A$ is
$(\LGR)^\sigma $-equivariant.
        \endproclaim

        \demo{Proof}
 We compute the left hand side of~\thetag{3.28} for $\beta $~of fixed energy
(after complexification) and act on $s\in \Sfin(A)$ of fixed energy.  Thus we
can replace~$Q_A'$ by~$(Q_A)\mstrut _N$ for $N$~sufficiently large and
compute $-[\frac 14(Q_A)\mstrut _N,\gamma (\beta ^*)]$ in~$\Cliff^c\bigl(\HfinsA
\bigr)$.  This computation involves only finite sums, and in the notation
of~\S{1} \thetag{3.28}~becomes the equation
  $$ -[\frac 1{24}\gamma ^a\gamma ^b\gamma ^c,g_{de}\gamma ^e] = \sigma _d,
     \tag{3.29} $$
which follows immediately from~\thetag{1.7}.  In our present context the
indices~$a,b,c,\dots $ now correspond to a basis of a finite dimensional
subspace of the loop algebra.  If $\beta $~has nonzero energy, then $\beta
^\sigma _A$~also has a definite nonzero energy (\theprotag{2.41} {Lemma}),
and this is true of the commutator~\thetag{3.28} as well.  In this case
shifting~$Q_A$ by an element of~$E_0(A)^1_{\le 1}$ does not affect the
commutator.  For $\beta $~of energy zero we restrict to~$\SS_0(A)\otimes
{\tsize\bigwedge} ^0$, where the commutator~$-[\frac 14(Q_A)_0,\gamma (\beta
^*)]$ gives the standard finite dimensional splitting $\frak{z}\mstrut
_A\to\frak{z}_A^\sigma $, by~\thetag{3.29}.  The splitting $\beta \mapsto
\beta _A^\sigma $ on the right hand side of~\thetag{3.28}
differs\footnote{For $P\to \cir$ trivial and $A=A_0$ the trivial connections
these splittings agree; see \theprotag{3.24(2)} {Lemma}.}  by a functional
on~$\frak{z}_A$, and this specifies the difference $Q_A - Q_A'$
in~$E_0(A)_{\le 1}^1$.
 
For the equivariance, if $\varphi \in (\LGR)^\sigma $, then $\chi (\varphi
)\in \Aut(\SS)$ maps~$\Sfin(A)$ to~$\Sfin\bigl((\varphi\inv )^*A  \bigr)$,
and also maps~$\ONA$ to~$\Omega _N\bigl((\varphi\inv )^*A  \bigr)$.  The
splitting $\beta \mapsto \beta _A^\sigma $ is mapped to the splitting
$\beta \mapsto \beta ^\sigma _{(\varphi \inv)^*A}$, and the equivariance
of~$A\mapsto Q_A$ follows: 
  $$ Q_{(\varphi \inv )^*A} = \chi (\varphi )\,Q_A\,\chi
     (\varphi )\inv ,\qquad \varphi \in (\LGR)^\sigma,  \tag{3.30} $$
as operators on~$\Sfin\bigl((\varphi\inv )^*A  \bigr)$.

        \enddemo

 \subhead \S{3.3}.  A family of cubic Dirac operators on loop groups; main
theorem  \endsubhead

Let $\LGRt$~be an admissible graded central extension which is {\it positive
definite\/}\footnote{In~\cite{FHT3} we consider more general central
extensions which instead satisfy a nondegeneracy condition.  Our assumption
here implies that $\tms$~is positive definite on~$L\mstrut _P\frak{z}$, where
$\frak{z}=\cent\frak{g}$.} in the sense that $\form_\tau $ is positive
definite on~$\LgR$.  Then we use~$\form=\form_\tau $ to build the spin
representation---see the text and footnote preceding~\thetag{3.7}.  Suppose
$V$~is a $\zt$-graded irreducible positive energy representation of~$\LGRts$.
Then $W=V\otimes \SS$ is a positive energy representation of~$\LGRt$.

Fix a connection~$A\in \AR$.  Assume the minimal energy on~$V$
is~$\emin$; then $\emin$~is the minimal energy on~$W$ and 
  $$ W_{\emin}(A) = V_{\emin}(A)\otimes \SS_0(A) \tag{3.31} $$
is the minimal energy subspace.  Since $V$~is irreducible, $W_{\emin}(A)$~is
finite dimensional.  Also, $\SS_0(A)$~is an irreducible $\zt$-graded
$\Cliff^c(\zA^*)$-module, by \theprotag{3.24} {Lemma}.  For ~$\beta \in \LgR$
let $\rdA(\beta )$~be the infinitesimal action of~$\beta _A^{\tau }\in
\LgR^{\tau}$ on~$W$.  Note $\rdA(\beta )$~is typically unbounded.  If $\beta
$~has finite energy, then $\rdA(\beta )$~preserves the dense
subspace~$\Wfin(A)\subset W$ of finite energy vectors.
 
We first define the Dirac operator~$\DA$ as an unbounded (formally)
skew-Hermitian operator on the (algebraic, incomplete) inner product
space~$\Wfin(A)$.  Fix a Hilbert space basis~$\{e\mstrut _p\}$ of~$H_{\CC}$
such that each~$e\mstrut _p$ has a definite energy---see~\thetag{3.7}.
Let~$\gamma ^p$ be the dual Clifford multiplication~\thetag{3.16}
on~$\Sfin(A)$.  Recall from \theprotag{2.18} {Lemma} that $A$~defines a
splitting of $\LgRts\to\LgR$, and so an infinitesimal action~$(R_p)\mstrut
_A\in \End\bigl(\Vfin(A)\bigr)$.  Set
  $$ \DA = i\gamma ^p(R_p)\mstrut _A + \frac i2 Q_A\qquad \text{on
     $\Wfin(A)$}, \tag{3.32} $$
where $Q_A$~is defined in \theprotag{3.27} {Proposition}.  The sum
in~\thetag{3.32} is finite on any element of~$\Wfin(A)$, so $\DA$~is
well-defined.  Let $E_A$~be the energy operator on~$W$.

        \proclaim{\protag{3.33} {Proposition}}
\roster 
 \item$\DA$~is odd (formally) skew-Hermitian on~$\Wfin(A)$ and
$[\DA,E_A]=0$.\newline 
\vskip-12pt
 \item The dependence of~$\DA$ on~$A$ is 
  $$ \Dirac_{A + \beta @,ds} = \Dirac_A + \gamma (\beta ^*),\qquad \beta
     \in \LgR. \tag{3.34} $$
 \item The operator 
  $$ \opA = \DA(1-\DAs)^{-1/2} \tag{3.35} $$
extends to a bounded odd skew-adjoint operator on~$W$ with closed
range.\newline \vskip-12pt
 \item The assignment $A\to \opA$ is $\LGRt$-invariant.\newline
\vskip-12pt
 \item We have $\Ker \opA\subseteq \Wmin(A)$.  In particular, $\opA$~is an
odd skew-adjoint Fredholm operator.
       \endroster
        \endproclaim

        \demo{Proof}
 For~(1) note from~\thetag{3.14} that $\gamma ^p$ has energy opposite to that
of~$(R_p)\mstrut _A$.  Equation~\thetag{3.34} follows from \theprotag{2.18}
{Lemma} and \theprotag{3.27} {Proposition}.  Now~(1) implies that
$\DA$~preserves the eigenspaces of~$E\mstrut_A$, and $1-\DAs$ is positive
self-adjoint on each eigenspace.  Hence \thetag{3.35}~defines a bounded odd
skew-adjoint operator on the dense subspace~$\Wfin(A)\subset W$, so extends
by continuity to~$W$.  The $\LGR^{\tau}$-invariance of~$A\mapsto\DA$ follows
from \theprotag{2.18} {Lemma}, \thetag{3.6}, and~\thetag{3.30}; whence
$A\mapsto \opA$ is $\LGRt$-invariant as well.  This proves (1)--(4).
 
The extension guaranteed in \theprotag{2.42(1)} {Definition} is unique up to
a character of~$\Trot$, so the energy operator~$E\mstrut_A$ is determined up
to a constant independent of~$A$.  We claim
  $$ \text{$D_A^2 + 2E\mstrut_A$ is constant on~$W$, independent of~$A$.}
     \tag{3.36} $$
Then (5)~follows immediately, since $\Ker \opA = \Ker D_A^2$ and $D_A^2$~is
nonpositive.  As a first step to proving~\thetag{3.36} we compute that for
any $\beta \in \LgR$
  $$ [\DA,\gamma (\beta ^*)] = -2i\rdA(\beta ). \tag{3.37} $$
To derive this use~\thetag{3.17} for the first term of~\thetag{3.32} and
\thetag{3.28}~for the second term.  Next, square~\thetag{3.34} and
combine~\thetag{3.37} with~\thetag{3.17} and~\thetag{2.46} to deduce that
$D_A^2 + 2E\mstrut_A$~is independent of~$A$.  Now we show that $D_A^2 +
2E\mstrut_A$ commutes with the action of~$\LGRt$ and with any Clifford
multiplication~$\gamma (\beta ^*)$.  In fact, $g\in \LGRt$ conjugates $D_A^2
+ 2E\mstrut_A$ to $D_{A'}^2 + 2E\mstrut_{A'}$ for $A'$~the transform of~$A$
under the image of~$g$ in~$\LGR$.  For $D_A^2$ this follows from~(4) in the
proposition; for~$E\mstrut_A$ we use the fact that the connections~$d_A$
and~$d_{A'}$ are conjugate and the nullity condition~\thetag{2.14} is
$\LGRt$-invariant, since $\form_\tau $~is $\LGRt$-invariant.  By the previous
argument, $D_A^2 + 2E\mstrut_A = D_{A'}^2 + 2E\mstrut_{A'}$.  As for the
commutator with Clifford multiplication, we first use the infinitesimal
version of the $\LGRt$-invariance of~$A\mapsto\DA$, obtained
from~\thetag{3.34}, which asserts that
  $$ [\DA,\rdA(\beta )] = -\gamma \bigl([d_A,\beta ]^{*} \bigr) =
     -i[E\mstrut_A,\gamma (\beta ^*)],\qquad \beta \in
     \LgR. \tag{3.38} $$
Iterating~\thetag{3.37} and~\thetag{3.38} we find
  $$ \aligned
      [\DAs,\rdA(\beta )] &= [-2E_A,\rdA(\beta )], \\
      [\DAs,\gamma(\beta ^*)] &= [-2E_A,\gamma(\beta ^*)].\endaligned
      $$
In particular, $[D_A^2 + 2E\mstrut_A,\gamma (\beta ^*)]=0$, as desired.
Finally, since $\SS$~is an irreducible Clifford module and $V$~an irreducible
representation of~$\LGRt$, it follows that $D_A^2 + 2E\mstrut_A$, which
commutes with all operators from~$\LGRt$ and the Clifford algebra, is a
constant on~$W=V\otimes \SS$.
        \enddemo

Recall that the graded central extension~$\LGRt$ of~$\LGR$ defines a
twisting~$\th$ of~$K_G(\GR)$---see~\thetag{2.9}.

        \proclaim{\protag{3.39} {Corollary}}
 Let $\LGRt$~be a positive definite admissible graded central extension
of~$\LGR$, and let $V$~be a finitely reducible $\zt$-graded positive energy
representation of~$\LGRts$.  Then
  $$ \aligned
      \Dirac(V)\:\AR &\longrightarrow \Fred(V\otimes \SS) \\
      A&\longmapsto \opA = \DA(1-\DAs)^{-1/2} \endaligned \tag{3.40} $$
represents an element of~$\tKGG$.  Furthermore, it only depends on~$V$ up to
isomorphism and is additive, so defines a map
  $$ \Phi \:\Rts\longrightarrow \tKGG. \tag{3.41} $$
        \endproclaim

\flushpar
 In the odd case the $\cco$-action comes from that on~$\SS_0(A)\subset \SS$,
since $\dim Z_A = \dim G\pmod2$.

        \demo{Proof}
 Recall from~\cite{FHT1,Appendix~A.5} that the family~$\{\opA\}$ of
skew-adjoint Fredholm operators represents an element of twisted $K$-theory
if $\opA^2+1$~is compact and $A\mapsto (\opA,\opA^2+1)$ is continuous as a
map into the product ~$\Cal{B}(W)\times \Cal{K}(W)$ of bounded and compact
operators.  Here $\Cal{B}(W)$~has the compact-open topology and
$\Cal{K}(W)$~the norm topology.  By the remarks in~\cite{AS,Appendix~1} for
this purpose we can replace the compact-open topology on~$\Cal{B}(W)$ with
the strong operator topology.  From~\thetag{3.36} we see that
$\opA^2+1=(1-D_A^2)\inv $~is the inverse of~$2E_A+C$ for some constant~$C$,
and the latter is a positive operator.  By \theprotag{2.45} {Lemma} the
energy operator of a finitely reducible positive energy representation~$V$
has eigenvalues of finite multiplicity tending to infinity.  The same is true
for the spin representation~$S$ by explicit construction~\thetag{3.15}, so it
is also true for the energy operator of the tensor product~$W$, whence the
inverse of~$2E_A+C$ is compact.
 
For the continuity we prove first that $A\mapsto(\opA^2+1)=(1-D_A^2)\inv $ is
norm continuous.  To ease the notation set~$y=D_A$ and $x=D_{A+\beta
\,ds}=y+b$, where $b=\gamma (\beta ^*)$ is a bounded operator continuously
varying with~$\beta $; see~\thetag{3.34}.  Then
  $$ (1-x^2)\inv - (1-y^2)\inv =\left[ (1-x^2)\inv x \right]b(1-y^2)\inv +
     (1-x^2)\inv b\left[ y(1-y^2)\inv \right].   $$
Consideration of finite energy vectors leads to the bounds $\|(1-x^2)\inv \|\le
1,\,\|x(1-x^2)\inv \|\le1$ on the operator norms.  The same estimates hold
for~$y$ replacing~$x$, so that 
  $$ \|(1-x^2)\inv - (1-y^2)\inv \| \le 2\|b\| \le C\|\beta \|  $$
for some constant~$C$, which proves $A\mapsto \opA^2+1$ is norm continuous
into~$\Cal{K}(W)$.  To show the strong continuity of $A\mapsto\opA$ we write
  $$ (1-x^2)^{-1/2}x - (1-y^2)^{-1/2}y = \left[ (1-x^2)^{-1/2} -
     (1-y^2)^{-1/2} \right]y + (1-x^2)^{-1/2}b.  \tag{3.42} $$
Now the square root map $z\mapsto z^{1/2}$ on bounded positive operators is
continuous (even analytic) in the uniform topology~\cite{Hi}.  Thus since
$A\mapsto (1-D_A^2)\inv $ is norm continuous, so is $A\mapsto
(1-D_A^2)^{-1/2}$.  It then follows directly from~\thetag{3.42} that
$A\mapsto \opA(v)$ is continuous for~$v\in \Wfin(A)$.  The strong continuity
of $A\mapsto \opA$ is a consequence of the following lemma, since
$\|\opA\|\le1$.

        \proclaim{\protag{3.44} {Lemma}}
  Let $W$~be a Hilbert space, $V\subset W$ a dense subspace, and $T\:X\to
B(W)$ a family of bounded operators parametrized by a space~$X$.  Assume that
$x\mapsto T_xv$ is continuous for all~$v\in V$ and that $\|T_x\|\le C$ for
all~$x\in X$.  Then $T$~is strongly continuous.  
        \endproclaim

        \demo{Proof}
 For $w \in W$ we must show that $x\mapsto T_xw $ is continous.  Fix
$x_0\in X$ and a sequence $v_n\to w $ with $v_n\in V$.  Now given~$\epsilon
>0$ choose~$N$ so that $\|w -v_N\|<\epsilon /4C$ and a neighborhood~$U$
of~$x_0$ so that $\|(T_x-T_{x_0})v_N\|<\epsilon /2$ for all~$x\in U$.  Then
  $$ \|(T_x - T_{x_0})w \| \le \|(T_x - T_{x_0})(w -v_N)\| + \|
     (T_x-T_{x_0})v_N\| < 2C\frac{\epsilon }{4C}+ \frac{\epsilon
     }{2}=\epsilon $$
if $x\in U$, which proves the continuity. 
        \enddemo

        \enddemo

We can now assert our main theorem. 

        \proclaim{\protag{3.43} {Theorem}}
 Let $\LGRt$~be a positive definite admissible graded central extension
of~$\LGR$.  The map $\Phi \:\Rts\to\tKGG$ is an isomorphism of graded free
abelian groups.
        \endproclaim

\flushpar 
 The proof of \theprotag{3.43} {Theorem} in the general case is
in~\cite{FHT3,\S13}.  In the next section we present the proof in case $G$~is
connected and $\pi _1G$~is torsion-free.

 \head
 \S{4} Proofs
 \endhead
 \comment
 lasteqno 4@ 30
 \endcomment

Let $G$~be a connected Lie group.  A principal $G$-bundle $P\to\cir$ is then
trivializable, so in the sequel we take it to be the trivial
bundle~$P=\cir\times G$.  Let $A_0\in \AR$ be the trivial connection.  Since
$P\to\cir$~is trivial the extended loop group is $\LGRR\cong LG\rtimes
\Troth$ for some finite cover of $\Troth\to\Trot$ of degree~$\delta $.  An
admissible central extension necessarily has the form $\LGRRt\cong
LG^\tau\rtimes \Troth$.  It follows that on any positive energy
representation the energy operator~$E_{A_0}$ satisfies $\exp(2\pi i\delta
E_{A_0})=\operatorname{id}$, so its eigenvalues are $1/\delta $~times
integers.  By tensoring with a character of~$\Troth$ we normalize the minimal
$E_{A_0}$-energy to be zero.  We allow the central extension~$LG^\tau $ to
have a nontrivial grading.  Now the stabilizer $Z_{A_0}\subset LG$ of the
trivial connection~$A_0$ is the group of constant loops~$G$.  Let $G^\tau \to
G$ be the restriction of the graded central extension $LG^\tau \to LG$ to the
constant loops.  It has a trivial grading since $G$~is connected.
 
Assume the admissible graded central extension $LG^\tau \to LG$ is positive
definite, i.e., the form $\form_\tau $ is positive definite on~$L\frak{g}$.
Thus it restricts to a positive definite form~$\langle \cdot ,\cdot
\rangle_\tau $ on constant loops~$\frak{g}\subset L\frak{g}$, which in turn
induces an isomorphism
  $$ \aligned 
      \kt\:\frak{g}&\longrightarrow \frak{g}^* \\
      \xi &\longmapsto \xi ^*, \endaligned\tag{4.1} $$
where $\xst(\eta ) = \langle \xi ,\eta \rangle_\tau $ for all~$\eta \in
\frak{g}.$   
 
Fix a maximal torus~$T\subset G$ with Lie algebra~$\frak{t}\subset \frak{g}$,
let $W$~be the associated Weyl group, and define the dual lattices~$\Pi
=\Hom(\TT,T)$ and~$\Lambda =\Hom(T,\TT)$.  Recall from~\S{1.5} the affine
space~$\aff T\tau $ and the subset~$\Lambda ^\tau \subset \aff T\tau $ which
is a $\Lambda $-torsor.  The lattice~$\Pi =\exp\inv (1)/2\pi \subset
\frak{t}$ acts by translation on~$\Lambda ^\tau $ via the map~$\kt$.  The
Weyl group~$W$ acts as well, and these actions fit together into an action of
the {\it extended affine Weyl group\/}~$\Waff=\Pi \rtimes W$.  Let
$\Lreg\subset \Lambda ^\tau $ be the subset of {\it affine
regular\/}\footnote{A weight~$\mu \in \Lambda ^{\tau }$ is called {\it affine
regular\/} if it has trivial stabilizer under the action of~$\Waff$.
Equivalently, by \theprotag{4.9(1)} {Lemma} below, if $\mu =\kappa ^{\tau
}(\xi )$ for~$\xi \in \frak{t}$, then $\mu $~is affine regular if and only if
$\exp(2\pi \xi )\in T$ is regular in the sense that no nontrivial element of
the Weyl group~$W$ fixes it.} weights.

To prove \theprotag{3.43} {Theorem} for $G$~connected with $\pi
_1G$~torsion-free we first identify the set of isomorphism classes of
irreducible positive energy representations of~$LG^{\tms} $ with certain
$\Waff$-orbits of projective weights.  This step is essentially a quotation
from~\cite{PS,\S9}.  There is a correspondence between conjugacy classes
in~$G$ and $\Waff$-orbits in~$\aff T{\tms} $, so each irreducible
representation determines a conjugacy class.  Collectively these are the {\it
Verlinde conjugacy classes\/}.  The main step in the proof is the assertion
that the kernel of the Dirac family~\thetag{3.40} associated to a given
irreducible representation is supported on the conjugacy class which
corresponds to the given irreducible.  Then the following results
from~\cite{FHT1} (Theorem~4.27 and Proposition~4.41) complete the proof.

        \proclaim{\protag{4.2} {Theorem}}
 \roster
 \item"\rom(i\rom)" For $G$~connected of rank~$n$ with $\pi _1G$~torsion-free, 
  $$ K_G^{\tau +q}(G)\cong \cases \Hom\mstrut
     _{\Waff}\bigl(\Ltil,H^n_c(\frak{t})\otimes \ZZ(\epsilon ) \bigr)
     ,&q\equiv n\pmod2;\\0,&q\not\equiv n\pmod2.\endcases \tag{4.3} $$
The $\Waff$-action on~$\Ltil$ is described in the previous paragraph;
$W\subset \Waff$ acts on~$H^n_c(\frak{t})$ naturally via the sign
representation; the grading of~$\tau $ determines a homomorphism $\epsilon
\:\Pi \to\zt$, and $\Pi \subset \Waff$~acts on~$\ZZ(\epsilon )$ via this
sign.
 \item"\rom(ii\rom)" If $k\:H\hookrightarrow G$ also has rank~$n$ and is also
connected with $\pi _1H$~torsion-free, then $k^*\:K_G^{\tau +q}(G)\to
K_H^{k^*\tau +q}(H)$ is, for $q\equiv n\pmod2$, the natural inclusion
  $$ \Hom\mstrut _{\Waff(G)}\bigl(\Ltil,H^n_c(\frak{t})\otimes \ZZ(\epsilon )
     \bigr)\longhookrightarrow \Hom\mstrut
     _{\Waff(H)}\bigl(\Ltil,H^n_c(\frak{t}) \otimes \ZZ(\epsilon )\bigr).
     \tag{4.4} $$
 \item"\rom(iii\rom)" Orient~$T$ and let $j\:\{t\}\hookrightarrow T$ be the
inclusion of a point.  Then the pushforward~$j_*$ for $q\equiv \dim T\pmod2$
is computed as the map
  $$ \xymatrix{K_T^{i^*k^*\tau +q}\bigl(\{t\} \bigr)\ar[r]^-{j_*} \ar[d] &
     K_T^{k^*\tau +q}(T)\ar[d]\\
      \ZZ[\Ltil] \ar[r] &\Hom\mstrut _\Pi \bigl(\Ltil,\ZZ(\epsilon )\bigr)}
     \tag{4.5} $$
which takes~$\mu \in \Ltil$ to the $\Pi $-equivariant function $f \:\Lambda
^\tau \to\ZZ$ which is supported on the $\Pi $-orbit of~$\mu $ and
satisfies~$f(\mu )=1$.
	\endroster
        \endproclaim

\flushpar
 Notice that an orientation of~$T$ identifies the compactly supported
cohomology~$H^n_c(\frak{t})$ with~$\ZZ$.  The nonzero abelian group
in~\thetag{4.3} is free, and the set of {\it regular\/} $\Waff$-orbits
in~$\Ltil$ provides a generating set. 

We carry out the proof separately for tori and simply connected groups, then
combine them to prove the general case.

 \subhead \S{4.1}. The proof for tori
 \endsubhead

Let $G=T$~be a torus with Lie algebra~$\frak{t}$ and~$\Pi =\Hom(\TT,T)$.
Recall the decomposition $LT\cong T\times \Pi \times U$ in~\thetag{2.26}.  By
\theprotag{2.27(1)} {Proposition} any admissible central extension~$LT^\tau $
is a product of Heisenberg central extensions~$(T\times \Pi )^\tau $
and~$U^\tau $.  We first assume the grading of~$\tau $ is trivial; below we
discuss the modifications for nontrivial grading.  Since $\tau $~is positive
definite, $L\frak{t}/\frak{t}$ is a symplectic vector space and its
Heisenberg extension~$U^\tau $ has a unique irreducible positive energy
representation~$\sH$.  (See~\cite{PS,\S9.5}.)  Fix an orbit~$\orb$ of the
$\Pi $-action on~$\Lambda ^\tau $ and consider the constant vector bundle
with fiber~$\sH$ over~$\orb$.  The space~$V_{\orb}$ of $L^2$~sections (with
respect to the measure which assigns unit mass to each point of~$\orb$) is an
irreducible representation~$V_{\orb}$ of~$LT^\tau $: the subgroup~$U^\tau $
acts on each fiber; $\Pi ^\tau $~permutes the fibers; and $T^\tau $~acts on
the fiber at~$\lambda \in \Lambda ^\tau \subset \Hom(T^\tau ,\TT)$ by scalar
multiplication, the character~$\lambda $ defining the multiplication.  These
are all of the irreducible positive energy representations up to isomorphism.
 
        \proclaim{\protag{4.6 \cite{PS,(9.5.11)}} {Proposition}}
 The isomorphism classes of irreducible positive energy representations
of~$LT^\tau $ are in 1:1~correspondence with the orbits of~$\Pi $ on~$\Lambda
^\tau $, i.e., with the points of~$\Lambda ^\tau/\kt(\Pi )$.
        \endproclaim
 
Turning to the Dirac family associated to the irreducible
representation~$V_{\orb}$, let $S$~be the spin representation of~$LT$ and set
$W_{\orb}=V_{\orb}\otimes S$.  Any connection on the trivial bundle
over~$\cir$ is gauge equivalent to $A_0 + \xi \,ds$ for some constant~$\xi\in
\frak{t}$.  Let~$D_\xi $ denote the corresponding Dirac
operator~\thetag{3.32} and $E_\xi $~the associated energy operator.
According to~\thetag{3.36} if $D_\xi (w)=0$ for some ~$w\in W_{\orb}$ of unit
norm, then $\langle E_\xi (w),w\rangle_{W}$~is a global minimum of~$E_\xi $
over all~$\xi \in \frak{t}$ and all unit norm~$w\in W_{\orb}$.  The energy
operator on~$\sH$ has discrete nonnegative spectrum and a one-dimensional
kernel spanned by a unit norm vector~$\Omega \in \sH$.  Let $\Omega_\lambda
$~denote the copy of~$\Omega $ in the fiber at~$\lambda \in \orb$.  Then
\thetag{2.32}, \thetag{2.33}, and~\thetag{2.46} imply
  $$ \bigl\langle E_\xi (\Omega _\lambda ),\Omega _\lambda \bigr\rangle_{\sH}
     = \frac 12\,\bigl| \lambda -\kt(\xi )\bigr|_\tau ^2 + C\tag{4.7} $$
for some real constant~$C$.  We can take~$C=0$ so that the global minimum of
energy is zero.  Let $S_0$~denote the zero energy subspace of~$S$; it is an
irreducible Clifford module for the algebra~$\Cliff(\frak{t}^*)$.  For
each~$\lambda _0\in \orb$ there is a unique~$\xi_0 \in \frak{t}$ which
makes~\thetag{4.7} vanish---it satisfies $\kt(\xi _0)=\lambda _0$---and
these~$\xi _0$ form a $\Pi $-orbit in~$\frak{t}$.  For each such~$\xi _0$ the
kernel of~$D_{\xi _0}$ is $\CC\cdot \Omega _{\lambda _0}\otimes S_0$.
Furthermore, \thetag{3.34}~shows that $D_\xi $~acts on $\CC\cdot \Omega
_{\lambda _0}\otimes S_0$ as Clifford multiplication by~$\kt(\xi -\xi _0)$.
Now $\Pi $-orbits in~$\frak{t}$ correspond to elements of~$T$, and from this
point of view the support of the kernel is a single element in~$T$.  More
precisely, the linear splitting of \theprotag{2.18} {Lemma} for the trivial
connection~$A=A_0$, restricted to constant loops, is a basepoint in the
affine space~$\aff T\tau $.  Hence \thetag{4.1}~may be regarded as an
isomorphism $\kt\:\frak{t}\to\aff T\tau $, which in turn induces an
isomorphism
  $$ \ktb\:T = \frak{t}/\Pi \longrightarrow \aff T\tau /\kt(\Pi
     ).  $$
The preceding proves the following.

        \proclaim{\protag{4.8} {Proposition}}
 Let $V_{\orb}$~be the irreducible positive energy representation of~$LT^\tau $
which corresponds to~$\orb\in \Ltil/\kt(\Pi )\subset \aff T\tau /\kt(\Pi )$.
Define~$t\in T$ by $\ktb(t)=\orb$.  Then the kernel of the Dirac family
associated to~$V_{\orb}$ is supported on~$\{t\}\subset T$.  Let $i\:\{t\}\to
T$ be the inclusion, $i_*\:K^{i^*\th}_T(\{t\})\to K^{\th + \dim T}_T(T)$ the
induced pushforward, and $K\to\{t\}$ the $T^\tau $-line with
character~$\lambda $, where $\lambda $~is some element of~$\orb\subset
\Ltil$.  Then the $K$-theory class of the Dirac family is~$i_*[K]$.
        \endproclaim

\theprotag{3.43} {Theorem} for~$G=T$ is now a direct consequence of
\theprotag{4.8} {Proposition} and \theprotag{4.2} {Theorem}.
 
A grading of~$LT$ is a homomorphism $LT\to\zt$, and it necessarily factors
through a homomorphism $\epsilon \:\Pi \to\zt$ since $\Pi \cong \pi _0(LT)$.
Suppose now the admissible central extension~$LT^\tau $ includes the
grading~$\epsilon $.  The representation~$V_{\orb}$ associated to a $\Pi
$-orbit~$\orb\subset \Ltil$ depends mildly on a basepoint~$\lambda \in \orb$.
Namely, define a $\zt$-graded vector bundle over~$\orb$ whose fiber at~$\pi
\cdot \lambda ,\; \pi \in \Pi,\, $ is~$\sH$ with grading~$\epsilon (\pi )$.
Then $V_{\orb,\lambda }$~is the $\zt$-graded Hilbert space of $L^2$~sections.
A change of basepoint may replace~$V_{\orb,\lambda }$ with the oppositely
graded representation, whose equivalence class in the ring~${R}^\tau (LT)$ is
negative to that of~$V_{\orb, \lambda }$.  (See \theprotag{2.50}
{Definition}.)  The analysis of the Dirac operator proceeds as before, but
now the kernel of the Dirac family is a $\zt$-graded line bundle over the
orbit~$\orb$.  The function~$f$ which represents~$i_*(\lambda )$
in~\thetag{4.5} is the rank of this kernel; it takes values~$\pm1$ due to the
nontrivial action of~$\Pi $ on~$\ZZ(\epsilon )$ in~$\Hom\mstrut _\Pi
\bigl(\Ltil,\ZZ(\epsilon )\bigr)$.

 \subhead \S{4.2}. The proof for simply connected groups
 \endsubhead

Let $G$~be a simply connected compact Lie group and fix a positive definite
admissible central extension~$LG^\tau $ of its loop group.  The
grading~$\epsilon $ is trivial since $LG$~is connected.  Fix a maximal
torus~$T$ and choice of Weyl chamber in~$\frak{t}$.  The set of conjugacy
classes in~$G$ is isomorphic to the set of $W$-orbits in~$T$.  Let $\Delta
$~be the set of roots of~$G$ and $\Delta ^+$~the set of positive roots.  The
following lemma is standard; see~\cite{PS,\S5.1} for example.

        \proclaim{\protag{4.9} {Lemma}}
 \roster
 \item $\exp(2\pi \;\cdot )\:\{\text{$\Waff$-orbits in
$\frak{t}$}\}\longrightarrow \{\text{$W$-orbits in $T$}\}$ is a
bijection.\newline\vskip-12pt
 \item Define the {\it alcove\/} 
  $$ \alcove=\{\xi \in \frak{t}: -1 <\alpha (\xi )< 0 \text{ for
     all $\alpha \in \Delta ^+$}\}. \tag{4.10} $$
Then $\alcove$~is contained in a Weyl chamber of~$\frak{t} $ and its
closure~$\calcove$ is a fundamental domain for the $\Waff$-action
on~$\frak{t}$.
 \endroster
        \endproclaim

We do not review here the construction of positive energy representations,
but merely quote the classification.  Recall that $LG^\sigma $~is the central
extension defined by the spin representation.  Let $E_0$~be the energy
operator of the connection~$A_0$.

        \proclaim{\protag{4.11~\cite{PS,(9.3.5)}} {Proposition}}
 \roster
 \item"(i)" If $\tau -\sigma $ is not positive semidefinite on~$L\frak{g}$,
then there are no nonzero positive energy representations of~$LG^{\tau
-\sigma }$.
 \item"(ii)" If $\tau -\sigma $ is positive semidefinite, then the isomorphism
classes of irreducible positive energy representations of~$LG^{\tau -\sigma
}$ are in 1:1~correspondence with the orbits of~$\Waff$ on~$\Lambda ^{\tms}
$, i.e., with the points of~$\Lambda ^{\tms} /\kts(\Waff)$.  Furthermore, the
representation~$V$ which corresponds to~$-\lambda \in \kts(\calcove)$ has
zero $E_0$-energy space~$V_0$ the irreducible representation of~$G$ with
lowest weight~$-\lambda $.
 \endroster
        \endproclaim

\flushpar
 Assertion~(i) is \theprotag{2.47} {Proposition}.  There is a $\rho $-shifted
restatement of~(ii) as follows.  First, $\sigma $~restricts on point loops to
the finite dimensional spin extension by \theprotag{3.24(i)} {Lemma}, and
\theprotag{1.15} {Lemma} shows~$\rho \in \Lambda ^{\sigma }$, where $\rho
$~is half the sum of the positive roots.  The map $-\lambda \mapsto -(\lambda
+\rho )$ maps $\Lambda ^{\tms}\cap \kts(\calcove)$ isomorphically to~$\Lambda
^\tau \cap \kt(\alcove)$, which consists of {\it affine regular\/}
weights~\cite{A,5.62}.  Then the isomorphism classes of irreducible positive
energy representations of~$LG^{\tms}$ are in 1:1~correspondence with the set
of affine regular $\Waff$-orbits in~$\Lambda ^{\tau }$, i.e., with the points
of~$\Lreg/\kt(\Waff)$.

We turn now to the Dirac family~\thetag{3.40}, assuming that $\tms$~is
positive semidefinite.  Let $V$~be a $\zt$-graded irreducible positive energy
representation~$V$ of~$LG^{\tau -\sigma }$ whose zero $E_0$-energy space has
lowest weight~$-\lambda _0$.  The minimal $E_0$-energy subspace~\thetag{3.31}
of $W=V\otimes \SS$ is
  $$ W_0 = V_0\otimes \SS_0,  $$
where $\SS_{0}$ is the zero $E_{0}$-energy subspaces.  Note that $V_0$~is a
finite dimensional $\zt$-graded representation of~$G^{\tau -\sigma }$ and
$\SS_0$~is a finite dimensional $\zt$-graded representation of~$G^\sigma $.
\theprotag{1.15} {Lemma} implies that a lowest weight of~$W_0$ is~$-\mu
_0=-(\lambda _0+\rho )$.  Write $-\mu _0=\kt(\xi _0)$ for~$\xi _0\in
\alcove$.  Since $\calcove$~is a fundamental domain for the $\Waff$-action
on~$\frak{t}$, any connection on the trivial bundle is gauge equivalent
to~$A_0+\xi \,ds$ for~$\xi \in \calcove$.  Denote the associated Dirac and
energy operators as~$D_\xi ,E_\xi $.

        \proclaim{\protag{4.12} {Proposition}}
 \roster
 \item $\Ker D_\xi \subset W_0$ if~$\xi \in \calcove$.
 \newline\vskip-12pt
 \item The restriction of~$D_\xi $ to~$W_0$ may be identified with the
operator~$D_\mu (V_0)$ of~\thetag{1.32} with~$\mu =\kt (\xi )$.
 \endroster
        \endproclaim

\flushpar 
 Recall that the operators~$D_\mu $ in~\thetag{1.32} are parametrized by~$\mu
\in \aff G\tau \cong \aff G{\tms}$, where the isomorphism is given by the
canonical splitting~\thetag{1.5} of the Lie algebra extension at
level~$\sigma $.  Observe that the linear splitting of \theprotag{2.18}
{Lemma} for the trivial connection~$A=A_0$, restricted to constant loops,
gives a basepoint in~$\aff G\tau $ and~$\aff G\tms$, and by \theprotag{3.24}
{Lemma} these correspond under the isomorphism $\aff G\tau \cong \aff
G\tms$.  We use the basepoint to regard~$\kt$ as a map~$\frak{g}\to\aff
G\tau$. 

        \demo{Proof}
As in the discussion preceding~\thetag{4.7} $D_\xi =0$ precisely on vectors
in~$W$ which realize the global minimum of~$E_\xi $ (over all~$\xi $).
Now $W$~is an irreducible representation of the semidirect product
of~$\Cliff^c(L\frak{g}^*)$ and~$L\frak{g}^\tau $, so the set of vectors of
the form
  $$ \multline
      w=\gamma \bigl((z^{-n_{r+s}}\eta _{r+s})^* \bigr)\dots \gamma 
     \bigl((z^{-n_{r+1}}\eta _{r+1})^* \bigr)\rtd_{A_0}(z^{n_r}\eta _r)\dots
     \rtd_{A_0}(z^{n_1}\eta _1)\Omega,\\ 
      r,s\ge0,\quad \eta _i\in \frak{g},\quad n_i>0,\quad
     \Omega \in W_0,\endmultline  $$
is dense in~$W$.  Here $\gamma $~is the Clifford action~\thetag{3.16} and
$\rtd_{A_0}$~the action of~$L\frak{g}$ defined after~\thetag{2.43}.  Assume
that $\eta _i$~lies in the $\alpha _i$-root space of~$\frak{g}$ or
in~$\frak{t}$, in which case~$\alpha _i=0$.  Since $\rdA(z^{n_i}\eta _i) =
\rtd_{A_0}(z^{n_i}\eta _i)$---see \theprotag{2.18} {Lemma} and observe $\ll\!
z^{n_i}\eta _i,\xi \!\gg=0$ since~$n_i>0$---we see that the $E_\xi $-energy
of~$w$ minus the $E_\xi $-energy of~$\Omega $ is
  $$ \sum\limits_{i=1}^{r+s}\bigl(n_i + \alpha _i(\xi ) \bigr). \tag{4.13} $$
For $\xi \in \calcove$ this is nonnegative, by~\thetag{4.10}, which implies:
For $\xi \in \calcove$, if $\Ker D_\xi \not= 0$ then $\Ker D_\xi \cap
W_0\not= 0$.  Suppose $\Omega $~lies in the $\mu $-weight space of~$W_0$.
Then by~\thetag{2.46} its $E_\xi $-energy is
  $$ \mu (\xi ) + \frac{\|\xi \|_\tau ^2}{2},  $$
which has a minimum value of $-\|\mu \|_\tau ^2/2$ uniquely realized at~$\xi
\in \frak{t}$ which satisfies~$\kt(\xi )=-\mu $.  So the global minimum
occurs for the lowest weight~$\mu =-\mu _0$ at the point ~$\xi =\xi
_0\in\calcove$.  Furthermore, since $|\alpha _i(\xi _0)|<1$ by~\thetag{4.10}
we see from~\thetag{4.13} that $E_{\xi _0}$ does not realize its minimum on
any $E_0$-eigenspace~$W_n$ for~$n>0$.
      
 For \therosteritem{2} consider first~$\xi =0$.  In the restriction of the
first term in~\thetag{3.32} to~$W_0$ only terms of zero energy contribute to
the infinite sum.  That this is the first term of~\thetag{1.32} follows from
the remarks preceding the proof.  As for the second term in~\thetag{3.32}, in
view of \theprotag{3.24(2)} {Lemma} the characterization~\thetag{3.28}
of~$Q_{A_0}$, restricted to~$\SS_0$, also characterizes the second term
in~\thetag{1.32}.  This proves~\therosteritem{2} for~$\xi =0$.  For $\xi
\not= 0$ it follows from~\thetag{3.34} that $D_\xi =D_0 + \gamma (\xst )$,
which matches~\thetag{1.33} for~$\nu =\xst$.
        \enddemo

\theprotag{1.19} {Proposition} and \theprotag{1.24} {Proposition}, or rather
their projective analogs described in~\S{1.5}, translate to an explicit
description of~$\Ker D_\xi ,\,\xi \in \calcove$, hence of the $K$-theory
class~$\Phi (V)$ in~\thetag{3.41}; we use \theprotag{4.12} {Proposition} to
reduce to consideration of the family of operators $D_\mu (V_0)$.  (The Dirac
operator~$D_\xi $ is invertible on the orthogonal complement~$W_0^\perp$
to~$W_0$, so the Dirac family restricted to~$W_0^\perp$ represents the zero
element of twisted $K$-theory.)  Denote by $k\:T\hookrightarrow G$ be the
inclusion of the maximal torus, and let $N\subset G$ be the normalizer
of~$T$.  There is a restriction map
  $$ k^*\:K_G^{\th+\bullet}(G)\longrightarrow
     K_T^{\ksth+\bullet}(T)  $$
computed topologically in \theprotag{4.2(ii)} {Theorem}.  The family of
operators~$D(V)$ parametrized by~$\AR$ represents a class $[D(V)]\in
K_G^{\tau + \dim G}(G)$ (see~\thetag{3.40}); its restriction to~$\hol\inv
(T)$---connections with holonomy in~$T$---is denoted~$k^*D(V)$ and it
represents $k^*\bigl[D(V)\bigr]\in K_T^{\ksth + \dim T}(T)$.  (We use
periodicity and $\dim G\equiv\dim T\pmod2$.)

        \proclaim{\protag{4.14} {Proposition}}
 Let $V$~be an irreducible positive energy representation of~$LG^{\tms}$
whose zero $E_0$-energy space has lowest weight~$-\lambda _0$.  Define~$\xi
_0\in \alcove$ by $\kt(\xi _0)=-(\lambda _0+\rho )$.  Then the kernel
of~$k^*D(V)$ is supported on~$\hol\inv (\orbit)$ for
$i\:\orbit\hookrightarrow T$ the {\rm regular\/} $W$-orbit which
contains~$t=\exp(\xi _0)$.  Also,
  $$ k^*\bigl[D(V)\bigr] = i_*[K\otimes L] \tag{4.15} $$
is the image of the $K$-theory class of an $N^\tau $-equivariant line bundle
$K\otimes L\to\orbit$ under the pushforward
  $$ i_*\:K^{i^*\ksth} _N(\orbit) \longrightarrow K^{i^*\ksth}
     _T(\orbit)\longrightarrow K_T^{\ksth+\dim T}(T) \tag{4.16} $$
The stabilizer~$T^{\tms} \subset N^{\tms }$ acts on the fiber~$K_t$ by the
character~$-\lambda _0$ and $T^\sigma $~acts on~$L$ by the character~$-\rho $.
        \endproclaim

        \demo{Proof}
 The pushforward~\thetag{4.16} is realized as multiplication by spinors on
the tangent space~$\frak{t}$, which is identified with~$\frak{t}^*$
via~$\kt$~\thetag{4.1}.  The description of the kernel follows directly from
\theprotag{1.19} {Proposition}, and \thetag{4.15}~follows
from~\thetag{1.30}.  For the last statement, see \theprotag{1.24}
{Proposition}. 
        \enddemo

        \demo{Proof of~\theprotag{3.43} {Theorem} \rom($G$~
simply connected\rom)}
 Consider the composition
  $$ {R}^{\tau -\sigma }(LG)@>\hphantom{K}\Phi\hphantom{K} >> K_G^{\tau +
     \dim G}(G) @> \hphantom{K}k^* \hphantom{K}>> K_T^{\ksth+\dim
     T}(T). \tag{4.17} $$
First, the distinguished basepoint of~$\orbit$ in~$\orbit\cap\kt(\alcove)$
leads to an isomorphism $K_N^{i^*k^*\tau }(\orbit)\cong \ZZ[\Ltil]$ and to a
function $\orbit\to\zt$ which is the sign representation of the Weyl group
after identifying $\orbit\cong W$ using the basepoint.  Then if~$t\in \orbit$
the restriction map $K_N^{i^*k^*\tau }\bigl(\orbit\bigr)\to K_T^{i^*k^*\tau
}\bigl(\{t\} \bigr)$ is $\mu \mapsto\pm\mu $ according to this sign, and so
by \theprotag{4.2(iii)} {Theorem} the map \thetag{4.16} is $\mu \mapsto f_\mu
$, where $f_\mu \:\Ltil\to\ZZ$ is the $\Pi $-invariant function supported on
the $\Waff$-orbit of~$\mu $ with $f_\mu (w\cdot \mu )= (-1)^{\sign(w)}$
for~$w\in W$.  \theprotag{4.11} {Proposition} identifies~$R^{\tau -\sigma
}(LG)$ as the free abelian group on $\kts(\calcove)\cap\Lambda^{\tms}$, and by
\theprotag{4.14} {Proposition} the composition~\thetag{4.17} sends $-\lambda
_0\in \kts(\calcove)\cap\Lambda ^{\tms}$ to~$f_\mu $ for $\mu =-(\lambda
_0+\rho )$.  It follows from \theprotag{4.2(ii)} {Theorem} that these~$f_\mu
$ form a basis for the image of~$k^*$ in~\thetag{4.17}, and so $\Phi $~is an
isomorphism.
        \enddemo

We remark that if $\tms$~is not positive semidefinite, then both sides
of~\thetag{3.41} vanish.

 \subhead \S{4.3}. The proof for connected~$G$ with $\pi _1$~torsion-free
 \endsubhead

If $G$~is connected with $\pi _1G$~free, then there is a finite cover
  $$ 1 @>>> A @>>> \Gt  @>>> G @>>> 1  $$
such that $\Gt=Z_1\times G'$ is the product of a torus and a simply connected
compact Lie group~$G'$.  This follows easily from the proof
of~\cite{FHT1,Lemma~4.1}. There is an induced exact sequence of groups
  $$ 1 @>>> A @>>> L\Gt @>i>> LG @>>> A @>>> 1, \tag{4.18} $$
where $A$~is included in~$L\Gt$ as the point loops, $L\Gt$~maps onto a union
of components $i(L\Gt)\subset LG$, and the last map factors as $LG \to \pi
_0LG\cong \pi _1G \to A$.  Note also $L\Gt\cong LZ_1\times LG'$.   

        \proclaim{\protag{4.19} {Lemma}}
 The pullback $L\Gt^{\ttl-\stl}$ factors as $\ttl - \stl = (\tZ)\times (\tau
'-\sigma ')$, i.e., the spin extension~$\sigma \mstrut _Z$ is trivial and
there are no cross terms.  Similarly, $\ttl=\tZ \times \tau '$.
        \endproclaim

        \demo{Proof}
 The vanishing of~$\sigma \mstrut _Z$ is immediate from the triviality of the
pointwise adjoint action~\thetag{3.11}.  The vanishing of the cross terms is
equivalent to the triviality of the homomorphism $LZ_1\to\Hom(LG',\TT)$
defined by conjugation in the central extension.  But a homomorphism $\phi
\:LG'\to\TT$ is determined by its derivative $\dot\pi \:L\frak{g}'\to i\RR$,
which vanishes on commutators. Since commutators span~$L\frak{g}'$, we
deduce that $\phi $~is trivial.
        \enddemo

Fix a maximal torus~$T'\subset G'$, set $\Tt=Z_1\times T'$ a maximal torus
of~$\Gt$, let $T\subset G$ be its image under the projection $\Gt\to G$, and
denote the extended affine Weyl groups as~$\Waff', \Wafft, \Waff$.  The
fundamental groups of~$T',\Tt,T$ are~$\Pi ',\Pit,\Pi$.  The weight lattices
are $\Lambda =\Hom(T,\TT)$ and $\Lt = \Hom(Z_1,\TT)\times \Hom(T',\TT)$.
Orient~$T$ and~$Z_1$.  Note the exact sequence $0 @>>> \Pit @>>> \Pi @>>> A
@>>> 0$.  For readability we delete the map~$\kappa \:\Pi \to\Lambda $ and
its superscripted cousins from the notation.

        \proclaim{\protag{4.20} {Lemma}}
 There is a commutative diagram 
  $$ \xymatrix{\ZZ\bigl[\Lambda ^{\tms}/\Waff \bigr] \ar@{^{(}->}[r] &
     \ZZ\bigl[\Lt^{\tilde\tau -\tilde\sigma }/\Wafft \bigr]\\
      {R}^{\tms}(LG)\ar[u]^{\cong }\ar[d]_-{\Phi_G}\ar@{^{(}->}[r]&
     {R}^{\ttl-\stl}(L\Gt)\ar[u]_{\cong } \ar[d]^-{\Phi_{\Gt}}\\
      K^{\tau + \dim G}_{G}(G)\ar@{^{(}->}[r]\ar[d]_{\cong } & K^{\ttl + \dim
     G}_{\Gt}(\Gt) \ar[d]^{\cong } \\
      \Hom\mstrut _{\Waff}\bigl(\Ltil,\ZZ(\epsilon )\bigr) \ar@{^{(}->}[r] &
     \Hom\mstrut _{\Wafft}\bigl(\Lt^{\tilde\tau },\ZZ(\epsilon )\bigr)}
      \tag{4.21} $$
in which the horizontal arrows are injective.
        \endproclaim

        \demo{Proof}
 The top right vertical isomorphism follows from \theprotag{4.6}
{Proposition} and \theprotag{4.11} {Proposition}.  The lower isomorphisms are
\thetag{4.3}.  For the isomorphism in the upper left we factorize the second
line as
  $$ {R}^{\tms}(LG)\longrightarrow
     {R}^{\tms}\bigl(i(L\Gt)\bigr)\longrightarrow
     {R}^{\ttl-\stl}(L\Gt) . \tag{4.22} $$
The exact sequences~\thetag{4.18} and $ 0 @>>> \Lambda @>>>\Lt @>>> \Adual
@>>> 0 $ imply
  $$ {R}^{\tms}\bigl(i(L\Gt)\bigr) \cong
     \ZZ\bigl[\Lambda ^{\tms}/\Wafft \bigr], \tag{4.23} $$
since a representation of~$(L\Gt)^{\ttl-\stl}$ drops to the quotient if and
only if its restriction to~$A\subset L\Gt$ is trivial.  Now $A\cong
LG^{\tms}/i(L\Gt)^{\tms}$~acts on the set of isomorphism classes of
irreducible positive energy representations of ~$i(L\Gt)^{\tms}$, which under
the isomorphism in~\thetag{4.23} corresponds to the $A$-action on $\Lambda
^{\tms} /\Wafft$ induced from the $\Pi $-action on~$\Lambda ^{\tms} $.  We
claim that this action is free. To see this recall first $\Wafft=\PZ\times
\Waff'$.  The $\PZ$ action on~$\Lambda ^{\tms} /\Waff'$ is free, as
$\PZ$~acts by translations on $\Lt^{\ttl-\tilde\sigma }/\Waff'\cong
\Lambda_Z^{\tZ}\times (\Lambda ')^{\tau '-\sigma '}/\Waff'\supset \Lambda
^{\tms} /\Waff'$.  Since $\Pi /\Pi '\cong \pi _1G$~is torsion-free, and $\Pi
/\Pi '\supset\PZ$ with finite quotient~$A$, it follows that $\Pi /\Pi '$~acts
freely on~$\Lambda ^{\tms} /\Waff'$ as well.  Hence $A$~acts freely
on~$\Lambda ^{\tms} /\Wafft$.  Since the set of $A$-orbits in $\Lambda
^{\tms}/\Wafft$ is isomorphic to~$\Lambda ^{\tms}/\Waff$, the upper left
isomorphism follows.\footnote{That is, an irreducible representation
of~$LG^{\tms}$ corresponds to an orbit of irreducible representations
of~$i(L\Gt)^{\tms}$.  It may be that $A$~carries a nontrivial grading even if
the restriction of~$\tms$ to~$i(L\Gt)$ does not---see \theprotag{4.28}
{Example} below---but the conclusion is the same since the $A$-action is
free: the irreducible representations in an orbit, when summed to construct
an irreducible representation of~$LG^{\tms}$, may be even or odd.}
Furthermore, this argument shows that each of the arrows in~\thetag{4.22} is
injective, and now the vertical isomorphisms in~\thetag{4.21} prove that the
top two horizontal arrows in that diagram are injective.
 
The bottom half of~\thetag{4.21} factors as 
  $$ \xymatrix{
      K^{\tau + \dim G}_{G}(G)\ar[r]\ar[d]_{\cong } & K^{\tau + \dim
     G}_{\Gt}(G)\ar[r]\ar[d]_{\cong } & K^{\ttl + \dim G}_{\Gt}(\Gt)
     \ar[d]^{\cong } \\
      \Hom\mstrut _{\Waff}\bigl(\Ltil,\ZZ(\epsilon )\bigr) \ar[r] &
     \Hom\mstrut _{\Wafft}\bigl(\Ltil,\ZZ(\epsilon )\bigr) \ar[r] &
     \Hom\mstrut _{\Wafft}\bigl(\Lt^{\tilde\tau },\ZZ(\epsilon )\bigr)}
     \tag{4.24} $$
The bottom left map is the obvious inclusion and the bottom right map is
extension by zero, which is also injective.  Hence the bottom two arrows
in~\thetag{4.21} are injective.

It remains to show that the center square commutes.  For this we observe that
the space of connections on the trivial $G$~bundle over~$\cir$ is naturally
isomorphic to the space of connections on the trivial $\Gt$~bundle
over~$\cir$.  Furthermore, the map~$i$ in~\thetag{4.18} induces an
isomorphism of Lie algebras.  A positive energy representation~$V$
of~$LG^{\tms}$ pulls back to a positive energy representation~$\tilde{V}$
of~$L\Gt^{\ttl-\stl}$, and using these isomorphism we see that the associated
Dirac families~\thetag{3.32} coincide.
       \enddemo

        \demo{Proof of~\theprotag{3.43} {Theorem} \rom($G$~connected, $\pi
_1G$~torsion-free\rom)} 
 The arguments in~\S{4.1} and~\S{4.2} prove that~$\Phi _{\Gt}$
in~\thetag{4.21} is an isomorphism.  It then follows immediately from
\theprotag{4.20} {Lemma} that $\Phi _G$~is injective.  The groups at the top
left and bottom left of~\thetag{4.21} are free of equal rank (the number of
regular $\Waff$-orbits in~$\Ltil$), and from the text following~\thetag{4.17}
we see that a natural set of generators of~$R^{\tms}(LG)$ maps under~$\Phi
_{\Gt}$ to an independent set of indivisible elements in~$K_{\Gt}^{\ttl +
\dim G}(\Gt)$.  It now follows that $\Phi _G$~is surjective, hence an
isomorphism. 
	\enddemo

 \subhead \S{4.4}.  Examples
 \endsubhead

        \example{\protag{4.25} {Example}}
 ($G$~simple, simply connected)\quad These have a well-known classification:
$G=\SU_n$, $\Spin_n$, $ \operatorname{Sp}_n$, $\operatorname{G}_2$, $ \Fr_4$,
$ \Er_6$, $ \Er_7$, $ \Er_8$.  A central extension~$LG^\tau $ is determined
up to isomorphism by the inner product~$\langle \cdot ,\cdot \rangle_\tau $
on~$\frak{g}$ (see~\S{2.3}), or equivalently by the induced
map~$\kt\:\frak{t}\to \frak{t}^*$.  The one-dimensional real vector space of
invariant inner products is naturally identified with~$H^3(G;\RR)$ using a
multiple of the formula in~\thetag{1.6}:
  $$ \spreadlines{6pt}
      \aligned
      \langle \cdot ,\cdot \rangle&\longleftrightarrow \frac{1}{8\pi
     ^2}\;\Omega ,\\
      \Omega (\xi _1,\xi _2,\xi _3&)= \langle [\xi _1,\xi _2],\xi _3
     \rangle\endaligned. \tag{4.26} $$
The forms which arise from central extensions are {\it integral\/}, i.e.,
correspond to elements of~$H^3(G;\ZZ)$.  Let $\langle \cdot ,\cdot
\rangle$~denote the distinguished positive integral generator.  Then any
other integral invariant form is an integer---the {\it level\/}---times the
generator.  For example, the level of the spin extension~$LG^\sigma $ is the
{\it dual Coxeter number\/}~$\dC G$.  Suppose $LG^\tau $ is a positive
definite admissible central extension, and let the level of~$LG^{\tms}$
be~$k$.  Then the positive definiteness means~$k>-\dC G$.  If~$k<0$ then
there are no nonzero positive energy representations at level~$k$
(\theprotag{4.12} {Proposition}).  Thus assume~$k\ge0$.  Let $\alpha\in
\frak{t}^*$ be the highest root.  Then the lowest weights~$-\lambda $ in the
image of~$\kts(\calcove)$ satisfy $\lambda $~is dominant with~$\langle
\lambda ,\alpha\rangle\le k$.  The kernel of any Dirac family is supported on
a finite set of {\it regular\/} conjugacy classes in~$G$, namely those whose
$k+\dC G$ power contains the identity element.
        \endexample

        \example{\protag{4.27} {Example}}
 ($G=\SU_2$)\quad This is the simplest case of \theprotag{4.25} {Example}, and
we spell out the details a bit.  First, the dual Coxeter number
is~$\dC{\SU_2}=2$.  Let $T\subset \SU_2$ be the standard maximal torus of
diagonal matrices, and identify its Lie algebra~$\frak{t}$ with~$i\RR$ as
usual: $\left(\smallmatrix ia&0\\0&-ia \endsmallmatrix\right)\leftrightarrow
ia$.  Choose the alcove $\alcove=\bigl(ia:-1/2<a<0 \bigr)\subset \frak{t}$.
The affine Weyl group is $\Waff\cong Z\rtimes \zt$, where the generator
of~$\ZZ$ acts on~$\frak{t}$ as translation by~$i$ and the nontrivial element
of~$\zt$ acts as the reflection~$ia\mapsto -ia$.  The generating integral
inner product is\footnote{Let $\zeta =\left(\smallmatrix i&0\\0&-i
\endsmallmatrix\right)$ be the coroot.  Then the integral of the 3-form
$\Omega (\xi _1,\xi _2,\xi _3)= \langle [\xi _1,\xi _2],\xi _3 \rangle$ is
$\int_{\SU_2}\Omega = 4\pi ^2|\zeta |^2$.  The properly normalized generator
has~$|\zeta |^2=2$.  The normalization may be seen from the Chern-Weil
formula for the second Chern class, for example, using the Chern-Simons form
to transgress~$[\Omega ]\in H^3(\SU_2;\RR)$ to a class in~$H^4(\BSU_2;\RR)$.
It is one explanation for the factor~$8\pi ^2$ in~\thetag{4.26}.} $\langle
A,A' \rangle = -\Tr(AA')$, which restricts on~$\frak{t}$ to~$\langle ia,ia'
\rangle = 2aa'$.  We use~$\langle \cdot ,\cdot \rangle$ to
identify~$\frak{t}^*\cong i\RR$.  The fundamental group~$\Pi $ is identified
with~$i\ZZ\subset i\RR\cong \frak{t}$ and the weight lattice~$\Lambda $
with~$\frac i2\ZZ$.  A positive definite admissible central
extension~$\LSU_2^\tau $ has $\langle \cdot ,\cdot \rangle _\tau =(k+2)\langle
\cdot ,\cdot \rangle$ for some~$k>-1$, and then $(\LSU_2)^{\tms}$ is the
extension at level~$k$.  The map $\kts\:\frak{t}\to\frak{t}^*$ is identified
with multiplication by~$k$, so the image of~$\calcove$ is the closed interval
in~$i\ZZ$ from~$-ik/2$ to~$0$, which includes $k+1$~elements of~$\Lambda $.
So there are $k+1$~irreducible positive energy representations
of~$(\LSU_2)^{\tms}$.  The map~$\kt$ is identified with multiplication
by~$k+2$, and in the induced $\Waff$-action on~$\frak{t}^*$ the generator
of~$\ZZ$ translates by~$i(k+2)$.  There are $k+1$~{\it affine regular\/}
$\Waff$-orbits on~$\Lambda $; a slice is given by the set~$i\{\frac
12,1,\dots ,\frac{k+1}{1}\}$.  The corresponding conjugacy classes in~$\SU_2$
consist of matrices with eigenvalues~$\pm\tfrac{2\pi i\ell }{2(k+2)}$, where
$\ell =1,2,\dots ,k+1$.  The kernel of the Dirac family~\thetag{3.40} for an
irreducible positive energy representation has support on a single such
conjugacy class.
        \endexample

        \example{\protag{4.28} {Example}}
 ($G=\Ur_2$)\quad The group of components of the loop group~$\LU_2$ is
isomorphic to~$\ZZ$: a generator is the loop $ \varphi (z)=
\left(\smallmatrix z&0\\0&1 \endsmallmatrix\right)$.  We claim that the spin
extension~$\LU_2^\sigma $, defined in~\thetag{3.12}, has the nontrivial
grading, which is the nontrivial homomorphism $\LU_2\to \pi _0(\LU_2)\to\zt$.
Furthermore, the loop group~$L\Gt$ of the double covering
group~$\Gt=\SU_2\times \TT$ also has group of components isomorphic to~$\ZZ$,
and the induced map $\pi _0(L\Gt)\to\pi _0(LG)$ is multiplication by~2.  So
the pullback~$L\Gt^{\stl}$ central extension has trivial grading.  Notice
that for any central extension~$\LU_2^\tau $, either $\tau $~ or $\tau -\sigma
$~has a nontrivial grading.
 
To verify the claim we write an arbitrary loop in~$L\frak{u}_2$ as
  $$ \pmatrix ix&\zeta \\-\bar\zeta &iy \endpmatrix \;+\;
     \sum\limits_{n>0}\pmatrix \alpha _n&\beta _n\\\gamma _n&\delta _n
     \endpmatrix z^n \;-\; \sum\limits_{n>0}\pmatrix \bar\alpha _n&\bar\gamma
     _n\\\bar\beta _n&\bar\delta _n \endpmatrix z^{-n}.  \tag{4.29} $$
The operator~\thetag{3.8} (for the trivial connection~$A=A_0$) defines a
polarization~$\scrJ$, and in that polarization~\thetag{3.3} we choose the
particular complex structure~$J$ which maps~\thetag{4.29} to
  $$ \pmatrix iy&i\zeta \\i\bar\zeta &-ix \endpmatrix \;+\;
     \sum\limits_{n>0}\pmatrix i\alpha _n&i\beta _n\\i\gamma _n&i\delta _n
     \endpmatrix z^n \;+\; \sum\limits_{n>0}\pmatrix i\bar\alpha _n&i\bar\gamma
     _n\\i\bar\beta _n&i\bar\delta _n \endpmatrix z^{-n}.  $$
The only change to~$J$ under conjugation by~$\varphi $ is to its action
on~$\zeta $, which changes sign.  It follows that $\varphi J\varphi \inv
$~and $J$~are in opposite components of~$\scrJ$, and so under~\thetag{3.11}
the loop~$\varphi $ maps to the non-identity component of~$\OJ(H)$.

        \endexample

        \example{\protag{4.30} {Example}}
 ($G=\SO_3$)\quad Since $\pi _1\SO_3$~is not torsion-free, this case is not
covered by our work in this paper.  There are new phenomena.  The main point
is that the centralizer subgroup of any rotation of~$\RR^3$ through
angle~$\pi $ is not connected, which modifies \theprotag{4.14} {Corollary} as
well as the twisted equivariant $K$-theory.  A closely related fact is that
the argument after~\thetag{4.23} fails: the action of~$F\cong \zt$ has fixed
points and so the horizontal maps in~\thetag{4.21} are not injective.
Furthermore, the spin extension $\Ur_2\to \SO_3$ is not split; its equivalence
class in~$H^3(BSO_3;\ZZ)$ is the universal integral third Stiefel-Whitney
class.  This example is covered by~\cite{FHT3}, which treats arbitrary
(twisted) loop groups and central extensions.  Detailed computations for
$G=\SO_3$ appear in~\cite{FHT2,Appendix~A}.
        \endexample

 \head
 Appendix: Central extensions in the semisimple case
 \endhead
 \comment
 lasteqno A@ 18
 \endcomment

Here we present the proof of~\theprotag{2.15} {Proposition}, which we repeat
for convenience.

        \proclaim{\protag{2.15} {Proposition}}
 Let $G$~be a compact Lie group with $[\frak{g},\frak{g}]=\frak{g}$, i.e.,
$\frak{g}$~semisimple.  Let $P\to\cir$ be a principal $G$-bundle.  Then any
central extension~$\LGRt$ is admissible.  Furthermore, for each~$\LGRRt$
which satisfies~\thetag{2.11} there exists a unique $\LGRRt$-invariant
symmetric bilinear form $\form\;=\;\form_\tau $ which
satisfies~\thetag{2.12}.
        \endproclaim

\flushpar 
 For $G$~connected and simply connected, see~\S{2.3}.  We remark that in the
main text we only used the statement about the bilinear form (in the proof of
\theprotag{3.13} {Proposition}).

        \demo{Proof}
 Central extensions of~$\LGR$ by~$\TT$ are classified by the smooth group
cohomology~$H^2_{\smooth}(\LGR;\TT)$, as defined in~\cite{S}.  For
$G$~semisimple $H^q_{\smooth}(\LGR;\RR)=0$ for~$q>0$ by~\cite{PS,\S14.6},
\cite{RW}, so the exponential sequence yields
  $$ H^2_{\smooth}(\LGR;\TT)\cong H^3_{\smooth}(\LGR;\ZZ)\cong
     H^3(B\LGR;\ZZ).  $$
In other words, the isomorphism class of a central extension of~$\LGR$ is
characterized topologically.  We must show that
  $$ H^3(B\LGRR;\ZZ)\longrightarrow H^3(B\LGR;\ZZ)\quad \text{is surjective}
     \tag{A.1} $$
for some finite cover~$\Troth\to\Trot$, and so construct~ $\LGRRt$ as
in~\thetag{2.11}.  For this consider the Leray spectral sequence
$H^p\bigl(\BTroth;H^q(B\LGR;\ZZ) \bigr)\Longrightarrow H^{p+q}(B\LGRR;\ZZ)$.
From the discussion surrounding~\thetag{2.3} we see $H^q(B\LGR;\ZZ)\cong
H^q_G\bigl(G[P];\ZZ \bigr)$.  The first differential of interest is
  $$ d_2\:H^3_G\bigl(G[P];\ZZ \bigr)\cong H^0\bigl(\BTroth;H^3_G(G[P];\ZZ)
     \bigr)\longrightarrow H^2\bigl(\BTroth;H^2_G(G[P];\ZZ)
     \bigr). \tag{A.2} $$
The hypothesis that $G$~has no torus factors implies $H^2_G\bigl(G[P];\ZZ
\bigr)$ ~is a finite group---use the Leray spectral sequence for the homotopy
quotient of~$G[P]$ fibered over~$BG$---and so for an appropriate cover the
codomain of~\thetag{A.2} vanishes.  The next nonzero differential is
  $$ d_4\: H^0\bigl(\BTroth;H^3_G(G[P];\ZZ) \bigr)\longrightarrow
     H^4\bigl(\BTroth;H^0_G(G[P];\ZZ) \bigr). \tag{A.3} $$
For a suitable cover the bundle $\LGRR\to\Troth$ has a section, in which case
\thetag{A.3}~vanishes.  This completes the proof of~\thetag{A.1}.

Now we assume given a group~$\LGRRt$ which fits into~\thetag{2.11}.  Fix a
connection~$A$ which is generic in the sense that its stabilizer~$Z_A\subset
\LGR$ has identity component a torus.  We use the basepoint to identify it
with a subgroup of~$G$, and so the complexification~$(\zA)_{\CC}$ of its Lie
algebra with an abelian subalgebra~$\frak{h}_0\subset \frak{g}_\CC$.  Let
$\frak{h}\subset \frak{g}_\CC$ be the centralizer of~$\frak{h}_0$; then
$\frak{h}$~is a Cartan subalgebra.\footnote{In other words, the real points
of~$\frak{h}$ form the Lie algebra of a maximal torus.  To see this, we first
show $\frak{h}_0\not= 0$.  Identify the holonomy automorphism
of~$\frak{g}_\CC$ with~$\Ad_g$ for some~$g\in G$.  Composing with~$\Ad_{g_1}$
for a suitable~$g_1$ in the identity component, we may assume that the
automorphism of~$\frak{g}_\CC$ fixes the Lie algebra of a maximal torus as
well as a Weyl chamber.  Then it permutes the positive roots, so fixes their
sum and the line generated by the sum of the coroot vectors.  Conjugating
back by~$\Ad_{g_1\inv }$ we conclude~$\frak{h}_0\not= 0$.  Now $g^n$~lies in
the identity component~$G_1$ for suitable~$n$, so may be written as the
exponential of a real element of~$\frak{h}_0$.  Multiply~$g$ by its inverse
to obtain an automorphism of~$\frak{g}_\CC$ of finite order.  Then
\cite{Ka,Lemma~8.1}~applies to prove that $\frak{h}$~is a Cartan subalgebra.
Note for the untwisted case ($P\to\cir$ trivializable) we
have~$\frak{h}=\frak{h}_0$.}  Decompose
  $$ \frak{g}_{\CC}\cong \frak{h}\;\oplus\;\bigoplus_{\lambda \in \Delta
     _0}\frak{g}_\lambda ,\qquad \lambda \in \frak{h}_0^*,  $$
into eigenspaces of~$\ad(\frak{h}_0)$, where $\frak{g}_\lambda \not= 0$.
Then $\lambda \not= 0$ for~$\lambda \in \Delta _0$, since
$\frak{h}_0$~contains regular elements.  The holonomy induces an automorphism
of~$\frak{h}$ which fixes~$\frak{h}_0$, and since $\frak{h}_0$~contains
regular elements we can choose a Weyl chamber which is invariant under the
holonomy.  Now $\ad(\frak{h})$ decomposes each~$\frak{g}_\lambda $ as a sum
of root spaces, and by our choice of Weyl chamber the roots which occur are
either all positive or all negative.  In this way we partition~$\Delta _0$
into a positive set and a negative set.  Use~\thetag{2.38} to
embed~$\frak{g}_\CC$ in the complexified loop algebra.  (We drop the
subscript~`$A$'.)  Note~$[d_A,\frak{h}_0]=0$, but $[d_A,\frak{h}]\not= 0$ in
general, and so $\frak{h}$~decomposes as a sum of eigenspaces of~$d_A$.  We
write $\frak{h}=\frak{h}_0\oplus \frak{h}'$, where $d_A$~has nonzero
eigenvalues on~$\frak{h}'$.  By adding to~$d_A$ a suitable regular element
of~$\frak{h}_0$ if necessary, and making an appropriate choice of the
logarithm of the holonomy (`$\SA$'~in the text preceding~\thetag{2.36}), we
can ensure
  $$ \gather
      \spec(-i \,d_A)\subset (0,1)\qquad \text{on $\frak{h}' \;\oplus
     \;\bigoplus\limits_{\lambda>0 }\frak{g}_\lambda $,} \tag{A.4} \\
       \text{$\spec(-i\,d_A)$ on $\frak{h}$ and $\bigoplus\limits_{\lambda
     \in \Delta _0}\frak{g}_\lambda $ are disjoint}. \tag{A.5} \endgather $$
Summarizing, we have a decomposition\footnote{We implicitly assume the joint
eigenspaces of~$\ad(\frak{h}_0)$ and~$d_A$ with nonzero eigenvalues have
dimension one, which may well be true in general.  If not, the argument is
only notationally more complicated.}
  $$ \frak{g}_\CC\cong \frak{h}_0 \;\oplus \; \bigoplus\limits_{j}\CC\cdot
     \cj\;\oplus \;\bigoplus\limits_{\lambda ,e>0}\CC\cdot \pa\;\oplus
     \;\bigoplus\limits_{\lambda ,e>0}\CC\cdot \pab,  $$
where the $\cj,\pa$~are determined up to a nonzero scalar by
  $$ \align
      [\xi ,\pa] &= i\lambda (\xi )\pa,\qquad \xi \in \frak{h}_0 
     \\
      [d_A ,\pa] &= ie\pa,  \\
      [d_A ,\cj] &= iE(\cj)\cj,  \endalign $$
and $0<e<1$ if~$\lambda >0$.  We arrange~$\pab = \overline{\pi _{\lambda
,e}}$.

Observe that if $\eta \in \LgR_\CC$ is an eigenvector of~$\ad(d_A)$ with
nonzero eigenvalue~$iE(\eta )$, then it has a lift~$\eta ^\tau _A$
to~$\LgRt_{\CC}$ characterized by
  $$ [d_A^\tau ,\eta ^\tau _A] = iE(\eta )\eta ^\tau _A. \tag{A.6} $$
Also, equation~\thetag{A.6} implies that if $\eta ,\eta
'$ are elements of~$\LgR_{\CC}$ of definite nonzero energy, and the sum of
the energies is nonzero, then
  $$ [\eta ,\eta ']^\tau _A = [\eta ^\tau _A,\lift{\eta '}],\qquad E(\eta
     ),E(\eta '),E([\eta ,\eta '])\not= 0. \tag{A.7} $$
Thus we have almost defined a decomposition of the finite energy vectors in
the central extension:
  $$ \bigl(\LgRRt_{\CC}\bigr)_{\operatorname{fin}}(A) \cong \CC\cdot d_A^\tau 
     \;\oplus\; \CC\cdot K \;\oplus\; \frak{h}_0\;\oplus
     \;\bigoplus\limits_j\CC\cdot \lift{\cj} \;\oplus\;
     \bigoplus\limits_{\lambda ,e \not= 0}\CC\cdot \lift{\pa} \;\oplus\;
     \bigoplus\limits_{n\not= 0}z^n\frak{g}_{\CC}. \tag{A.8} $$
It remains to define the lift of~$\frak{h}_0$ to the central extension, which
we do below.
 
The invariance of the desired form~$\form_\tau $ implies that it is nonzero
only if the total energy of its arguments vanishes.  (We implicitly use
\theprotag{2.41} {Lemma} here as well.)  Similarly, eigenvectors
for~$\frak{h}_0$ pair nontrivially only if their eigenvalues (as linear
functionals on~$\frak{h}_0$) sum to zero.  If $\eta ,\eta '\in \LgR_{\CC}$,
then by semisimplicity write $\eta =\sum [\xi _i,\xi '_i]$ for some~$\xi
_i,\xi '_i\in \LgR_{\CC}$.  Any invariant form satisfies
  $$ \split
      \ip{z^n\eta }{z^{-n}\eta '}_\tau \; &= \sum\limits_{i}\ip{[\xi
     _i,z^n\xi '_i]}{z^{-n}\eta '}_\tau  \\
      &= \sum\limits_{i}\ip{\xi _i}{[z^n\xi '_i,z^{-n}\eta ']}_\tau  \\
      &= \sum\limits_{i}\ip{\xi _i}{[\xi '_i,\eta ']}_\tau  \\
      &= \sum\limits_{i}\ip{[\xi _i,\xi '_i]}{\eta '}_\tau  \\
      &= \;\ip\eta {\eta '}_\tau .\endsplit \tag{A.9} $$
Therefore, it suffices to define~$\form_\tau $ on~$\frak{g}_\CC$; all other
pairings are determined by~\thetag{2.12}, \thetag{2.13}, \thetag{2.14},
\thetag{2.19}, and~\thetag{A.9}. 
 
Set 
  $$ \za=i[\pa,\pab]\in \frak{h}_0.  $$
The~$\za$ span~$\frak{h}_0$, but in general there are linear relations among
them.  Define a sequence of lifts of~$\za$ by
  $$ \zan n= i\bigl[ \lift{z^n\pa},\lift{z^{-n}\pab}\bigr]\,\in \LgRt,\qquad
     n\in \ZZ. \tag{A.10} $$
We verify below that $\zan {n+1} - \zan{n}\in \RR\cdot K$ is independent
of~$n$, and we use it to define \break$\ip{\pa}{\pab}_\tau $:
  $$ \delta (\za) = \zan {n+1} - \zan{n} =\;\ip{\pa}{\pab}_\tau \,K,\qquad
     n\in \ZZ. \tag{A.11} $$
This definition is forced by invariance: 
  $$ \split
      \ip{d_A^\tau }{\delta (\za)}_\tau \;&= i\;\ip{d_A^\tau }{\bigl[
     \lift{z^{n+1}\pa},\lift{z^{-(n+1)}\pab}\bigr]}_\tau - \,i\ip{d_A^\tau }{\bigl[
     \lift{z^{n}\pa},\lift{z^{-n}\pab}\bigr]}_\tau \\
      &= -(n+1+e)\,\ip{z^{n+1}\pa}{z^{-(n+1)}\pab}_\tau +\;
     (n+e)\,\ip{z^{n}\pa}{z^{-n}\pab}_\tau \\
      &=-\;\ip{\pa}{\pab}_\tau ,\endsplit  $$
where at the last stage we use~\thetag{A.9}.  We also define the lift 
  $$ \lift{\za} = \zan 0 - e \delta (\za). \tag{A.12} $$
This is also determined by invariance (and~\thetag{2.19}): if $\lift{\za} =
\zan0 + cK$ for some constant~$c$, then 
  $$ 0 =\;\ip{d_A^\tau }{\lift{\za}}_\tau  \;=i\ip{d_A^\tau }{\bigl[
     \lift{\pa},\lift{\pab}\bigr]}_\tau  - \;c \;= -e \ip{\pa}{\pab}_\tau  -\;
     c.  $$
For any~$\xi \in \frak{h}_0$ we set
  $$ \ip{\za}\xi_\tau  \;= - \lambda (\xi )\ip{\pa}{\pab}_\tau , \tag{A.13} $$
as determined by invariance: 
  $$ \ip{\za}\xi_\tau \;=i\ip{[\pa,\pab]}{\xi }_\tau \;=i\ip{\pa}{[\pab,\xi
     ]}_\tau \;=-\lambda (\xi )\ip{\pa}{\pab}_\tau .  $$
Finally, by semisimplicity we can write (non-uniquely) 
  $$ \cj = \sum\limits_{k} [\ejk,\epjk], \tag{A.14} $$
where each $\ejk,\epjk$ is a multiple of some~$z^n\pa$, and
$E(\ejk)+E(\epjk)=E(\cj)$.  Then define 
  $$ \ip{\cj}{z\inv \chi \mstrut _{j'}}_\tau \;=\sum\limits_{k,\ell
     }\ip{[[\ejk,\epjk],\eta _{j'}^{(\ell )}]}{z\inv {\eta '_{j'}}^{(\ell
     ')}}_\tau , \tag{A.15} $$
which is forced by invariance.  Condition \thetag{A.5}~implies that the
triple bracket is a multiple of some~$z^n\pa$, and so the right hand side
of~\thetag{A.15} is determined by~\thetag{A.11} and~\thetag{A.9}.

The invariance arguments prove that $\form_\tau $~is unique.  We must check,
though, that \thetag{A.12}, \thetag{A.13}, and~\thetag{A.15} are consistent.
In particular, this will show that \thetag{A.12}~completes the definition of
the splitting~\thetag{A.8} (and so the form: $\ip{d_A^\tau }{\zeta ^\tau
_A}\;=0_\tau $ for all~$\zeta \in \frak{h}_0$).
 
First, we verify that \thetag{A.11}~is independent of~$n$.  Choose~$\xi \in
\frak{h}_0$ regular, i.e., with~$\lambda (\xi )\not= 0$ for all~$\lambda $.
Then using~\thetag{A.7} repeatedly we find
  $$ \split
      \zan{n+1} &= \frac{1}{\lambda (\xi )}\bigl[ [z\xi ,z^n\pa]
     _A^\tau ,\lift{z^{-(n+1)}\pab}\bigr] \\
      &= \frac{1}{\lambda (\xi )}\bigl[ [\lift{z\xi} ,\lift{z^{-(n+1)}\pab}]
     ,\lift{z^n\pa}\bigr] + \frac{1}{\lambda (\xi )}\bigl[ \lift{z\xi}
     ,[\lift{z^n\pa},\lift{z^{-(n+1)}\pab}]\bigr] \\
      &= \zan n - \frac{i}{\lambda (\xi )}[\lift{z\xi },\lift{z\inv
     \za}],\endsplit  $$
so 
  $$ \delta (\za) = \zan {n+1} - \zan{n} = -\,\frac{i}{\lambda (\xi
     )}[\lift{z\xi },\lift{z\inv \za}] \tag{A.16} $$
is independent of~$n$, as claimed.   
 
Next, suppose $c^{\lambda,e} \za=0$ for some~$c^{\lambda,e} \in \RR$.  Then
we find using~\thetag{A.11} and~\thetag{A.16} that
  $$ \split
      \sum\limits_{\lambda >0} c^{\lambda,e} \lambda (\xi
     )\ip{\pa}{\pab}_\tau K
      &=-\, \sum\limits_{\lambda >0}ic^{\lambda,e} \bigl[ \lift{z\xi
     },\lift{z\inv \za}\bigr]\\
      &=-\, \sum\limits_{\lambda >0 }i \bigl[ \lift{z\xi },\lift{z\inv
     c^{\lambda,e} \za}\bigr]\\
      &= 0. \endsplit $$
This demonstrates that \thetag{A.13}~is a consistent definition
of~$\form_\tau $ on~$\frak{h}_0$.
 
To check that \thetag{A.15}~is consistent, we must check that any
expression~\thetag{A.14} leads to the same right hand side of~\thetag{A.15},
or equivalently that if $\sum\limits_{k}[\eta ^{(k)},{\eta '}^{(k)}]=0$ then
the right hand side of~\thetag{A.15} vanishes.  But this is obvious. 
 
To see that the lift of~\thetag{A.12} is consistent we must show that any
linear relation among the~$\za$ is also satisfied by the lifts.  Let $\lambda
_1,\dots ,\lambda _s\in \frak{h}_0^*$ be the restrictions of the simple roots
to~$\frak{h}_0$, after eliminating duplicates; call these the {\it
simple\/}~$\lambda $.  Then $\{\zeta _{\lambda _1,e_1},\dots ,\zeta _{\lambda
_s,e_s}\}$ is a basis of~$\frak{h}_0$, where we choose $0<e_1,\dots ,e_s<1$
to be minimal.  There are two types of relation which generate them all.
First, if $\lambda _k$~is simple, then any~$\pi _{\lambda _k,e}$ which occurs
is of the form~$\pi _{\lambda _k,e_k+E(\cj)}$ for some~$\cj\in \frak{h}'$ of
definite (positive) energy, and for a suitable choice we can take $\pi\mstrut
_{\lambda _k,e} = [\cj,\plek]$.  Then for some constant~$c$ we have
$[\pi\mstrut _{\lambda _k,e},\bar\chi \mstrut _j]=c\plek$, and for any~$n\in
\ZZ$ we find
  $$ \split
      \zeta \mstrut _{\lambda _k,e} &= i[z^n \pi \mstrut _{\lambda
     _k,e},z^{-n}[\bar\chi \mstrut _j,\plekb]] \\
      &= i[[ z^{n}\pi \mstrut _{\lambda _k,e},\bar\chi\mstrut
     _j],z^{-n}\plekb]]] \;+\;
       i[\bar\chi \mstrut _j,[z^n\pi \mstrut _{\lambda _k,e},z^{-n}\plekb]]
     \\
      &= c\,\zet{k}.\endsplit \tag{A.17} $$
Now in the central extension an easy argument using invariance and
\thetag{A.15} shows that  
  $$ i[\lift{\bar\chi \mstrut _j},[\lift{z^n\pi \mstrut _{\lambda
     _k,e}},\lift{z^{-n}\plekb}]] = cE(\cj)\delta (\zet k).  $$
By repeatedly applying~\thetag{A.7}, and using the definition~\thetag{A.12}
of the lifts, we verify the relation~\thetag{A.17} for the lifts.  The second
type of relation comes from writing any~$\lambda >0$ as a sum of
simple~$\lambda _k$, and we proceed by induction on the length of such a
relation.  In the inductive step we take
  $$ (\lambda ,e) = (\lambda ',e') + (\lambda _k,e_k) - (0,\nu ),
      $$
where $0<e,e',e_k<1$, which determines $\nu = 0$ or~$\nu =1$, and we already
know the lift of~$\zlep$ is in the linear span of the lifts of~$\zet1,\dots
,\zet s$.  Take 
  $$ \pa = z^{-\nu }[\plep,\plek]. \tag{A.18} $$
Then $[\plep,\plekb]$ is a (possibly vanishing) multiple of~$\pi \mstrut
_{\lambda '-\lambda _k,e'-e_k}$.  Substituting~\thetag{A.18} and expanding we
find a relation
  $$ \split
      \za &= i[z^\nu \pa, \zeta ^{-\nu }\pab] \\
      &= \lambda '(\zet k)\zlep + \lambda _k(\zlep)\zet k \\
      &\qquad \qquad ic\left[ \ad(\plep)\ad(\plekb)\pi \mstrut _{-(\lambda
     '-\lambda _k),-(e'-e_k)} - \ad(\plek)\ad(\plepb) \pi \mstrut _{\lambda
     '-\lambda _k,e'-e_k}\right] \\
      &= \left[ \lambda '(\zet k)+c_1 \right]\zlep \;+\;\left[ \lambda
     _k(\zlep) + c_2 \right]\zet k \endsplit  $$
for some (possibly vanishing) constants~$c_1,c_2$.  Using~\thetag{A.10}
and~\thetag{A.12} we see that the lifts satisfy the same relation if and only
if
  $$ (e+\nu )\delta (\za) = \left[ \lambda '(\zet k)+c_1 \right]e'\delta
     (\zlep) \;+\;\left[ \lambda _k(\zlep) + c_2 \right]e_k\delta (\zet
     k).  $$
But this follows from~\thetag{A.16} after choosing~$\xi \in \frak{h}_0$ such
that $e_k=\lambda _k(\xi )$ and~$e'=\lambda '(\xi )$. 

This completes the proof that $\form_\tau $~is well-defined and invariant
under the adjoint action of~$\LgRRt$, so under the adjoint action of the
identity component of~$\LGRRt$.  Finally, the uniqueness shows that it is
invariant under the entire group $\LGRRt$.
        \enddemo

\widestnumber\key{SSSSSSSSS}   

\enddocument